%% file: 2025_ApelZilk_CMAME_preprint.tex
\documentclass[preprint,3p]{elsarticle}

\usepackage{amssymb}

\usepackage[english]{babel}
\usepackage[utf8]{inputenc}
\usepackage[T3,OT1]{fontenc}
\usepackage{lmodern}
\usepackage{latexsym}
\usepackage{array}
\usepackage{xcolor}
\usepackage{pdfcolmk} 
\usepackage{amsmath}
\usepackage{amssymb}
\usepackage{mathtools}
\usepackage{dsfont}
\usepackage{upgreek}
\usepackage{amsthm}    
\usepackage{amsfonts}  
\usepackage{relsize}   
\usepackage{graphicx}
\usepackage{enumitem}
\usepackage{esvect}
\usepackage{url}
\usepackage{float}                        
\usepackage{siunitx} 
\usepackage{colortbl} 
\usepackage{color,soul}   
\usepackage{todonotes}
\usepackage{subcaption}
\usepackage{tikz-cd}

\DeclareUnicodeCharacter{0306}{\v}

\usepackage{tikz, pgfplots}
\tikzstyle{bag} = [align=center]
\usetikzlibrary{shapes.misc}
\tikzset{cross/.style={cross out, draw=black, minimum size=2*(#1-\pgflinewidth), inner sep=0pt, outer sep=0pt},
	cross/.default={1pt}}
\usetikzlibrary{arrows}
\usetikzlibrary{patterns}
\pgfplotsset{compat = newest}

\usepackage[
colorlinks=true,
urlcolor=blue,
linkcolor=red,
citecolor=teal
]{hyperref}

\newcommand{\N}{{\mathbb N}}

\newcommand{\R}{{\mathbb R}}

\newcommand{\uproman}[1]{\uppercase\expandafter{\romannumeral#1}}

\def\.{\cdot}

\newcommand*{\pol}{\mathrel{\mathsmaller{\mathsmaller{\mathsmaller{\bs \odot}}}}}
\DeclareSymbolFont{tipa}{T3}{cmr}{m}{n}
\DeclareMathAccent{\invbreve}{\mathalpha}{tipa}{16}

\newcommand{\dt}{\,\mathrm{d}t}
\newcommand{\dr}{\,\mathrm{d}r}
\newcommand{\dphi}{\,\mathrm{d}\varphi}

\newcommand{\bs}{\boldsymbol}
\newcommand{\mr}{\mathring}

\newcommand{\abs}[1]{\left\lvert#1\right\rvert}
\newcommand{\norm}[1]{\left\lVert#1\right\rVert}

\newcommand{\Dd}{\,\textnormal{d}}

\newtheorem{Theorem}{Theorem}[section]
\newtheorem{lemma}[Theorem]{Lemma}
\newtheorem{corollary}[Theorem]{Corollary}
\theoremstyle{definition}
\newtheorem{definition}[Theorem]{Definition}
\theoremstyle{definition}
\newtheorem{example}[Theorem]{Example}

\newtheoremstyle{remarkstyle} 
{}{}{}{}{\bfseries}{.}{.5em}{{\thmname{#1 }}{\thmnumber{#2}}{\thmnote{ (#3)}}}
\theoremstyle{remarkstyle}
\newtheorem{remark}[Theorem]{Remark}

\numberwithin{equation}{section}

\graphicspath{{Bilder/}}

\newcommand{\logLogSlopeTriangleUpsideDownNeg}[6]
{
	
	\pgfplotsextra
	{
		\pgfkeysgetvalue{/pgfplots/xmin}{\xmin}
		\pgfkeysgetvalue{/pgfplots/xmax}{\xmax}
		\pgfkeysgetvalue{/pgfplots/ymin}{\ymin}
		\pgfkeysgetvalue{/pgfplots/ymax}{\ymax}
		
		\pgfmathsetmacro{\xArel}{#1}
		\pgfmathsetmacro{\yArel}{#3}
		\pgfmathsetmacro{\xBrel}{#1-#2}
		\pgfmathsetmacro{\yBrel}{\yArel}
		\pgfmathsetmacro{\xCrel}{\xArel}
		
		\pgfmathsetmacro{\lnxB}{\xmin*(1-(#1-#2))+\xmax*(#1-#2)} 
		\pgfmathsetmacro{\lnxA}{\xmin*(1-#1)+\xmax*#1} 
		\pgfmathsetmacro{\lnyA}{\ymin*(1-#3)+\ymax*#3} 
		\pgfmathsetmacro{\lnyC}{\lnyA-#4*(\lnxA-\lnxB)}
		\pgfmathsetmacro{\yCrel}{\lnyC-\ymin)/(\ymax-\ymin)} 
		
		\coordinate (A) at (rel axis cs:\xArel,\yArel);
		\coordinate (B) at (rel axis cs:\xBrel,\yBrel);
		\coordinate (C) at (rel axis cs:\xCrel,\yCrel);
		
		\draw[#5]   (A)-- node[pos=0.5,anchor=north] {1}
		(B)--
		(C)-- node[pos=0.5,anchor=west] {#6}
		cycle;
	}
}

\newcommand{\logLogSlopeTriangleNeg}[6]
{
	
	\pgfplotsextra
	{
		\pgfkeysgetvalue{/pgfplots/xmin}{\xmin}
		\pgfkeysgetvalue{/pgfplots/xmax}{\xmax}
		\pgfkeysgetvalue{/pgfplots/ymin}{\ymin}
		\pgfkeysgetvalue{/pgfplots/ymax}{\ymax}
		
		\pgfmathsetmacro{\xArel}{#1}
		\pgfmathsetmacro{\yArel}{#3}
		\pgfmathsetmacro{\xBrel}{#1+#2}
		\pgfmathsetmacro{\yBrel}{\yArel}
		\pgfmathsetmacro{\xCrel}{\xArel}
		
		\pgfmathsetmacro{\lnxB}{\xmin*(1-(#1-#2))+\xmax*(#1-#2)} 
		\pgfmathsetmacro{\lnxA}{\xmin*(1-#1)+\xmax*#1} 
		\pgfmathsetmacro{\lnyA}{\ymin*(1-#3)+\ymax*#3} 
		\pgfmathsetmacro{\lnyC}{\lnyA+#4*(\lnxA-\lnxB)}
		\pgfmathsetmacro{\yCrel}{\lnyC-\ymin)/(\ymax-\ymin)} 
		
		\coordinate (A) at (rel axis cs:\xArel,\yArel);
		\coordinate (B) at (rel axis cs:\xBrel,\yBrel);
		\coordinate (C) at (rel axis cs:\xCrel,\yCrel);
		
		\draw[#5]   (A)-- node[pos=0.5,anchor=south] {1}
		(B)-- 
		(C)-- node[pos=0.5,anchor=east] {#6}
		cycle;
	}
}

\numberwithin{equation}{section}

\journal{}

\begin{document}

\begin{frontmatter}



\title{Error Estimates and Graded Mesh Refinement for Isogeometric Analysis on Polar Domains with Corners}

\author[unibw]{Thomas Apel}
\ead{thomas.apel@unibw.de}
 
\author[unibw]{Philipp Zilk\corref{cor1}}
\ead{philipp.zilk@unibw.de}
 
\affiliation[unibw]{organization={Institute for Mathematics and Computer-Based Simulation, Universität der Bundeswehr München},
             addressline={\\Werner-Heisenberg-Weg 39},
             city={Neubiberg},
             postcode={85577},
             country={Germany}}

\cortext[cor1]{Corresponding author}

%

\begin{abstract}
Isogeometric analysis (IGA) enables exact representations of computational geometries and higher-order approximation of PDEs. In non-smooth domains, however, singularities near corners limit the effectiveness of IGA, since standard methods typically fail to achieve optimal convergence rates. These constraints can be addressed through local mesh refinement, but existing approaches require breaking the tensor-product structure of splines, which leads to increased implementation complexity.

This work introduces a novel local refinement strategy based on a polar parameterization, in which one edge of the parametric square is collapsed into the corner. By grading the standard mesh toward the collapsing edge, the desired locality near the singularity is obtained while maintaining the tensor-product structure. A mathematical analysis and numerical tests show that the resulting isogeometric approximation achieves optimal convergence rates with suitable grading parameters.

Polar parameterizations, however, suffer from a lack of regularity at the polar point, making existing standard isogeometric approximation theory inapplicable. To address this, a new framework is developed for deriving error estimates on polar domains with corners. This involves the construction of polar function spaces on the parametric domain and a modified projection operator onto the space of $C^0$-smooth polar splines. The theoretical results are verified by numerical experiments confirming both the accuracy and efficiency of the proposed approach.
%
\end{abstract}



\begin{keyword}
Isogeometric analysis \sep Error estimates \sep Corner singularities \sep Local mesh refinement \sep Graded meshes \sep Polar parameterizations



\MSC[2010]  65N30 \sep 65N25 \sep 65D07 \sep 65N50 \sep 35A21

\end{keyword}

\end{frontmatter}




\section{Introduction}
Isogeometric analysis was introduced in \cite{HughesCottrellBazilevs2005} as an advanced numerical method for solving PDEs, offering precise representations of CAD-type geometries and efficient higher-order approximation of functions. However, on non-smooth domains with corners, standard Galerkin methods fail to achieve higher-order convergence due to the reduced regularity of the PDE solution. To address this issue, several strategies have been proposed in the literature, with local mesh refinement being among the most widely used.

Corner singularities are well understood and can be handled effectively by using a priori mesh grading schemes. Yet, graded refinement of tensor product meshes on smoothly parameterized domains results in non-local refinement along the mapped edges of the parametric domain, see, for instance, previous works on mesh grading in IGA  \cite{BeiraodaVeigaChoSangalli2012,LangerMantzaflarisMooreToulopoulos2015}. To achieve true locality, various alternative spline technologies have been introduced, including T-splines, hierarchical splines, LR-splines, and PHT-splines. This list is not exhaustive, and a more comprehensive overview is given in \cite{BuffaGantnerGiannelliPraetoriusVazquez2022}. These spline constructions, however, require breaking the tensor-product structure and are therefore more complex to implement.

We propose a new strategy to construct locally refined tensor-product spline meshes near singularities by combining mesh grading with polar parameterizations. Classical two-dimensional model domains with corners -- such as circular sectors or L-shaped domains -- can be represented conveniently in IGA by collapsing one edge of the parametric square into a conical point. Since such parameterizations resemble a transformation from polar to Cartesian coordinates, the resulting geometries are referred to as polar domains with corners. By grading the parametric mesh toward the collapsed edge, true local refinement of the B\'ezier mesh near the singularity is achieved while preserving the tensor-product structure and without introducing additional degrees of freedom. In this way, optimal approximation of singular solutions can be obtained for suitable grading parameters, as demonstrated numerically for circular sectors in our recent work \cite{ApelZilk2024}. In the present paper, we extend our approach to more general domains with corners. More importantly, we provide a mathematical analysis of the method, which includes error estimates and a theoretical explanation of the required grading parameters.

Polar parameterizations play a significant role in various applications of isogeometric analysis and have been the focus of several recent studies \cite{ToshniwalSpeleersHiemstraHughes2017,SpeleersToshniwal2021}. For instance, the well-known Pacman domain is frequently used as a benchmark in the development and evaluation of adaptive refinement algorithms \cite{FaliniGiannelliKanducSampoliSestini2019,FeischlGantnerHaberlPraetorius2016}. Polar domains also arise in practical applications, such as the modeling of cylindrical sector geometries in electric rotors \cite{WiesheuKomannMerkelSchoepsUlbrichCortes2024}. In addition, they appear in the context of scaled boundary parameterizations in IGA \cite{ArioliShamanskiyKlinkelSimeon2019}, which have also been extended to the simulation of plates and shells \cite{ArfReichleKlinkelSimeon2023,ReichleArfSimeonKlinkel2023,ReichleKlassenLiKlinkel2024}.

However, polar mappings suffer from a lack of regularity, and standard isogeometric approximation spaces -- defined as push-forwards of the corresponding parametric spaces -- are not suitable \cite{TakacsJuettler2011,TakacsJuettler2012,ToshniwalSpeleersHiemstraHughes2017}. For this reason, these cases are excluded from standard isogeometric approximation theory, which instead assumes that the parameterization is sufficiently regular  \cite{BazilevsBeiraoDaVeigaCottrellHughesSangalli2006,BeiraodaVeigaChoSangalli2012,BeiraodaVeigaBuffaSangalliVazquez2014}. A complete theoretical understanding remains lacking when this assumption is not satisfied. To our knowledge, approximation properties have only been established for a specific type of singularly parameterized triangles \cite{Takacs2015preprint} and, more recently, for subdivision based isogeometric analysis \cite{Takacs2025}. In this paper, we address this gap in the literature by introducing a framework for deriving error estimates for isogeometric analysis on polar domains with corners. To this end, we construct a projector onto the well-known space of $C^0$-smooth polar splines \cite{TakacsJuettler2011,TakacsJuettler2012,ToshniwalSpeleersHiemstraHughes2017}, which serves as the approximation space for the model problems considered. Based on the concept of weighted Sobolev spaces, we define a novel class of polar function spaces on the parametric square, which provide a natural environment for pull-backs of functions defined on polar domains with corners. Subsequently, we derive error estimates in these spaces and map them again to the physical domain.

Finally, we include the effect of mesh grading into our numerical analysis. The singular functions to be approximated can be described using weighted Sobolev spaces, which are naturally combined with the previously introduced polar Sobolev spaces. This allows us to establish optimal approximation properties of the considered graded polar splines under a simple condition on the grading parameter. Furthermore, we present a series of numerical results to validate the theory and demonstrate the practical efficiency of our graded mesh refinement approach.

The outline of this paper is as follows. Section \ref{section: preliminaries and notation} provides an overview of the mathematical foundations of IGA, which are then applied in Section \ref{section: main concepts and results} to introduce polar parameterizations, the model problem under consideration and the proposed graded mesh refinement scheme. In particular, in Section \ref{subsec: main result}, the key result of our work, an error estimate on graded polar meshes, is stated, establishing optimal convergence rates for suitable mesh grading parameters. In Section \ref{sec: polar function spaces}, we construct a projector on the space of $C^0$-smooth polar splines and define polar function spaces on the parametric domain. Subsequently, in Section \ref{section: Error estimates}, we derive the corresponding projection error estimates, which are needed to prove our main theorem. The theoretical findings are then confirmed numerically in Section \ref{section: Numerical results}. Finally, we conclude the presented contributions and suggest directions for further research in Section \ref{section: Conclusion}.

In the sequel, the symbol $C$ is used for a generic positive constant, which may be different at each occurrence. Moreover, for $a,b >0$ we use the notation $a \sim b$ if there are positive constants $C_1$ and $C_2$ such that $C_1 b \leq a \leq C_2 b$. The mentioned constants are independent of the mesh parameter $h$ and the function under consideration.

\section{Fundamentals of B-splines and NURBS}
\label{section: preliminaries and notation}
This section contains a brief overview of the fundamentals of B-Splines, NURBS and related quasi-interpolants, largely following the framework presented in the review paper \cite{BeiraodaVeigaBuffaSangalliVazquez2014}. For more details, the interested reader may consult the books \cite{Schumaker2007, CottrellHughesBazilevs2009, LycheManniSpeleers2018}.

\subsection{Univariate B-splines and NURBS}
First, we introduce the concept of B-splines and NURBS in the univariate case and provide a brief overview of common refinement procedures. 

Let $p \in \N_0$ and $n \in \N$. We call $\Xi := \{\xi_1, \xi_2, \dots, \xi_{n+p+1}\}$ a $p$-open knot vector if
\begin{align*}
	\xi_1 = \xi_2 = \dots = \xi_{p+1} < \xi_{p+2} \leq \xi_{p+3} \leq \dots \leq \xi_{n-1} \leq \xi_{n} < \xi_{n+1} = \xi_{n+2} = \dots = \xi_{n+p+1} ,
\end{align*}
where $\xi_i \in \R$ for $i =1, \dots, n+p+1$ is called the $i$-th knot which is allowed to occur repeatedly. Without loss of generality, we assume that $\xi_1=0$ and $\xi_{n+p+1}=1$. Furthermore, we define the vector $Z = \{\zeta_1, \dots, \zeta_N\}$ of knots without repetitions, also called breakpoints, with
\begin{align*}
	\Xi = \{\underbrace{\zeta_1, \dots, \zeta_1}_{m_1 \text{ times}}, \underbrace{\zeta_2, \dots, \zeta_2}_{m_2 \text{ times}}, \dots, \underbrace{\zeta_N, \dots, \zeta_N}_{m_N \text{ times}}\} ,
\end{align*}
where $N \in \N$ is the total number of pairwise different knots and $m_j \in \N$ denotes the multiplicity of the breakpoint $\zeta_j$ such that $\sum_{j=1}^N m_j = n + p + 1$. For $p$-open knot vectors, $m_1=m_N=p+1$ always holds, and we assume $m_j \leq p+1$ for the internal knot multiplicities. The entries of $Z$ define a mesh on the unit interval $[0,1]$.

From the given knot vector, B-spline basis functions of degree $p$, denoted by
\begin{align*}
	\widehat{B}_{i,p} \colon [0,1] \to \R , \quad \zeta \mapsto \widehat{B}_{i,p}(\zeta) , \quad  i= 1, 2, \dots, n ,
\end{align*}
can be constructed using the iterative scheme as explained, for instance, in \cite[Section 2.1]{CottrellHughesBazilevs2009}. They build a basis of the space of splines on the subdivision $Z$, that is, piecewise polynomials of degree $p$ with $p-m_j$ continuous derivatives at the internal breakpoints $\zeta_j$, $j=2, \dots, N-1$, where $p-m_j=-1$ stands for a discontinuity at $\zeta_j$. Besides other characteristics, the B-spline basis functions are non-negative and form a partition of unity. Finally, we denote by
\begin{align*}
	S_p(\Xi) = \text{span}\left\{\widehat{B}_{i,p} : i=1,2,\dots,n \right\}
\end{align*}
the univariate spline space that is spanned by the B-splines resulting from the knot vector $\Xi$.

The support of each basis function is given by $\text{supp}(\widehat{B}_{i,p}) = [\xi_i, \xi_{i+p+1}]$. For each interval $I_j = (\zeta_j, \zeta_{j+1})$, which can also be written as $I_j=(\xi_i, \xi_{i+1})$ for a unique $i$, we further define the support extension $\widetilde{I}_j$ by
\begin{align}
	\widetilde{I}_j = (\xi_{i-p}, \xi_{i+p+1}) ,
	\label{eq: univariate support extension}
\end{align}
which represents the interior of the union of the supports of basis functions whose support intersects $I_j$, i.e.,
\begin{align}
	\overline{\widetilde{I}_j} = \bigcup_{\{k \, : \, \text{supp}(\widehat{B}_{k,p}) \cap I_j \neq \emptyset\}} \text{supp}(\widehat{B}_{k,p}) = [\xi_{i-p}, \xi_{i+p+1}] .
\end{align}

Classical splines face restrictions in representing essential geometries like conic sections. To overcome this limitation, non-uniform rational B-splines (NURBS) are introduced, see \cite{PieglTiller1995} for more details. Therefore, the weight function
\begin{align*}
	W(\zeta) = \sum_{l=1}^{n} w_l \, \widehat{B}_{l,p} (\zeta)
\end{align*}
is determined by choosing positive constants $w_l>0$, $l=1,\dots, n$, which are called weights. The NURBS basis functions are then defined by
\begin{align*}
	\widehat{N}_{i,p} (\zeta) = \frac{w_i  \widehat{B}_{i,p} (\zeta)}{\sum_{l=1}^n w_l \, \widehat{B}_{l,p} (\zeta)} = \frac{w_i  \widehat{B}_{i,p} (\zeta)}{W(\zeta)} , \quad i = 1, \dots, n .
\end{align*}
We denote the corresponding NURBS space by
\begin{align*}
	N_p(\Xi, W) = \text{span}\left\{ \widehat{N}_{i,p} : i = 1, \dots, n \right \} .
\end{align*}
The B-spline and NURBS spaces can be refined through knot insertion and degree elevation. In total, three refinement schemes can be constructed by combining the algorithms \cite{HughesCottrellBazilevs2005}.

Next, we introduce a projection operator onto the space of univariate splines that is used in standard literature on isogeometric approximation theory \cite{BazilevsBeiraoDaVeigaCottrellHughesSangalli2006,BeiraodaVeigaChoSangalli2012,BeiraodaVeigaBuffaSangalliVazquez2014}. We define
\begin{align}
	\label{eq: univariate quasi-interpolant}
	\Pi_{p,\Xi} : L^1([0,1]) \to S_p(\Xi), \quad \Pi_{p,\Xi} (v) : = \sum_{i=1}^n \lambda_{i,p}(v) \, \widehat{B}_{i,p} ,
\end{align}
where $\lambda_{j,p}$ are dual functionals with
\begin{align}
	\label{eq: univariate dual functional property}
	\lambda_{i,p}(\widehat{B}_{k,p}) = \delta_{ik} ,
\end{align}
and $\delta_{ik}$ is the Kronecker symbol. In general, $\Pi_{p,\Xi}$ is not an interpolation operator, but it is typically called a quasi-interpolant in the literature. 
The dual basis $\{\lambda_{i,p}\}_{i=1,\dots,n}$ can be chosen in multiple ways \cite{LeeLycheMorken2001}. Here, we stick to the classical construction from \cite[Section 4.6]{Schumaker2007},
\begin{align}
	\lambda_{i,p}(v) = \int_{\xi_j}^{\xi_{j+p+1}} v(s) \, D^{p+1} \psi_i(s) \, \Dd s ,
	\label{eq: dual basis functionals}
\end{align}
where $\psi_i(\zeta) = G_i(\zeta) \Phi_i(\zeta)$ with
\begin{align}
	\label{eq: dual functional definition of Phi}
	\Phi_i(\zeta) = \frac{(\zeta-\xi_{i+1}) \cdots (\zeta - \xi_{i+p})}{p!}
\end{align}
and
\begin{align}
	\label{eq: dual functional definition of G}
	G_i(\zeta) = g\left(\frac{2 \zeta - \xi_i - \xi_{i+p+1}}{\xi_{i+p+1} - \xi_i} \right) ,
\end{align}
where $g$ is the transition function described in \cite[Theorem 4.37]{Schumaker2007}. Note that the dual basis functionals \eqref{eq: dual basis functionals} are well-defined for $v \in L^1([0,1])$, see \cite[Theorem 4.41]{Schumaker2007}.

The presented quasi-interpolant is commonly adapted further to deal with boundary conditions \cite{BeiraodaVeigaBuffaSangalliVazquez2014}, 
\begin{align}
	\label{eq: univariate quasi-interpolant with boundary conditions}
	&\widetilde{\Pi}_{p, \Xi} : C([0,1]) \to S_p(\Xi), \quad \widetilde{\Pi}_{p, \Xi}(v) :=  \sum_{i=1}^n \widetilde{\lambda}_{i,p}(v) \, \widehat{B}_{i,p} \quad \text{ with } \\
	&\widetilde{\lambda}_{1,p}(v) = v(0), \ \widetilde{\lambda}_{n,p}(v) = v(1) , \ \widetilde{\lambda}_{i,p}(v) = \lambda_{i,p}(v), \quad i=2,\dots,n-1 \nonumber .
\end{align}
Since we only consider open knot vectors, the first and last B-spline basis function are interpolatory in $0$ and $1$, respectively, that is,
\begin{align}
	\label{eq: dual functional property boundary univariate 1}
	\widetilde{\lambda}_{1,p}(B_{1,p}) = B_{1,p}(0) = 1 \quad \text{ and } \quad\widetilde{\lambda}_{n,p}(B_{n,p}) = B_{n,p}(1) = 1 .
\end{align}
Moreover, we have
\begin{align}
	\label{eq: dual functional property boundary univariate 2}
	\widetilde{\lambda}_{1,p}(B_{i,p}) = B_{i,p}(0) = 0 \quad \text{ and } \quad\widetilde{\lambda}_{n,p}(B_{i,p}) = B_{i,p}(1) = 0 \quad \text{ for } i = 2,3,\dots,n-1 ,
\end{align}
hence, the modified functionals $\widetilde{\lambda}_{i,p}(v)$ also satisfy the dual functional property \eqref{eq: univariate dual functional property}.

\subsection{Bivariate B-splines and NURBS}
Bivariate B-splines and NURBS are defined as tensor products of their univariate counterparts. This section provides a summary of the main concepts and notation.

Let $n_l \in \N$, the degrees $p_l\in \N$ and the $p_l$-open knot vectors $\Xi_l = \{\xi_{l,1}, \xi_{l,2} \dots , \xi_{l,n_l+p+1}\}$ be given for $l=1,2$. We define the polynomial degree vector $\bs p = (p_1, p_2)$ and the bivariate knot vector $\bs \Xi = \Xi_1 \times \Xi_2$. Further, let $N_l \in \N$ be the number of knots without repetition in the $l$-th direction such that the corresponding univariate knot vectors of breakpoints are given by $Z_l = \{\zeta_{l,1}, \zeta_{l,2} \dots , \zeta_{l,N_l}\}$ for $l=1,2$. They form a Cartesian grid in the parametric domain $\widehat{\Omega} = (0,1)^2$, which defines the parametric B\'ezier mesh $\widehat{\mathcal{M}}$, 
\begin{align}
	\label{eq: parametric Bezier mesh}
	\widehat{\mathcal{M}} := \left \{Q_{\bs j} \subset \widehat{\Omega} : Q_{\bs j} = Q_{(j_1,j_2)} = (\zeta_{1,j_1} , \zeta_{1,j_1 + 1}) \times (\zeta_{2,j_2} , \zeta_{2,j_2 + 1}), \bs j \in \bs J \right \} , 
\end{align}
where we introduce the set of multi-indices $\bs J= \{\bs j = (j_1, j_2) : 1 \leq j_l \leq N_l - 1,  l=1,2\}$. The diameter of each element $Q \in \widehat{\mathcal{M}}$ is denoted by $h_Q$, and the global mesh size of $\widehat{\mathcal{M}}$ is defined as $h = \max\{h_{Q} : Q \in \widehat{\mathcal{M}}\}$.

Bivariate B-spline functions are then defined by
\begin{align*}
	\widehat{B}_{\bs i, \bs p} : [0,1]^2 \to \R , \quad \widehat{B}_{\bs i, \bs p} (\bs \zeta) = \widehat{B}_{i_1,p_1}(\zeta_1) \,  \widehat{B}_{i_2,p_2}(\zeta_2)
\end{align*}
for $\bs i \in \bs I = \{\bs i = (i_1, i_2) : 1 \leq i_l \leq n_l,  l=1,2\}$. The corresponding bivariate NURBS basis functions read
\begin{align}
	\label{eq: relation between NURBS and B-splines}
	\widehat{N}_{\bs i, \bs p} (\bs \zeta) = \frac{w_{\bs i} \, \widehat{B}_{\bs i, \bs p} (\bs \zeta)} {W(\bs \zeta)}
\end{align}
using the weight function
\begin{align}
	W(\bs \zeta) = \sum_{\bs l \in\bs I}  w_{\bs l} \, \widehat{B}_{\bs l, \bs p} (\bs \zeta),
	\label{eq: weight function}
\end{align}
where we choose weights $w_{\bs l} > 0$ for all $\bs l \in \bs I$. The space of NURBS on the parametric domain is finally denoted by
\begin{align}
	\label{eq: parametric NURBS space}
	N_{\bs p}(\bs \Xi, W) = \text{span} \left \{\widehat{N}_{\bs i, \bs p} (\bs \zeta) , \bs i \in \bs I \right \}.
\end{align}

NURBS parameterizations of bivariate geometries are defined as linear combinations of the tensor product functions introduced above,
\begin{align}
	\bs F(\bs \zeta) = \sum_{\bs i \in \bs I} \bs c_{\bs i} \widehat{N}_{\bs i, \bs p} (\zeta) ,
	\label{eq: parameterization definition}
\end{align}
where each basis function is associated with a control point $\bs c_{\bs i} \in \R^2$, $\bs i \in \bs I$. The $\bs F$-image $\Omega = \bs F(\widehat{\Omega})$ of the parametric domain $\widehat{\Omega} = (0,1)^2$ is commonly referred to as the physical domain. Using this construction, exact parameterizations $\bs F \colon  \widehat{\Omega} \to \Omega$ of various types of domains, including geometric primitives like conic sections, can be obtained  \cite{deBoor2001,Farin1995,PieglTiller1995,CottrellHughesBazilevs2009}.

For every B\'ezier element $Q_{\bs j} \in \widehat{\mathcal{M}}$, we introduce its support extension
\begin{align}
	\label{eq: bivariate support extension}
	\widetilde{Q}_{\bs j}:= \widetilde{I}_{1,j_1} \times \widetilde{I}_{2,j_2} ,
\end{align}
where $\widetilde{I}_{l,j_l}$ is the univariate support extension given by \eqref{eq: univariate support extension} for $l=1,2$.

To define a mesh in $\Omega$, we consider the image under $\bs F$ of the partition given by the knot vectors without repetitions, i.e., each element $Q_{\bs j} \in \widehat{\mathcal{M}}$ of the parametric B\'ezier mesh \eqref{eq: parametric Bezier mesh} is mapped to an element $K_{\bs j} =\bs F(Q_{\bs j} )$ in the physical domain. We set
\begin{align*}
	\mathcal{M}:= \left \{ K_{\bs j}  \subset \Omega : K_{\bs j}  = \bs F(Q_{\bs j}), \bs j \in \bs J \right\} ,
\end{align*}
which is commonly known as the physical B\'ezier mesh, or simply B\'ezier mesh.

\label{subsec: bivariate quasi-interpolant}
The univariate quasi-interpolant \eqref{eq: univariate quasi-interpolant with boundary conditions} can be generalized to the bivariate case by a tensor product construction \cite{BeiraodaVeigaBuffaSangalliVazquez2014},
\begin{align*}
	\Pi_{\bs p, \bs \Xi}(v) = \left( \Pi_{p_1, \Xi_1} \otimes \Pi_{p_2, \Xi_2} \right) .
\end{align*}
It can also be expressed in terms of a dual basis,
\begin{align}
	\label{eq: bivariate quasi-interpolant}
	\Pi_{\bs p,\bs \Xi} : L^1([0,1]^2) \to S_{\bs p}( \bs \Xi) , \quad \Pi_{\bs p, \bs \Xi}(v) = \sum_{\bs i \in \bs I} \lambda_{\bs i, \bs p}(v) \widehat{B}_{\bs i, \bs p} ,
\end{align}
where the dual functionals are again given by tensor products \cite[Chapter XVII]{deBoor2001},
\begin{align*}
	\lambda_{\bs i, \bs p}= \lambda_{i_1, p_1} \otimes \lambda_{i_2, p_2} 
\end{align*}
and satisfy
\begin{align}
	\label{eq: dual functional property}
	\lambda_{\bs i, \bs p}(\widehat{B}_{\bs k, \bs p}) = \delta_{\bs i \bs k} .
\end{align}
Here, $\delta_{\bs i \bs k}$ stands for a multi-index Kronecker symbol, i.e.,
\begin{align*}
	\delta_{\bs i \bs k} = \begin{cases} 1 \quad \text{ if } \bs i = (i_1,i_2) = (k_1,k_2) = \bs k , 
		\\ 0 \quad \text{ else} .
	\end{cases}
\end{align*}

This operator is typically modified further for the consideration of boundary conditions. Using the univariate construction \eqref{eq: univariate quasi-interpolant with boundary conditions}, we define
\begin{align*}
	\widetilde{\Pi}_{\bs p, \bs \Xi} = \left( \widetilde{\Pi}_{p_1, \Xi_1} \otimes \widetilde{\Pi}_{p_2, \Xi_2} \right) :  C([0,1]^2) \to  S_{\bs p}( \bs \Xi) .
\end{align*}
Due to \eqref{eq: dual functional property boundary univariate 1} and \eqref{eq: dual functional property boundary univariate 2}, the modified projection satisfies a strong locality property on the boundary $\partial \widehat{\Omega}$ of $\widehat{\Omega} = (0,1)^2$. In more detail, $(\widetilde{\Pi}_{\bs p, \bs \Xi}v)\big|_{\Gamma}$ only depends on $v\big|_{\Gamma}$ for every face $\Gamma$ of $\partial \widehat{\Omega}$. Moreover, we obtain a representation
\begin{align*}
	\widetilde{\Pi}_{\bs p, \bs \Xi}(v) = \sum_{\bs i \in \bs I} \widetilde \lambda_{\bs i, \bs p}(v) \widehat{B}_{\bs i, \bs p}
\end{align*}
with modified functionals $\widetilde \lambda_{\bs i, \bs p} = \widetilde \lambda_{i_1, p_1} \otimes \widetilde \lambda_{i_2, p_2}, \bs i \in \bs I$, which also satisfy the dual functional property \eqref{eq: dual functional property} since all considered knot vectors are open.

\section{Main concepts and results}
\label{section: main concepts and results}
In this section, the main concepts and results of the paper are presented. First, we define polar parameterizations of domains with corners, introduce the model problem, and investigate regularity properties of the solutions. To address the corner singularity, we then explain the proposed graded mesh refinement scheme. Finally, in Section \ref{subsec: main result}, the key results of our work are provided, including error estimates on graded polar meshes that establish optimal convergence rates for suitable grading parameters. 

\subsection{Polar parameterizations}
\label{subsection: Parameterization}

Polar parameterizations in IGA have been subject to multiple research works and several definitions can be found in the literature \cite{TakacsJuettler2011,TakacsJuettler2012, ToshniwalSpeleersHiemstraHughes2017,SpeleersToshniwal2021}. To point out the properties of interest for this paper and fix the notation, we impose the following assumptions:
\begin{enumerate}[label=\textbf{(P.\arabic*)},ref=(P.\arabic*), leftmargin = *]
	\item \label{assumption: knot vectors}The knot vectors are given by a $1$-open knot vector of the form
	\begin{align}
		\label{eq: knot vector Xi_1}
		\Xi_{1} = \left\{0,0,1,1\right\}
	\end{align}
	and a $p_2$-open knot vector of the form
	\begin{align}
		\label{eq: knot vector Xi_2}
		\quad \Xi_2 = \left\{\underbrace{0,\dots,0}_{p_2+1 \text{ times}}, \underbrace{\frac{1}{N_2}, \dots , \frac{1}{N_2}}_{p_2 \text{ times}}, \underbrace{\frac{2}{N_2}, \dots , \frac{2}{N_2}}_{p_2 \text{ times}}, \dots, \underbrace{\frac{N_2-1}{N_2}, \dots , \frac{N_2-1}{N_2}}_{p_2 \text{ times}},\underbrace{1,\dots,1}_{p_2+1 \text{ times}} \right\} 
	\end{align}
	for some $N_2 \in \N$. The required control points and weights can thus be numbered by the index set 
	\begin{align*}
		\bs I = \{ \bs i = (i_1,i_2): 1 \leq i_1 \leq n_1:= 2, 1\leq i_2 \leq n_2 := p_2+1+p_2 \cdot(N_2-1)\} .
	\end{align*}
	\item \label{assumption: control points}The control points $\bs c_{\bs i}$ associated to the left edge of the parametric domain are located at the same point, which is called the polar point $\bs P$ and, w.l.o.g., chosen to be $\bs P = (0,0)^T$. Thus, we have
	\begin{align}
		\label{eq: collapsing control points}
		\bs c_{\bs i_{\pol}} = \bs P = (0,0)^T \text{ for all } \bs i_{\pol} \in \bs I_{\pol} := \{ \bs i =(i_1,i_2) \in \bs I : i_1=1\} ,
	\end{align}
	where $\bs I_{\pol}$ is called the polar index set.
	\item \label{assumption: weights}The weights $w_{\bs i}$ corresponding to the same $i_2$-index are identical,
	\begin{align}
		\label{eq: weights identical}
		w_{1,i_2} = w_{2,i_2} \quad \text{ for all } i_2 = 1,2, \dots, n_2 .
	\end{align}
\end{enumerate}

\begin{definition}[Polar parameterization]
	\label{def: polar parameterization}
	We say that the isogeometric mapping $\bs F: \widehat{\Omega} \to \Omega$ defined by equation \eqref{eq: parameterization definition} is a polar parameterization if it is constructed from knot vectors, control points and weights satisfying assumptions \ref{assumption: knot vectors}, \ref{assumption: control points} and \ref{assumption: weights}. In that context, the resulting domain $\Omega$ is called a polar domain with corner.
\end{definition}

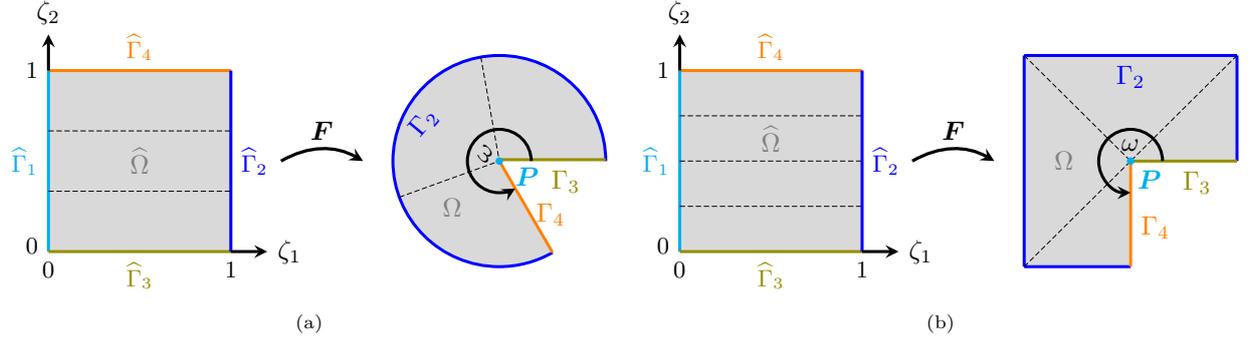
\begin{figure}
	\begin{center}
		\begin{subfigure}{0.495\textwidth}
			\begin{center}
				\input{param_circ.tikz}
				\hspace{0mm}
				\input{phys_circ.tikz}
				\caption{}	
				\label{fig: sketch of circular sector}
			\end{center}
		\end{subfigure}
		\hfill
		\begin{subfigure}{0.495\textwidth}
			\begin{center}
				\input{param_L.tikz}
				\hspace{0mm}
				\input{phys_L.tikz}
				\caption{}	
				\label{fig: sketch of L-shape}
			\end{center}
		\end{subfigure}
	\end{center}
	\vspace{-5mm}
	\caption{Polar parameterizations of exemplary domains with corners and corresponding boundary notation. (a) Circular sector with angle $\omega = \frac53 \pi$, also known as Pacman domain. (b) L-shaped domain.}
	\label{fig: sketch of model domain}
\end{figure}

In Figure \ref{fig: sketch of model domain}, we illustrate exemplary polar parameterizations of a circular sector and an L-shaped domain, which will be discussed in more detail in Examples \ref{ex: circular sector} and \ref{ex: L-Shape}. Beforehand, we point out some useful general conclusions resulting from Definition \ref{def: polar parameterization}. We denote the four boundary edges of the parametric domain by $\widehat{\Gamma}_i$, $i=1,\dots,4$, following the ordering proposed in \cite{Vazquez2016}. Formula \eqref{eq: parameterization definition} and assumptions \ref{assumption: knot vectors} and \ref{assumption: control points} yield that every polar parameterization $\bs F$ is of the form
\begin{align*}
	\bs F(\zeta_1,\zeta_2) 
	= \sum_{i_2=1}^{n_2}  \widehat{N}_{2,p_1}(\zeta_1) \widehat{N}_{i_2,p_2}(\zeta_2) \bs c_{(2,i_2)} 
	= \zeta_1 \sum_{i_2=1}^{n_2} \widehat{N}_{i_2,p_2}(\zeta_2) \bs c_{(2,i_2)} 
	=: \zeta_1 \, \bs{F}^{(2)}(\zeta_2)
\end{align*}
with a function $ \bs{F}^{(2)} = \left({F}^{(2)}_1,{F}^{(2)}_2\right)^T \colon [0,1] \to \Gamma_{2}$ that only depends on $\zeta_2$ and maps into the boundary part away from the polar point, $\Gamma_{2}:= \bs F(\widehat{\Gamma}_2) = \bs F^{(2)}([0,1]) \subset \partial \Omega$. The Jacobian of $\bs F$ satisfies
\begin{align}
	\label{eq: Jacobian of parameterization}
	\det(J_{\bs F}(\zeta_1,\zeta_2)) = \zeta_1 \left(F_1^{(2)}(\zeta_2) {F_2^{(2)}}'(\zeta_2) - {F_1^{(2)}}'(\zeta_2) F^{(2)}_2(\zeta_2) \right) \sim \zeta_1 .
\end{align}
The left edge of the parametric domain $\widehat{\Gamma}_1$ is collapsed into the polar point,
\begin{align}
	\label{eq: collapsing edge}
	\bs F (0,\zeta_2) = \bs P = (0,0)^T \quad \text{ for all } \zeta_2 \in [0,1],
\end{align} 
which, at the same time, represents a corner of the domain, and the corresponding inner angle is denoted by $\omega \in (0,2\pi]$. The remaining parametric boundary edges are mapped onto the edges adjacent to the corner, $\Gamma_3:=\bs F(\widehat{\Gamma}_3)$, and $\Gamma_4:=\bs F(\widehat{\Gamma}_4)$. The parameterization has a polar singularity in the sense that
\begin{align*}
	\frac{\partial \bs F}{\partial \zeta_1} (0,\zeta_2) = (0,0)^T \quad \text{ for all } \zeta_2 \in [0,1]
\end{align*}
and thus 
\begin{align}
	\label{eq: Jacobian of polar parameterization collapses}
	\det(J_{\bs F} (0,\zeta_2) )= 0 \quad \text{ for all } \zeta_2 \in [0,1].
\end{align}
In that context, the edge $\widehat{\Gamma}_1$ will sometimes be called the singular edge. Further, Assumption \ref{assumption: weights} has the effect that the weight function \eqref{eq: weight function} depends solely on $\zeta_2$ and is given by
\begin{align}
	\label{eq: weight function calculated}
	W(\bs \zeta)
	&= \sum_{i_1=1}^{n_1} \sum_{i_2 = 1}^{n_2} w_{(i_1,i_2)} \widehat B_{i_1,p_1}(\zeta_1) \widehat{B}_{i_2, p_2}(\zeta_2) \nonumber \\
	&= \sum_{i_1=1}^{n_1}\widehat B_{i_1,p_1}(\zeta_1)  \sum_{i_2 = 1}^{n_2} w_{(1,i_2)}\widehat{B}_{i_2, p_2}(\zeta_2) \nonumber \\
	&= \sum_{i_2 = 1}^{n_2} w_{(1,i_2)}\widehat{B}_{i_2, p_2}(\zeta_2) \\
	&=: w(\zeta_2) , \nonumber
\end{align}
where we use that the univariate B-splines form a partition of unity. The weight function $W$ and its inverse $W^{-1}$ as well as their derivatives are bounded \cite{BeiraodaVeigaBuffaSangalliVazquez2014}.

The name polar parameterization comes from the similarity to the classical transformation from polar to Cartesian coordinates. In more detail, each polar domain with corner can be parameterized alternatively by a mapping of the form
\begin{align}
	\label{eq: reference transformation}
	\bs G : \widehat{\Omega} \to \Omega , \quad \bs G(r,\varphi): = r R(\Phi(\varphi)) (\cos (\Phi(\varphi)) , \sin(\Phi(\varphi)) )^T ,
\end{align}
which will be called the reference parameterization in polar coordinates in the following. Here, $\Phi:[0,1] \to [0,\omega]$ is a piecewise linear transformation of the isogeometric angle $\varphi \in [0,1]$ to the corresponding polar angle $\Phi(\varphi) \in [0, \omega]$ and $R:[0,\omega] \to \R$ is a function that describes the radius of the domain at each angle. The radius is bounded from below and above, i.e., $R \sim 1$. Hence, the Jacobian of $\bs G$ satisfies
\begin{align}
	\label{eq: Jacobian of reference transformation}
	\det(J_{\bs G}(r,\varphi)) = r (R(\Phi(\varphi)))^2 \Phi'(\varphi) \sim r .
\end{align}

For better understanding, we investigate in more detail the two examples of polar domains with corners depicted in Figure \ref{fig: sketch of model domain}.

\begin{example}[Circular sector]
	\label{ex: circular sector}
	A circular sector of angle $\omega = \frac53 \pi$ can be constructed using the knot vectors $\Xi_1 = \{0,0,1,1\}$ and $\Xi_2 = \left \{0,0,0, \frac13, \frac13, \frac23, \frac23, 1, 1, 1 \right \}$ together with the usual control points and weights from the literature, see for instance \cite{Lu2009,CottrellHughesBazilevs2009}. The geometry mapping from the parametric to the physical domain is illustrated in Figure \ref{fig: sketch of circular sector} and satisfies assumptions \ref{assumption: knot vectors}, \ref{assumption: control points} and \ref{assumption: weights}. An explicit representation of the isogeometric parameterization as well as a comparison with the corresponding reference transformation \eqref{eq: reference transformation} with $R \equiv 1$ and $\Phi(\varphi) = \omega \varphi$ have been derived in \cite{ApelZilk2024} for an exemplary circular sector of angle $\omega=2\pi$.
\end{example}

\begin{example}[L-Shape]
	\label{ex: L-Shape}
	An L-shaped domain can be represented by a polar parameterization using the knot vectors $\Xi_1 = \{0,0,1,1\}$ and $\Xi_2 = \left \{0,0,\frac14, \frac12, \frac34, 1,1 \right \}$, the control points as shown in Table \ref{tab: control points for polar parameterization of L-shape} and weights $w_{\bs i} = 1$ for all $\bs i \in \bs I$. The mapping is illustrated in Figure \ref{fig: sketch of L-shape} and satisfies assumptions \ref{assumption: knot vectors}, \ref{assumption: control points} and \ref{assumption: weights}. The corresponding reference parameterization in polar coordinates \eqref{eq: reference transformation} is defined by
	\begin{align*}
		\Phi(\varphi) = \begin{cases}
			\frac23 \omega \varphi  &\quad \text{ if } \varphi \in \left[0,\tfrac{1}{4}\right]  \\
			\omega\left(\frac43\varphi-\frac16\right) &\quad \text{ if } \varphi \in \left[\tfrac{1}{4}, \tfrac{3}{4}\right] \\
			\omega \left(\frac23\varphi+\frac13\right)  &\quad \text{ if } \varphi \in \left[\tfrac{3}{4},1\right]
		\end{cases}
		\quad \text{ and } \quad   	
		R(\Phi) = \begin{cases} 
			\frac{1}{\cos(\Phi)} &\quad \text{ if } \Phi \in \left[0, \frac{1}{4} \pi\right]   \\
			\frac{1}{\sin(\Phi)} &\quad \text{ if } \Phi \in \left[\frac{1}{4}\pi, \frac{3}{4} \pi\right] \\
			- \frac{1}{\cos(\Phi)} &\quad \text{ if } \Phi \in \left[\frac{3}{4}\pi, \frac{5}{4} \pi\right]  \\
			- \frac{1}{\sin(\Phi)} &\quad \text{ if } \Phi \in \left[\frac{5}{4}\pi, \frac{3}{2} \pi\right]
		\end{cases}
	\end{align*}
	with $\omega = \frac32 \pi$. In Figure \ref{fig: reference transfomation L-shape comparison}, we plot the functions $\varphi \mapsto  R(\Phi(\varphi))\cos(\Phi(\varphi))$ and $\varphi \mapsto R(\Phi(\varphi)) \sin(\Phi(\varphi))$ versus $F^{(2)}_1$ and $F^{(2)}_2$ and in Figure \ref{fig: reference transfomation L-shape difference} we show the corresponding difference. We can see that the two parameterizations coincide up to a maximal difference of $C \approx 0.1$,
	\begin{align*}
		\abs{F^{(2)}_1(\varphi) - R(\Phi(\varphi))\cos(\Phi(\varphi))} \leq C \quad \text{ and } \quad \abs{F^{(2)}_2(\varphi) - R(\Phi(\varphi))\sin(\Phi(\varphi))} \leq C \quad \text{ for all } \varphi \in [0,1] ,
	\end{align*}
	and are even identical in some parts,
	\begin{align*}
		&F^{(2)}_1(\varphi) = R(\Phi(\varphi)) \cos(\Phi(\varphi)) \quad \text{ for } \varphi \in \left[0,\tfrac16 \right] \cup \left[\tfrac12,\tfrac56\right] , \\
		&F^{(2)}_2(\varphi) = R(\Phi(\varphi)) \sin(\Phi(\varphi)) \quad \text{ for } \varphi \in \left[ \tfrac16, \tfrac12 \right] \cup \left[\tfrac56,1\right] .
	\end{align*} 
	
	\begin{table}
		\begin{center}
			\bgroup
			\def\arraystretch{1.1}%
			\setlength\tabcolsep{5pt}
			\begin{tabular}{ c | c c c c c } 
				$\bs c_{(i_1,i_2)}$ & $i_2=1$ & $i_2=2$ & $i_2=3$ & $i_2=4$ & $i_2=5$ \\[0.5ex]
				\hline & & & & \\[-2.5ex]
				$i_1 = 1$ & $(0,0)^T$ & $(0,0)^T$ & $(0,0)^T$ & $(0,0)^T$ & $(0,0)^T$  \\
				$i_1 = 2$ & $(1,0)^T$ & $(1,1)^T$ & $(-1,1)^T$ & $(-1,-1)^T$ & $(0,-1)^T$ \\
			\end{tabular}
			\egroup
			\caption{Definition of control points $\bs c_{\bs i}$, $\bs i \in \bs I$, for a polar parameterization of an L-shaped domain.}
			\label{tab: control points for polar parameterization of L-shape}
		\end{center}
	\end{table}
	\begin{figure}
		\begin{subfigure}{0.495 \linewidth}
			\includegraphics[width=\linewidth, trim=2.2cm 0.7cm 2.2cm 1.4cm, clip]{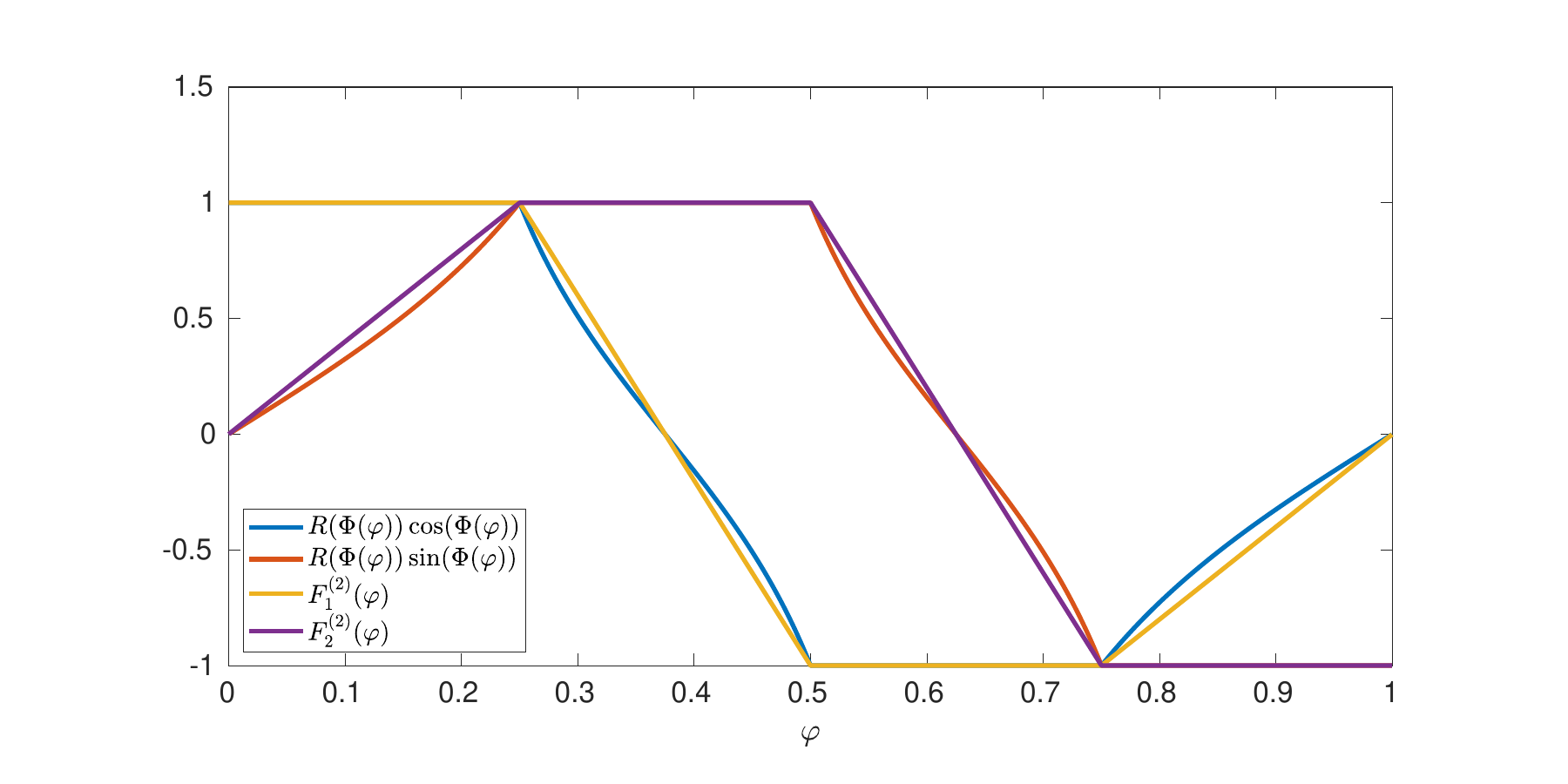}
			\caption{}
			\label{fig: reference transfomation L-shape comparison}
		\end{subfigure}
		\hfill
		\begin{subfigure}{0.495 \linewidth}
			\includegraphics[width=\linewidth, trim=2.2cm 0.7cm 2.2cm 1.4cm, clip]{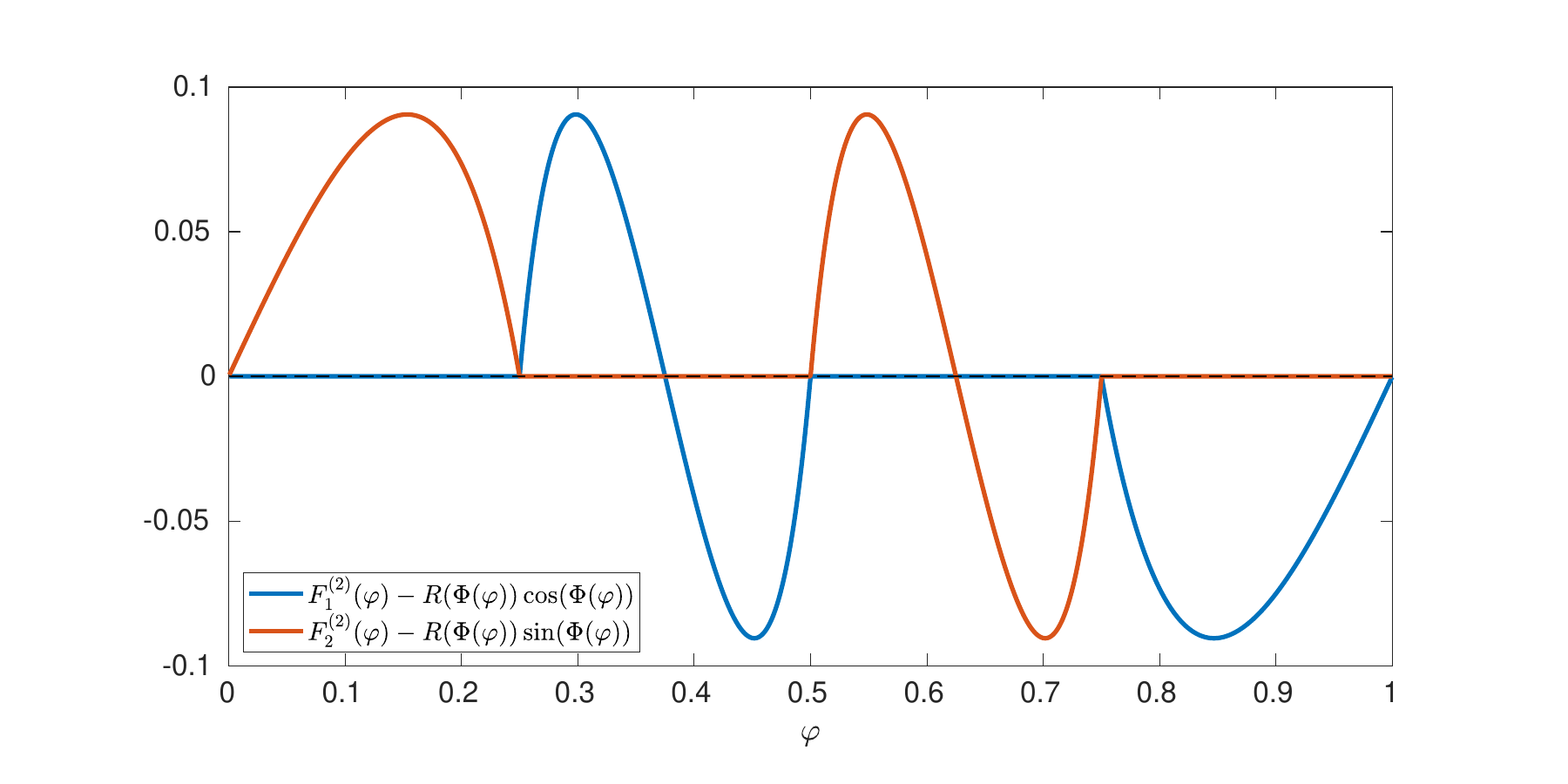}
			\caption{}
			\label{fig: reference transfomation L-shape difference}
		\end{subfigure}
		\vspace{-2mm}
		\caption{(a): Isogeometric mappings $F^{(2)}_1$ and $F^{(2)}_2$ versus scaled polar angular mappings. (b): Difference between the corresponding functions.}
		\label{fig: reference transfomation L-shape}
	\end{figure}	
\end{example}

Finally, we discuss a crucial connection between the polar parameterization $\bs F$ and its corresponding reference transformation $\bs G$. We use standard multi-index notation,
\begin{align*}
	D^{\bs \alpha}:= \frac{\partial ^{\alpha_1}}{\partial x_1 ^{\alpha_1}} \frac{\partial ^{\alpha_2}}{\partial x_2 ^{\alpha_2}} 
	\quad \text{ and } \quad 
	\abs{\bs \alpha} := \alpha_1 + \alpha_2 \quad \text{ with } \bs \alpha = \begin{pmatrix}
		\alpha_1 \\ \alpha_2\end{pmatrix}\in \N_0^2 .
\end{align*}
Given a sufficiently regular function $v \colon \Omega \to \R$, let the derivatives of the pull-backs $\widehat{v} = v \circ \bs F \colon \widehat{\Omega} \to \R$ and $\invbreve{v} = v \circ \bs G \colon \widehat{\Omega} \to \R$ be denoted by $\widehat{D}^{\bs \alpha} \widehat{v}$ and $\invbreve{D}^{\bs \alpha} \invbreve{v}$, respectively. We then write
\begin{align}
	\widehat{D}^{\bs \alpha} \widehat{v} \sim \invbreve{D}^{\bs \alpha} \invbreve{v} , \quad \bs \alpha \in \N_0^2 ,
	\label{eq: equivalence relation isogeometric and polar coordinates}
\end{align}
meaning that $\widehat{D}^{\bs \alpha} \widehat{v}$ and $\invbreve{D}^{\bs \alpha} \invbreve{v}$ can be used in an equivalent way. Similarly, polar and isogeometric coordinates can always be interchanged, and we write $r \sim \zeta_1$ and $\varphi \sim \zeta_2$.

\subsection{Model problem and regularity properties}
\label{subsec: model problems and regularity properties}
In the following, we describe our main model problem and investigate reduced regularity properties of the resulting solutions, which lead to inaccurate approximation with standard methods. Let $\Omega \subset \R^2$ be a polar domain with corner as introduced in Section \ref{subsection: Parameterization} and $\partial \Omega = \Gamma = \Gamma_N \cup \Gamma_D$ its boundary, where the Dirichlet and Neumann boundary are denoted by $\Gamma_D $ and $\Gamma_N$, respectively, with $\Gamma_N \cap \Gamma_D = \emptyset$. We assume non-empty Dirichlet boundary, $\Gamma_{D} \neq \emptyset$, and require each boundary edge to be completely contained in either the Dirichlet or the Neumann boundary, that is, $\Gamma_{\gamma} \subset \Gamma_D$ or $\Gamma_{\gamma}  \subset \Gamma_N$ for $\gamma = 2,3,4$. Our model problem is the Poisson equation
\begin{align}
	\label{eq: poisson equation}
	- \Delta u &= f \quad \text{ in } \Omega , \nonumber \\
	u &= 0 \quad \text{ on } \Gamma_D ,\\
	\frac{\partial u}{\partial n}& = 0 \quad \text{ on } \Gamma_N . \nonumber
\end{align}
For the sake of simplicity, we assume homogeneous boundary conditions. The numerical analysis of inhomogeneous boundary data adds further complexity which is not in the scope of this paper. The weak formulation of problem \eqref{eq: poisson equation} reads: Find $u \in V_0$ such that
\begin{align}
	a( u,  v)	=  b(v)  \quad \forall v \in V_0 ,
	\label{eq: weak formulation Poisson equation}
\end{align}
where we define the spaces
\begin{align}
	V:= H^1(\Omega) \quad \text{ and } \quad
	V_0 := \{v \in  V : v = 0 \text{ on } \Gamma_D \}
	\label{eq: continuous space with bc Poisson equation}
\end{align} 
and the bilinear and linear form
\begin{align*}
	a : V \times V \to \R, \quad & a( u,  v)  := \int_\Omega \nabla u \cdot \nabla v \, \mathrm{d} \bs x, \\
	b : V \to \R, \quad & b(v)  := \int_\Omega f v \,  \mathrm{d} \bs x . \nonumber
\end{align*}
Since both forms are bounded and $a$ is coercive on $V_0$, the variational problem \eqref{eq: weak formulation Poisson equation} has a unique solution $u \in V_0$. Let further $V_{0h} \subset V_0$ be a subspace of finite dimension. The corresponding discrete problem reads: Find $u_h \in V_{0h}$ such that
\begin{align}
	\label{eq: discrete weak formulation Poisson equation}
	a( u_h,  v_h)=  b(v_h)  \quad \forall v_h \in V_{0h} ,
\end{align}
which yields a unique discrete solution $u_{h} \in V_{0h}$.

The regularity of the solution $u$ can be described conveniently by using the concept of weighted Sobolev spaces, which we introduce in the following. For a given domain $G\subset \R^2$, we denote the usual Sobolev spaces on $G$ by $H^{s}(G)$, $s \in \N_0$, with $L^2(G) = H^{0}(G)$. Moreover, let $\mathcal{D}^\prime(G)$ be the space of distributions on $G$. 

\begin{definition}[Weighted Sobolev spaces on the physical domain]
	\label{def: weighted Sobolev norm physical domain}
	Let $G\subset \R^2$ be a bounded domain and let $r = r(\bs x) = \abs{\bs x} = \sqrt{x_1^2+x_2^2}$ be the distance of every point $\bs x =(x_1,x_2) \in G$ to $\bs 0$. We define the weighted Sobolev spaces $H^{s}_{\beta}(G)$ and $V^{s}_{\beta}(G)$ for $s \in \N_0$ and $\beta \in \R$ by
	\begin{align*}
		H^{s}_{\beta}(G) := \{v \in \mathcal{D}^\prime(G) : \norm{v}_{
			H^{s}_{\beta}(G)} < \infty \}
		\quad \text{ and } \quad 
		V^{s}_{\beta}(G) := \{v \in \mathcal{D}^\prime(G) : \norm{v}_{
			V^{s}_{\beta}(G)} < \infty \} ,
	\end{align*}
	respectively, where
	\begin{align*}
		\norm{v}_{H^{s}_{\beta}(G)} 
		&:= \left(\sum_{\abs{\bs \alpha} \leq s} \norm{r^{\beta} D^{\bs \alpha}v}^2_{L^2(G)} \right)^{1/2} 
		= \left(\sum_{\abs{\bs \alpha} \leq s} \int_G \abs{r(\bs x)^{\beta} D^{\bs \alpha}v(\bs x)}^2 \Dd \bs x \right)^{1/2}, \\
		\norm{v}_{V^{s}_{\beta}(G)} 
		&:= \left(\sum_{\abs{\bs \alpha} \leq s} \norm{r^{\beta-s+\abs{\bs \alpha}} D^{\bs \alpha}v}^2_{L^2(G)} \right)^{1/2} 
		= \left(\sum_{\abs{\bs \alpha} \leq s} \int_G \abs{r(\bs x)^{\beta-s+\abs{\bs \alpha}} D^{\bs \alpha}v(\bs x)}^2 \Dd \bs x \right)^{1/2}.
	\end{align*}
	Corresponding seminorms are defined by
	\begin{align*}
		\abs{v}_{H^{s}_{\beta}(G)} := 	\abs{v}_{V^{s}_{\beta}(G)} := 
		\left( \sum_{\abs{\bs \alpha} = s} \norm{r^{\beta} D^{\bs \alpha}v}^2_{L^2(G)} \right)^{1/2}.
	\end{align*}
	The Sobolev spaces on $G$ without weight are denoted by $H^{s}(G) = H^{s}_{0}(G)$ and $ V^{s}(G) =  V^{s}_{0}(G)$ and we write $ L^{2}_{\beta} (G) =  H^{0}_{\beta}(G) =  V^{0}_{\beta}(G)$.
\end{definition}
Such spaces have been studied both in the analytic setting \cite{MazjaPlamenevskii1978,Grisvard1985,Dauge1988} and in the context of finite element methods \cite{ApelSaendigWhiteman1996,MercierRaugel1982}. It is well-known that they form Hilbert spaces and corresponding embedding theorems have been derived. Moreover, for all $s \in \N$ and $\beta \in \R$, the embedding
\begin{align}
	\label{eq: embedding V and H spaces physical domain}
	V^{s}_{\beta}(G) \hookrightarrow H^{s}_{\beta}(G) . 
\end{align}
follows directly.

In general, the solution $u$ of the Poisson equation \eqref{eq: poisson equation} has a singularity of type $r^{\nu}$ in the vicinity of the corner, where the singular exponent $\nu\in \R$ is defined by the boundary conditions and the corner angle $\omega$. Under adequate regularity assumptions on the right-hand side $f$, there is a number $s_0 \geq 2$, where $s_0= \infty$ is allowed, such that
\begin{align*}
	u \in H^{s}_\beta(\Omega) \quad \text{ for } s\leq s_0 \text{ and } s-1> \beta>s- 1- \nu ,
\end{align*}
in particular, $u \in H^s(\Omega) = H^s_0(\Omega)$ for $s \leq s_0$ and $s<1+\nu$. In more detail, it holds $\nu = \frac{\pi}{\omega}$ for Dirichlet or Neumann boundary conditions on both edges adjacent to the corner and $\nu = \frac{\pi}{2\omega}$ for mixed boundary conditions at $r=0$. For sufficiently smooth $f$, the even stronger result
\begin{align}
	\label{eq: weighted regularity of singular part V spaces}
	u \in V^{s}_\beta(\Omega) \quad \text{ for } s \leq s_0 \text{ and }  s-1> \beta>s- 1- \nu
\end{align}
can be established for Dirichlet and mixed boundary conditions at the corner \cite{ApelSaendigWhiteman1996}. In the case of Neumann boundary conditions, the regularity \eqref{eq: weighted regularity of singular part V spaces} holds under the additional assumption that $u$ vanishes at the corner.

\subsection{Graded isogeometric mesh refinement for polar domains with corners}
\label{subsec: construction of graded mesh}
Next, we describe how isogeometric meshes can be refined in a graded way to numerically resolve the corner singularities. The idea of mesh grading toward corners is inspired by comparable approaches for finite elements \cite{ApelSaendigWhiteman1996,ApelNicaise1998,Babuska1970,OganesjanRukhovets1968}. Similar concepts have also been proposed in isogeometric literature \cite{BeiraodaVeigaChoSangalli2012,LangerMantzaflarisMooreToulopoulos2015}, but only for smoothly parameterized domains, leading to non-local refinement along the mapped edges of a patch. In our recent work \cite{ApelZilk2024}, we introduced the grading scheme for polar parameterizations of circular sectors, which will be generalized in the following.

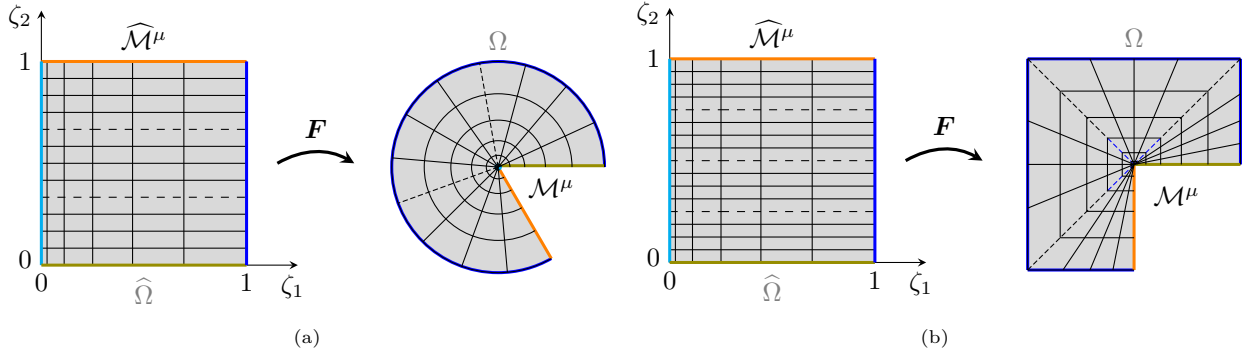
\begin{figure}
	\begin{center}
		\begin{subfigure}{0.495\textwidth}
			\begin{center}
				\input{grad_mesh_circ_par.tikz}
				\hspace{1mm}
				\input{grad_mesh_circ_phys.tikz}
				\caption{}
				\label{fig: graded mesh circular sector}
			\end{center}
		\end{subfigure}
		\begin{subfigure}{0.495\textwidth}
			\begin{center}
				\input{grad_mesh_L_par.tikz}
				\hspace{2mm}
				\input{grad_mesh_L_phys.tikz}
				\caption{}	
				\label{fig: graded mesh L-shape}
			\end{center}
		\end{subfigure}
		\caption{Graded parametric and physical B\'ezier meshes. (a): Circular sector. (b): L-shaped domain.}
		\label{fig: graded meshes}
	\end{center}
\end{figure}

The central idea of our method is to modify the typical uniform $h$-refinement procedure for the knot vector \eqref{eq: knot vector Xi_1}. In this way, we achieve a grading of the resulting B\'ezier mesh towards the corner of the polar domain. In more detail, let
\begin{align*}
	\Xi_1^{h_1} = \{\underbrace{0,\dots,0}_{p_1+1 \text{ times}}, h_1, 2 h_1, \dots, (N_1-3)h_1, (N_1-2)h_1,\underbrace{1,\dots,1}_{p_1+1 \text{ times}}\}
\end{align*} 
be a uniform refinement of $\Xi_1$ obtained by standard degree elevation and knot insertion for $p_1,N_1 \in \N$ with $N_1>1$ and $h_1 = \frac{1}{N_1-1}$. Then, we choose a grading parameter $\mu \in (0,1]$ and define the graded knot vector
\begin{align}
	\label{eq: knot vector graded mesh refinement}
	\Xi_1^{h_1,\mu} = \{\underbrace{0,\dots,0}_{p_1+1 \text{ times}}, (h_1)^{\frac{1}{\mu}}, (2 h_1)^{\frac{1}{\mu}} , \dots, ((N_1-3)h_1)^{\frac{1}{\mu}}, ((N_1-2)h_1)^{\frac{1}{\mu}},\underbrace{1,\dots,1}_{p_1+1 \text{ times}} \} ,
\end{align} 
which we refer to as graded $h$-refinement. Hence, the graded vector of knots without repetitions is given by
\begin{align}
	\label{eq: graded vector of breakpoints}
	Z_1^{h_1,\mu} = \{ \zeta_{1,j_1} := \left((j_1-1) h_1\right)^{\frac{1}{\mu}} : 1 \leq j_1 \leq N_1\} .
\end{align}
In $\zeta_2$-direction, no adjustments are needed, and we employ the standard knot vector $\Xi_2^{h_2}$ and its corresponding vector of breakpoints 
\begin{align}
	\label{eq: uniform vector of breakpoints 2}
	Z_2^{h_2} = \{ \zeta_{2,j_2} := (j_2-1) h_2 : 1 \leq j_2 \leq N_2\} ,
\end{align}
where $h_2 = \frac{1}{N_2-1}$ for $N_2 \in \N$ with $N_2>1$. Consequently, we obtain a graded parametric B\'ezier mesh
\begin{align*}
	\widehat{\mathcal{M}}^\mu
	= \left\{Q_{\bs j} \subset \widehat{\Omega}: Q_{\bs j} =Q_{(j_1,j_2)} =\left(\zeta_{1,j_1}, \zeta_{1,j_1+1}\right) \times \left(\zeta_{2,j_2}, \zeta_{2,j_2+1}\right) , \bs j \in \bs J \right\} ,
\end{align*}
which is locally refined towards the singular edge. The corresponding physical B\'ezier mesh
\begin{align*}
	\mathcal{M}^\mu
	& = \left\{K_{\bs j} \subset \Omega : K_{\bs j} = \bs F (Q_{\bs j}), \bs j \in \bs J \right\}
\end{align*}
is locally refined towards the corner of the polar domain. In Figure \ref{fig: graded mesh circular sector} and \ref{fig: graded mesh L-shape}, we illustrate the graded meshes exemplarily for a circular sector and an L-shaped domain. The following two characteristics distinguish our method from other local refinement approaches:
\begin{enumerate}
	\item Graded meshes allow for local refinement without introducing additional degrees of freedom. The total number of degrees of freedom in the uniform mesh is identical to that of the graded mesh.
	\item Graded meshes towards a polar point enable local refinement near a corner while preserving the tensor product structure of the parametric mesh. To our knowledge, this can not be achieved by any other boundary-fitted approach in IGA, see for instance the overview on adaptive IGA in \cite{BuffaGantnerGiannelliPraetoriusVazquez2022}.
\end{enumerate}
Furthermore, setting $\mu=1$ corresponds to uniform refinement, thus the proposed scheme is a generalization of standard knot insertion. To simplify the notation, the superscript $\mu$ will sometimes be omitted throughout the paper.

\subsection{Main results}
\label{subsec: main result}
In this section, we state the main results of our paper, given by error estimates on polar domains with corners, which will then be proven in the subsequent sections. In particular, we establish optimal approximation rates on graded meshes with suitable grading parameters.

Let $\Omega$ be a polar domain with corner and let $\widehat{V}_h = N_{\bs p}(\bs \Xi, W)$ be a standard NURBS space on $\widehat \Omega$, which has been obtained by graded refinement of the initial polar discretization, following an isoparametric approach. For simplicity, we assume $h_1 \sim h_2$ for the refinement parameter in both parametric directions and set $h=\max\{h_1,h_2\}$. Moreover, $\bs p = (p_1,p_2)$ denotes the bivariate NURBS degree of the refined approximation space and $p = \min\{p_1,p_2\}$. The lack of regularity in the polar parameterization results in a variational crime when standard isogeometric approximation spaces
\begin{align}
	\label{eq: standard approximation space}
	V_h = \left \{f \circ \bs F^{-1} : f \in N_{\bs p}(\bs \Xi, W) \right \} = \text{span} \left\{ N_{\bs i, \bs p} (\bs x) := \widehat{N}_{\bs i, \bs p} \circ \bs F^{-1}(\bs x) , \bs i \in \bs I \right\}
\end{align}
are used since the basis functions associated to the control points collapsing in the polar point,
\begin{align}
	\label{eq: singular basis functions}
	\left\{ N_{\bs i_{\pol}, \bs p}  = \widehat{N}_{\bs i_{\pol}, \bs p} \circ \bs F^{-1}, \bs i_{\pol} \in \bs I_{\pol} \right \} \quad \text{ with } \bs I_{\pol} \text{ from } \eqref{eq: collapsing control points} ,
\end{align}
are not in $H^1(\Omega)$ and thus $V_h \not \subset V = H^1(\Omega)$, see \cite{TakacsJuettler2011,TakacsJuettler2012}. Besides, the functions \eqref{eq: singular basis functions} are not well-defined at the polar point since the corresponding parametric functions are not constant on the singular edge $\widehat{\Gamma}_1$ and consequently, it also holds $V_h \not \subset C^0(\overline{\Omega})$. Instead, the modified polar approximation space
\begin{align}
	\label{eq: modified approximation space}
	V_h^{\pol} : = V_h \cap H^1(\Omega) 
\end{align}
which was initially proposed in \cite{TakacsJuettler2011}, see also \cite{ToshniwalSpeleersHiemstraHughes2017}, is employed. It is well-known that definition \eqref{eq: modified approximation space} yields a space of NURBS that are $C^0$-continuous at the polar point $\bs P$. Different approaches have been proposed in the literature to generate higher continuity as well \cite{ToshniwalSpeleersHiemstraHughes2017,SpeleersToshniwal2021}, but they are not required by the second-order model problem and thus not considered in this paper.

The main result of our paper is the optimal approximation power of the space $V_h^{\pol}$ with respect to the mesh parameter $h$ for the approximation of possibly singular functions on $\Omega$, provided that the grading parameter $\mu$ is chosen correctly. We formulate this in the following theorem.

\begin{Theorem}
	\label{theorem: main result}
	Let $q \in \{0,1\}$ and $s, s_0 \in \N$ with $2 \leq s \leq s_0 \leq p +1$. Further, let $v \in V^{s}_{\beta}(\Omega)$ for all $s-1>\beta >s-1-\nu$, recall \eqref{eq: weighted regularity of singular part V spaces}, and let the mesh grading parameter $\mu \in (0,1]$ satisfy the condition
	\begin{align}
		\label{eq: mesh grading parameter condition}
		\mu <  \frac{\nu-q+1}{s-q}.
	\end{align}
	Then, for every $h>0$ and $s\leq s_0$, there is a NURBS function $v_h \in V_h^{\pol}$ on the corresponding graded mesh $\mathcal{M}^{\mu}$ such that
	\begin{align*}
		\norm{v - v_h}_{H^q(\Omega)} \leq C h^{s-q} \norm{v}_{V^{s}_{\beta}(\Omega)} .
	\end{align*}
\end{Theorem}

Since the approximation error in the $H^1$-norm decreases with the same order as the projection error, we can immediately deduce approximation error estimates in the $L^2$- and $H^1$-norm for the model problem introduced in Section \ref{subsec: model problems and regularity properties} by requiring the grading parameter to satisfy \eqref{eq: mesh grading parameter condition} for $q=1$. In fact, it is even possible to establish optimal convergence of the Galerkin approximation in the $L^2$-norm under a weaker condition on the grading parameter, but this is beyond the scope of the present paper.

\begin{corollary}
	\label{cor: approximation error}
	Let $u$ be the solution of the variational model problem \eqref{eq: weak formulation Poisson equation} such that the regularity property \eqref{eq: weighted regularity of singular part V spaces} is satisfied. Further, let $u_h \in V_h^{\pol}$ be the solution of the discrete weak problem \eqref{eq: discrete weak formulation Poisson equation}, where the discretization is based on the graded mesh refinement scheme presented in Section \ref{subsec: construction of graded mesh}. If the mesh grading parameter $\mu \in (0,1]$ satisfies
	\begin{align*}
		\mu < \frac{\nu}{p} ,
	\end{align*}
	we obtain optimal convergence with respect to the mesh size $h$ in the energy norm and the $L^2$-norm,
	\begin{align*}
		\norm{u - u_h}_{H^q(\Omega)} \leq C h^{p+1-q} \norm{u}_{V^{p+1}_{\beta}(\Omega)} \quad \text{for all } p \leq s_0-1 \text{ and } q \in \{0,1\} .
	\end{align*}
\end{corollary}
\begin{proof}
	For $q=1$, the assertion follows from C\'ea's lemma and Theorem \ref{theorem: main result}. For $q=0$, a standard duality argument can be applied.
\end{proof}

We finish this chapter with a remark on the explicit choice of grading parameter $\mu$. Equation \eqref{eq: mesh grading parameter condition} in Theorem \ref{theorem: main result} only states a necessary condition for optimal convergence, but for an actual computation it arises the question how to choose $\mu$ precisely. This has been discussed several times in the literature, see for instance \cite{ApelMelenk2017}. Based on the numerical studies in \cite{ApelZilk2024}, where the Laplace eigenfunctions of circular sectors, which possess comparable regularities as the solutions of our model problem \eqref{eq: poisson equation}, are approximated using different grading parameters, we choose
\begin{align}
	\label{eq: choice of grading parameter}
	\mu = \min \left\{ 0.9 \frac{\nu}{p} ,1 \right\}.
\end{align}

\subsection{Properties of graded meshes}
In this section, we deduce essential properties of the constructed graded isogeometric meshes. The parametric and physical B\'ezier mesh both contain anisotropic elements \cite{ApelZilk2024}. Therefore, instead of using the seminal approximation framework described in \cite{BazilevsBeiraoDaVeigaCottrellHughesSangalli2006}, which is restricted to uniform meshes, we orientate ourselves on the anisotropic theory from \cite{BeiraodaVeigaChoSangalli2012,BeiraodaVeigaBuffaSangalliVazquez2014}. To begin with, we show that graded parametric meshes satisfy a local quasi-uniformity assumption, see \cite[Assumption 4.10]{BeiraodaVeigaBuffaSangalliVazquez2014}. In the following, we denote the edge lengths of each element $Q_{\bs j} \in \widehat{\mathcal{M}}^{\mu}$ by $h_{1,j_1}$ and $h_{2,j_2}$. We recall that $h_1$ and $h_2$ are the refinement parameters in each direction defined in Section \ref{subsec: construction of graded mesh}, and $h = \max\{h_{Q} : Q \in \widehat{\mathcal{M}}^{\mu}\}$ is the global mesh size of $\mathcal{M}^{\mu}$. 

\begin{lemma}
	\label{lemma: local quasi-uniformity of the graded mesh}
	The graded parametric B\'ezier mesh $\widehat{\mathcal{M}}^{\mu}$ is locally quasi-uniform for all $\mu \in (0,1]$, i.e., there are constants $\theta_l \geq 1$, $l=1,2$, such that
	\begin{align}
		\label{eq: locally quasi-uniform mesh assumption}
		\theta_l^{-1} \leq h_{l,j_l} / h_{l,j_l+1} \leq \theta_l, \quad j_l = 1,2, \dots, N_l-2 ,
	\end{align} 
	in particular, $\theta_1 = 2^{1/\mu}-1$ and $\theta_2=1$.
\end{lemma}
\begin{proof}
	In $\zeta_2$-direction, the uniform vector of breakpoints \eqref{eq: uniform vector of breakpoints 2} leads to a constant edge length,
	\begin{align*}
		h_{2,j_2} &= \zeta_{2,j_2+1} - \zeta_{2,j_2} = (j_2 \, h_2) - ((j_2-1) \, h_2) = h_2, 
		\quad j_2= 1,2, \dots, N_2-1, 
	\end{align*}
	and thus condition \eqref{eq: locally quasi-uniform mesh assumption} is fulfilled for $l=2$ with $\theta_2 = 1$. In $\zeta_1$-direction, the graded vector of knots without repetition \eqref{eq: graded vector of breakpoints} yields
	\begin{align}
		\label{eq: edge length parametric element}
		h_{1,j_1} =  \zeta_{1,j_1+1} -  \zeta_{1,j_1}  
		= \left(j_1^{1/\mu} - (j_1-1)^{1/\mu}\right) \, h_1^{1/\mu} ,
		\quad j_1=1, \dots, N_1-1,
	\end{align}
	and we define
	\begin{align*}
		\Theta_1(j_1) := \frac{h_{1,j_1}}{h_{1,j_1+1}} = \frac{\left(j_1^{1/\mu} - (j_1-1)^{1/\mu}\right) \, h_1^{1/\mu}}{\left((j_1+1)^{1/\mu} - j_1^{1/\mu}\right) \, h_1^{1/\mu}}
		= \frac{j_1^{1/\mu} - (j_1-1)^{1/\mu}}{(j_1+1)^{1/\mu} - j_1^{1/\mu}}.
	\end{align*}
	In Lemma \ref{lemma: function monotonous} of the Appendix, we show that the sequence $(\Theta_1(j_1))_{j_1 = 1}^\infty$ is monotonously increasing and satisfies $\lim_{j_1 \to \infty}\Theta_1(j_1) = 1$. Furthermore, we have $\Theta_1(1)=\frac{1}{2^{1/\mu}-1}\in(0,1]$, and we obtain the desired local quasi-uniformity property \eqref{eq: locally quasi-uniform mesh assumption} for $l=1$ with $\theta_1 = 2^{1/\mu} - 1$.
\end{proof}

Throughout this paper, we will refer to the parameter $\bs \theta = (\theta_1,\theta_2)$ in the sense of Lemma \ref{lemma: local quasi-uniformity of the graded mesh}. In the following, we compute further quantities to describe the graded meshes in more detail. Due to property \eqref{eq: locally quasi-uniform mesh assumption}, it holds
\begin{align*}
	h_{1,j_1} \leq  \theta_1 h_{1,j_1-1} \leq \dots \leq \theta_1^{j_1-1} h_{1,1}.
\end{align*}
Hence, for elements $Q_{\bs j}  \in \widehat{\mathcal{M}}^{\mu}$ that are located close to the singular edge, e.g.\ for $j_1\leq p_1+1$, it holds 
\begin{align*}
	h_{1,j_1} \leq C h_{1,1} = C h_1^{1/\mu} .
\end{align*}
If $Q_{\bs j} \in \widehat{\mathcal{M}}^{\mu}$ is not adjacent to the singular edge, e.g.\ $j_1\geq 2$, the mean value theorem yields the existence of a point $x_0 \in ((j_1-1) h_1,j_1 h_1)$ such that
\begin{align*}
	h_{1,j_1} = (j_1 h_1)^{\frac{1}{\mu}} - ((j_1-1) h_1)^{\frac{1}{\mu}} 
	= \left(\frac{1}{\mu} -1\right) h_1 (x_0)^{\frac{1}{\mu} -1}.
\end{align*}
As the function $x \to x^{1/\mu}$ is monotonously increasing, it follows by definition \eqref{eq: graded vector of breakpoints} and Lemma \ref{lemma: local quasi-uniformity of the graded mesh} that
\begin{align*}
	h_{1,j_1} = C h_1 (x_0)^{\frac{1}{\mu} -1}  
	\leq C h_1  (j_1 h_1)^{\frac{1}{\mu} (1- \mu)}  
	= C h_1  (\zeta_{1,j_1+1})^{1 - \mu} 
	\leq C h_1 (\zeta_{1,j_1})^{1 - \mu} \quad \text{ for } j_1\geq 2 .
\end{align*}
A similar way to prove this property is demonstrated in \cite[page 143]{Apel1999}.
Furthermore, we denote the edge lengths of the support extension $\widetilde{Q}_{\bs j}$ of $Q_{\bs j} \in \widehat{\mathcal{M}}^{\mu}$ by $\widetilde{h}_{1,j_1}$ and $\widetilde{h}_{2,j_2}$. Due to the local quasi-uniformity of the mesh, we obtain similar relations,
\begin{alignat}{2}
	\label{eq: edge lengths parametric support extension}
	&\widetilde{h}_{1,j_1} \leq C \, h_1^{1/\mu} &&\quad \text{for elements close to the singular edge}, \nonumber \\
	&\widetilde{h}_{1,j_1} \leq C \, h_1 \, (\zeta_{1,j_1})^{1 - \mu} &&\quad \text{for elements away from the singular edge}, \\
	&\widetilde{h}_{2,j_2}  \leq C \, h_2 \, \nonumber &&\quad \text{for all elements}.
\end{alignat} 

\section{Function spaces and projectors for polar parameterizations}
\label{sec: polar function spaces}
In this section, we describe suitable isogeometric approximation spaces on polar domains with corners and introduce a corresponding projector. Furthermore, we present the concept of polar Sobolev spaces on the parametric domain, which are defined as pull-backs of classical Sobolev spaces on the physical domain.

\subsection{Polar splines and NURBS}
We start by introducing polar spline and NURBS spaces, which have been studied extensively in isogeometric literature \cite{TakacsJuettler2011,TakacsJuettler2012,ToshniwalSpeleersHiemstraHughes2017,SpeleersToshniwal2021}, and define a suitable projector into these spaces.

\subsubsection{Basis functions of the modified approximation space}
\label{subsection: modified approximation spaces}	
\begin{figure}
	\begin{center}
		\begin{subfigure}{0.79\linewidth}
			\includegraphics[width=0.24\linewidth, trim=6.5cm 3cm 6.5cm 6cm, clip]{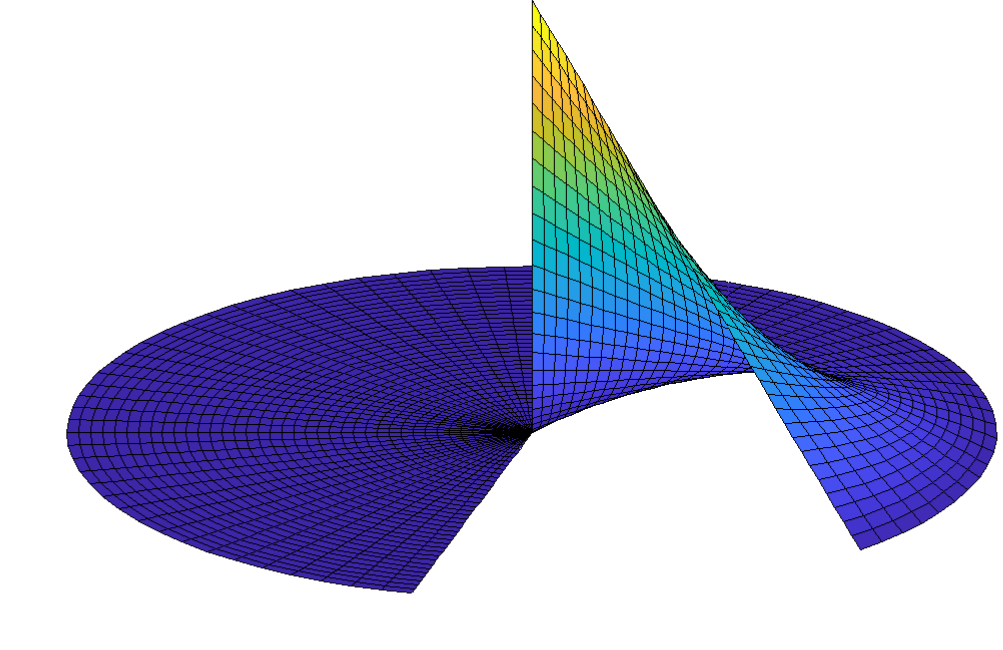}%
			\includegraphics[width=0.24\linewidth, trim=6.5cm 3cm 6.5cm 6cm, clip]{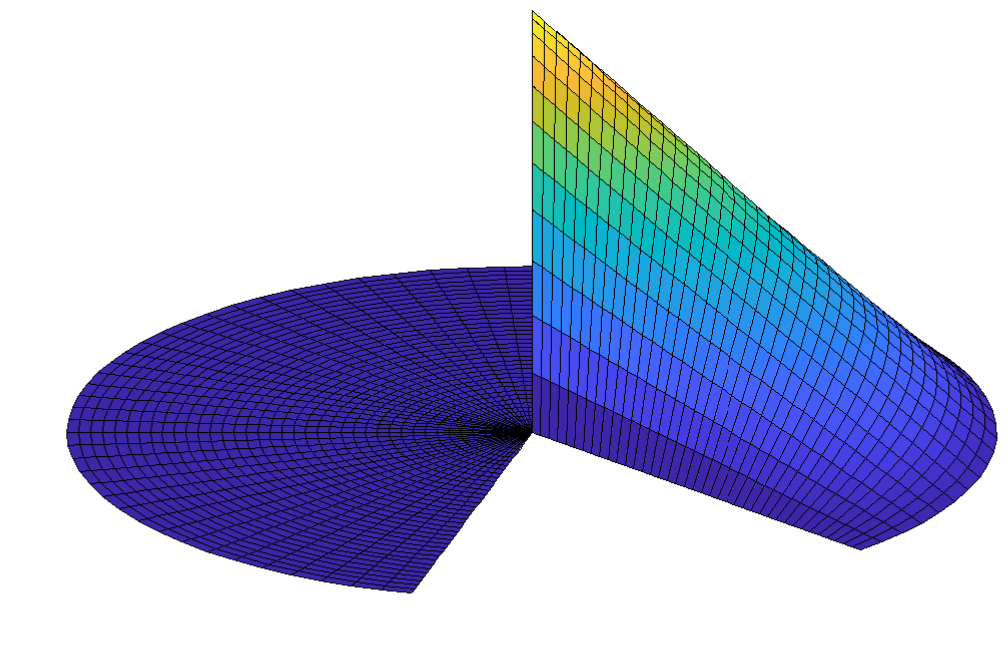}%
			\includegraphics[width=0.24\linewidth, trim=6.5cm 3cm 6.5cm 6cm, clip]{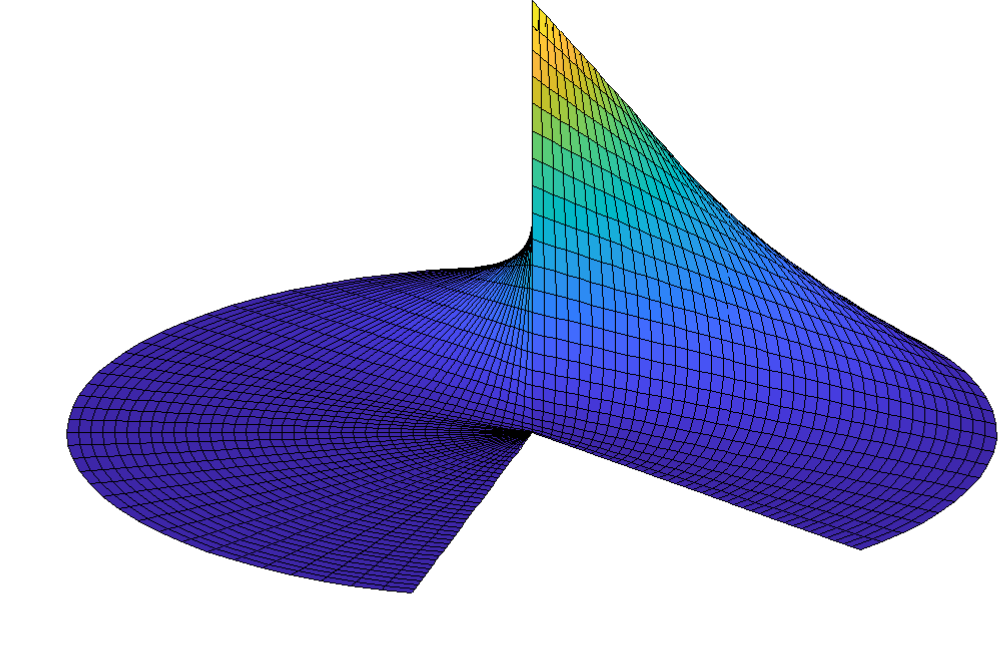}%
			\includegraphics[width=0.24\linewidth, trim=6.5cm 3cm 6.5cm 6cm, clip]{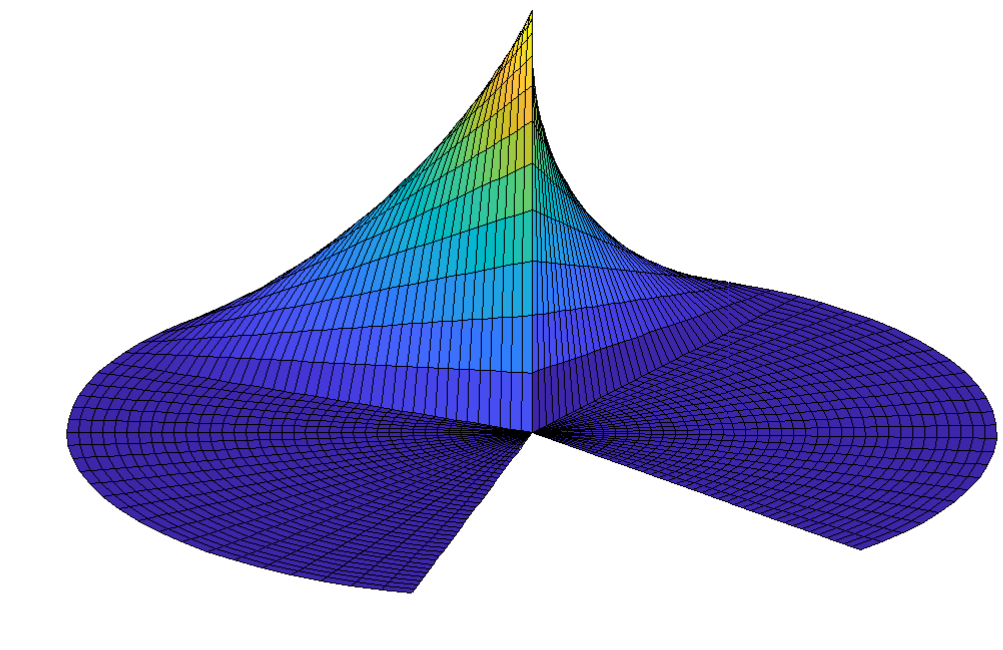}%
			
			\includegraphics[width=0.24\linewidth, trim=6.5cm 3cm 6.5cm 6cm, clip]{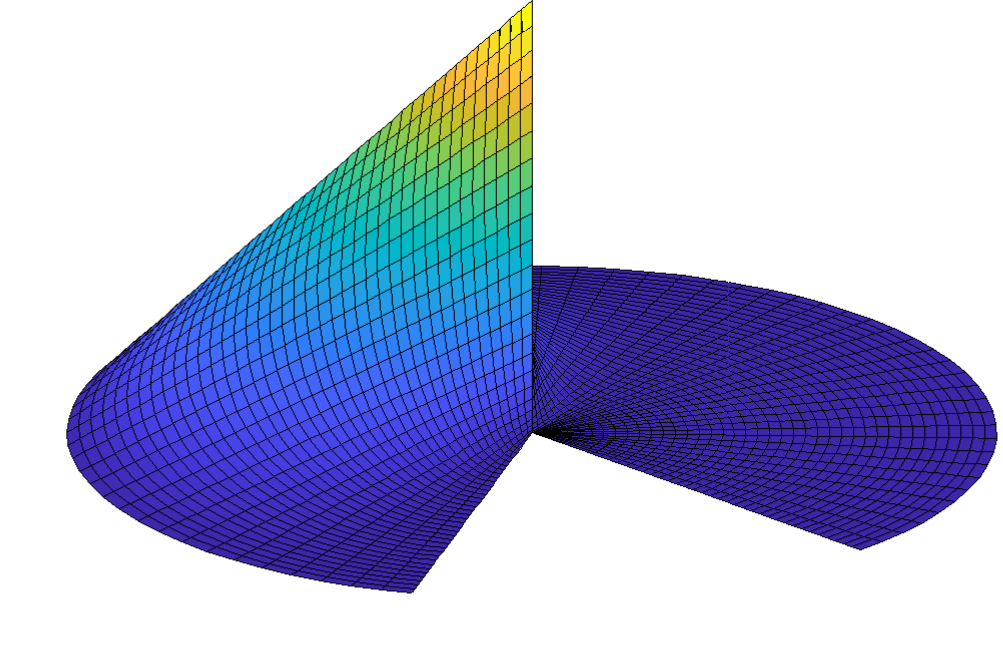}%
			\includegraphics[width=0.24\linewidth, trim=6.5cm 3cm 6.5cm 6cm, clip]{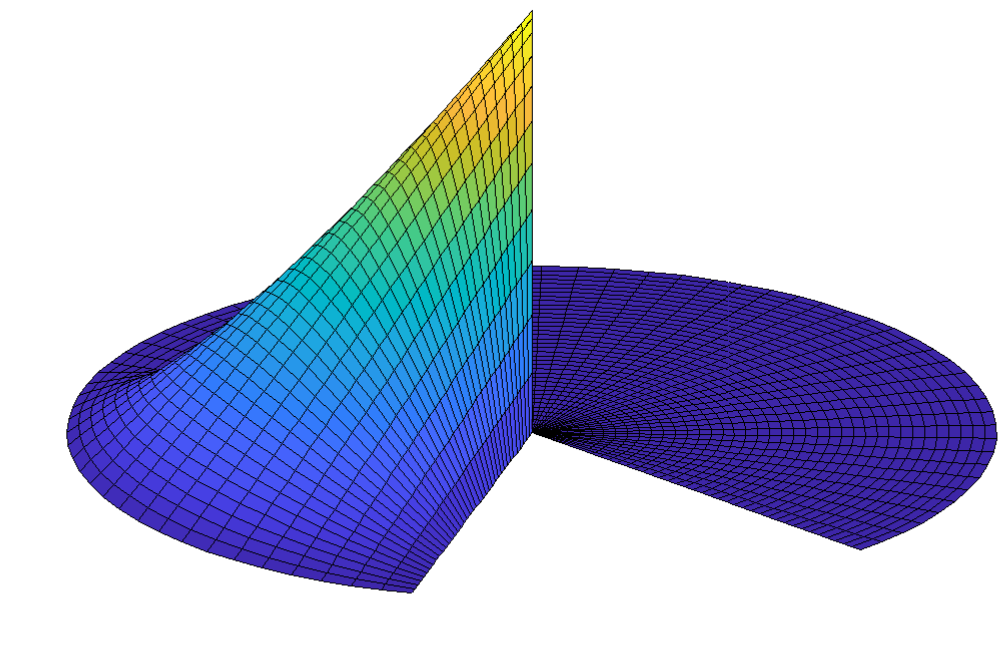}%
			\includegraphics[width=0.24\linewidth, trim=6.5cm 3cm 6.5cm 6cm, clip]{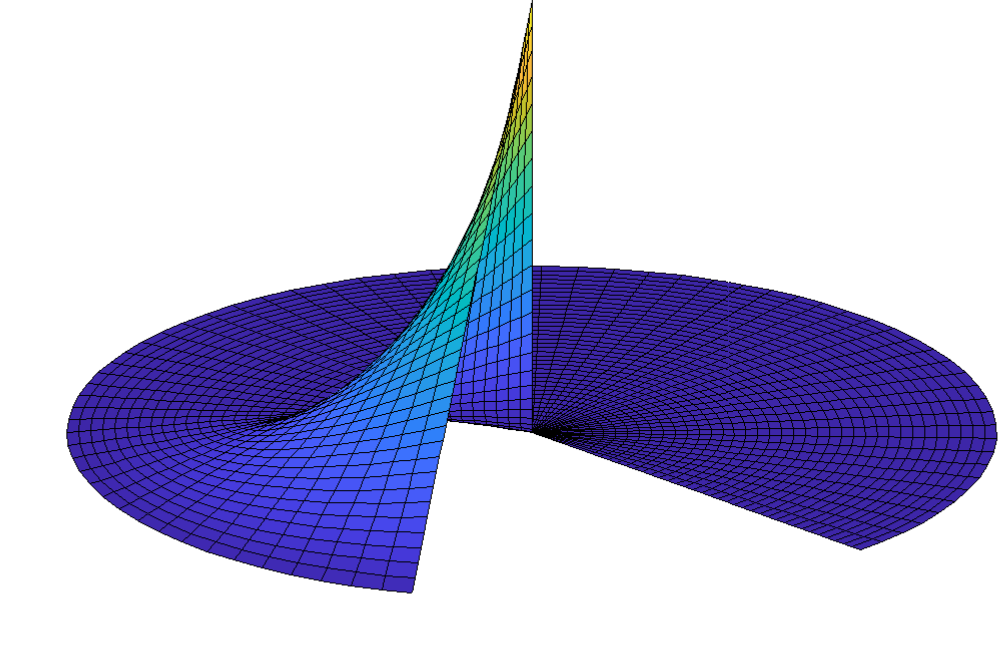}		
			\caption{}
			\label{fig: singular basis functions}
		\end{subfigure}%
		\hfill
		\begin{subfigure}{0.2\textwidth}
			\begin{center}
				\includegraphics[width=\linewidth,trim=6.5cm 3cm 6.5cm 6cm, clip]{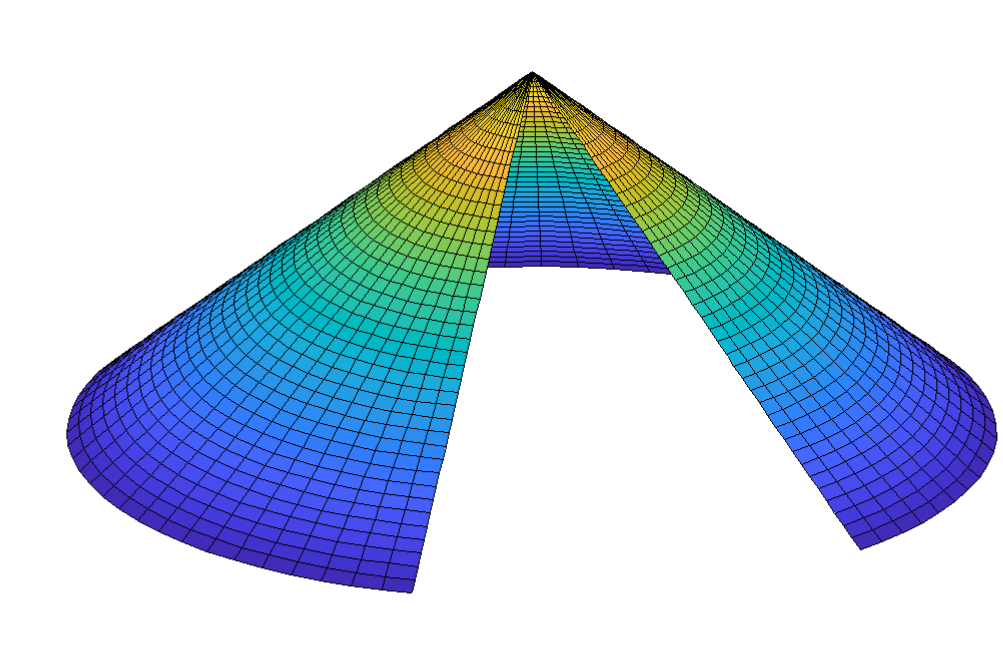}
			\end{center}
			\vspace{-3mm}
			\caption{}
			\label{fig: modified basis function}
		\end{subfigure}
	\end{center}
	\vspace{-5mm}
	\caption{Basis functions for a coarse discretization of a circular sector sector with angle $\omega = \frac{5}{3} \pi$. Some mesh lines that do not correspond to the coarse discretization are inserted for better visualization. (a) Standard singular basis functions \eqref{eq: singular basis functions}. (b): Modified basis function \eqref{eq: modified basis function physical domain}.}
	\label{fig: singular vs modified basis function}
\end{figure}

A basis of the polar spline space \eqref{eq: modified approximation space} can be constructed by replacing all the basis functions \eqref{eq: singular basis functions} with a single function consisting of their sum \cite{TakacsJuettler2011,TakacsJuettler2012},
\begin{align}
	\label{eq: modified basis function physical domain}
	N_{\pol,\bs p} := \sum_{\bs i_{\pol} \in \bs I_{\pol}} N_{\bs i_{\pol}, \bs p} = \sum_{i_2 = 1}^{n_2} N_{(1,i_2), \bs p},
\end{align}
which satisfies $N_{\pol,\bs p} \in H^1(\Omega)$.
In Figure \ref{fig: singular basis functions} and \ref{fig: modified basis function}, we exemplarily display the excluded basis functions \eqref{eq: singular basis functions} and the modified basis function \eqref{eq: modified basis function physical domain}, respectively, for a coarse discretization of a circular sector with angle $\frac53 \pi$. The parametric B-spline and NURBS functions corresponding to the modified physical basis function \eqref{eq: modified basis function physical domain} are given by
\begin{align}
	&\widehat{B}_{\pol,\bs p} \colon [0,1] \to \R, \ \widehat{B}_{\pol,\bs p} (\bs \zeta)_ := \sum_{\bs i_{\pol} \in \bs I_{\pol}} \widehat{B}_{\bs i_{\pol}, \bs p} (\bs \zeta)  \quad \text{ and } \nonumber \\
	&\widehat{N}_{\pol,\bs p} \colon [0,1] \to \R, \ \widehat{N}_{\pol,\bs p} (\bs \zeta)_ := \sum_{\bs i_{\pol} \in \bs I_{\pol}} \widehat{N}_{\bs i_{\pol}, \bs p} (\bs \zeta) ,
	\label{eq: modified parametric NURBS basis function}
\end{align}
respectively. In fact, due to the specific form \eqref{eq: weight function calculated} of the weight function, both functions reduce to the first univariate B-spline basis function in $\zeta_1$-direction,
\begin{align*}
	\widehat{B}_{\pol,\bs p} (\bs \zeta)
	&= \sum_{i_2 = 1}^{n_2} \widehat{B}_{(1,i_2), \bs p}  (\bs \zeta)
	= \widehat B_{1,p_1} (\zeta_1)\sum_{i_2 = 1}^{n_2} \widehat{B}_{i_2, p_2} (\zeta_2)
	= \widehat B_{1,p_1} (\zeta_1), \\
	\widehat{N}_{\pol, \bs p} (\bs \zeta)
	&= \frac{\sum_{i_2 = 1}^{n_2} w_{(1,i_2)} \widehat{B}_{(1,i_2), \bs p}  (\bs \zeta)}{W}
	= \frac{\widehat B_{1,p_1} (\zeta_1) \sum_{i_2 = 1}^{n_2} w_{(1,i_2)} \widehat{B}_{i_2, p_2}(\zeta_2)}{\sum_{i_2 = 1}^{n_2} w_{(1,i_2)} \widehat{B}_{i_2, p_2}(\zeta_2)}
	= \widehat B_{1,p_1} (\zeta_1) .
\end{align*}
Since univariate B-splines constructed from open knot vectors are interpolatory at the boundary, $\widehat{N}_{\pol, \bs p}$ is constant on $\widehat{\Gamma}_1$, with
$\widehat{N}_{\pol, \bs p}(0,\cdot) = \widehat B_{1,p_1} (0) = 1$. Hence, the push-forward of $\widehat{N}_{\pol, \bs p}$ is well-defined and yields the modified NURBS basis function \eqref{eq: modified basis function physical domain},
\begin{align}
	\label{eq: modified basis function as push forward}
	{N}_{\pol, \bs p} = \sum_{\bs i_{\pol} \in \bs I_{\pol}} N_{\bs i_{\pol}, \bs p} 
	= \sum_{\bs i_{\pol} \in \bs I_{\pol}} \widehat{N}_{\bs i_{\pol}, \bs p} \circ \bs F ^{-1} 
	= \widehat{N}_{\pol, \bs p} \circ \bs F ^{-1} = \widehat{B}_{1, p_1} \circ (\bs F ^{-1})_1 ,
\end{align}
where $(\bs F ^{-1})_1$ is the first component of $\bs F ^{-1}$. In that sense, we will write the basis of our approximation space throughout the paper as
\begin{align}
	\label{eq: basis of modified space as push-forwards}
	V_h^{\pol} = \text{span}\left(\left\{N_{\pol, \bs p} =\widehat{B}_{1, p_1} \circ \bs F ^{-1} \right\} \cup \left\{N_{\bs i, \bs p} =\widehat{N}_{\bs i, p_1} \circ \bs F ^{-1} : \bs i \in \bs I \setminus \bs I_{\pol} \right\} \right) \, 
\end{align}
noting that all basis functions are continuous in $\bs P$ with $N_{\pol, \bs p} (\bs P) = 1$ and $  N_{\bs i, \bs p}(\bs P) = 0$ for all $\bs i \in \bs I \setminus \bs I_{\pol}$.

Lastly, we define the space
\begin{align*}
	V_{0h}^{\pol} =  V_{h}^{\pol} \cap V_{0} = \{v \in V_h^{\pol} : v = 0 \text{ on } \Gamma_D\}
\end{align*}
with $V_0$ from \eqref{eq: continuous space with bc Poisson equation}, which also takes homogeneous Dirichlet boundary conditions into account.

\subsubsection{A projector into the modified approximation space}
\label{subsec: projectors in physical domain}
In the next step, we define a projector on the modified isogeometric approximation space $V_h^{\pol}$. In standard literature \cite{BazilevsBeiraoDaVeigaCottrellHughesSangalli2006,BeiraodaVeigaChoSangalli2012,BeiraodaVeigaBuffaSangalliVazquez2014}, a projector on the usual approximation space $V_h$ is defined based on the multivariate quasi-interpolant on the parametric spline space \eqref{eq: bivariate quasi-interpolant}, the parameterization $\bs F$ and the weight function $W$,
\begin{align}
	\label{eq: standard projector physical domain}
	\Pi_{V_h} : L^1(\Omega) \to V_h, \quad \Pi_{V_h} v := \frac{\Pi_{\bs p , \bs \Xi} (W( v \circ \bs F)) }{W} \circ \bs F^{-1} .
\end{align}
However, the isogeometric mapping $\bs F$ is typically assumed to satisfy certain regularity properties, which are not fulfilled by polar parameterizations. In what follows, we propose a slightly modified projector $\Pi_{V_h^{\pol}}$ that is adapted to the polar setting and maps into the subspace $V_h^{\pol} \subset V_h$.

Since we also consider boundary conditions in the polar point, the projection $\Pi_{V_h^{\pol}} v$ of a function $v:\Omega \to \R$ has to be $C^0$-continuous in $\bs P$. In general, the standard projection $\Pi_{V_h} v$ from \eqref{eq: standard projector physical domain} is not well-defined in $\bs P$. Instead, we need a modified quasi-interpolant $\Pi^{\pol}_{\bs p,\bs \Xi}$ that satisfies
\begin{align}
	\label{eq: condition on modified quasi-interpolant}
	\Pi^{\pol}_{\bs p,\bs \Xi} (W( v \circ \bs F)) (0, \cdot) = (W(v \circ \bs F))(0, \cdot) = cW(0,\cdot)
\end{align}
with $c:=v(\bs P) \in \R$. In other words, we require $\Pi^{\pol}_{\bs p,\bs \Xi}$ to satisfy Dirichlet boundary conditions on the singular edge $\widehat{\Gamma}_1$. Motivated by the univariate construction for modified projectors \eqref{eq: univariate quasi-interpolant with boundary conditions}, we define modified dual functionals $\left\{\lambda^{\pol}_{i_1,p_1}\right\}_{\{i_1 = 1, \dots, n_1\}}$ in $\zeta_1$-direction by
\begin{align}
	\label{eq: modified dual functionals}
	\lambda^{\pol}_{1,p_1}(w) = w(0) \quad \text{ and } \quad \ \lambda^{\pol}_{i_1,p_1}(w) = \lambda_{i_1,p_1}(w), \ i_1=2,\dots,n_1,
\end{align}
for $w \in C([0,1])$ and with that the modified quasi-interpolant
\begin{align}
	\label{eq: modified quasi-interpolant}
	\Pi^{\pol}_{\bs p,\bs \Xi} : C([0,1]^2) \to S_{\bs p}( \bs \Xi) , \quad \Pi_{\bs p, \bs \Xi}(v) = \sum_{\bs i \in \bs I} \lambda^{\pol}_{\bs i, \bs p}(v) \widehat{B}_{\bs i, \bs p} ,
\end{align}
where $\lambda^{\pol}_{\bs i, \bs p} = \lambda^{\pol}_{i_1, p_1} \otimes \lambda_{i_2,p_2}$ for all $\bs i \in \bs I$. Note that we have $ \lambda^{\pol}_{\bs i, \bs p} = \lambda_{\bs i, \bs p}$ for all $\bs i \in \bs I \setminus \bs I_{\pol}$. Finally, we can define a suitable projector onto $V_h^{\pol}$ by
\begin{align}
	\label{eq: modified projector physical domain}
	\Pi_{V_h^{\pol}} : C(\overline{\Omega}) \to V^{\pol}_h, \quad 
	\Pi_{V_h^{\pol}} v = \frac{\Pi^{\pol}_{\bs p , \bs \Xi} (W( v \circ \bs F)) }{W} \circ \bs F^{-1} .
\end{align}

\begin{lemma}
	The operator $\Pi_{V_h^{\pol}}$ is a well-defined projector onto $V_h^{\pol}$, that is, $\Pi_{V_h^{\pol}} v \in V_h^{\pol}$ for all $v \in C(\overline{\Omega})$ and $\Pi_{V_h^{\pol}} v = v$ for all $v \in V_h^{\pol}$.
\end{lemma}
\begin{proof}
	Let $v \in C(\overline{\Omega})$. We set $\widehat{v} = v \circ \bs F$ and $c = v(\bs P) = \widehat{v}(0,\cdot)$ and compute the modified dual functionals for every $\bs i_{\pol} = (1,i_2) \in \bs I_{\pol}$ as defined in \eqref{eq: modified dual functionals}, 
	\begin{align*}
		\lambda^{\pol}_{\bs i_{\pol}, \bs p}(W \widehat v) 
		&= (\lambda^{\pol}_{1, p_1} \otimes \lambda_{i_2, p_2})(W \widehat v) = \lambda_{i_2, p_2}(W \widehat v(0,\cdot)) 
		=  \lambda_{i_2, p_2}(W c) = c \,\lambda_{i_2, p_2}(W) \\
		&= c \, \lambda_{i_2, p_2} \left(\sum_{l_2=1}^{n_2} w_{(1,l_2)}\widehat{B}_{l_2,p_2} \right)
		= c \, \sum_{l_2=1}^{n_2} w_{(1,l_2)} \lambda_{i_2, p_2} (\widehat{B}_{l_2,p_2}) 
		= c \, w_{(1,i_2)} = c  \, w_{\bs i_{\pol}} ,
	\end{align*}
	where we used the form \eqref{eq: weight function calculated} of the weight function and the dual basis property \eqref{eq: dual functional property}. It follows
	\begin{align}
		\Pi^{\pol}_{\bs p , \bs \Xi}  (W( v \circ \bs F)) 
		&= \Pi^{\pol}_{\bs p , \bs \Xi}  (W \widehat{v})  = \sum_{\bs i \in \bs I} \lambda^{\pol}_{\bs i, \bs p}(W \widehat v) \widehat{B}_{\bs i, \bs p} \nonumber \\
		&= \sum_{\bs i_{\pol} \in \bs I_{\pol}} \lambda^{\pol}_{\bs i_{\pol}, \bs p}(W \widehat v) \widehat{B}_{\bs i_{\pol}, \bs p}+ \sum_{\bs i \in \bs I \setminus \bs I_{\pol}} \lambda_{\bs i, \bs p}(W \widehat v) \widehat{B}_{\bs i, \bs p} \nonumber \\
		&= \ c \sum_{\bs i_{\pol} \in \bs I_{\pol}} w_{\bs i_{\pol}}\widehat{B}_{\bs i_{\pol}, \bs p} + \sum_{\bs i \in \bs I \setminus \bs I_{\pol}} \lambda_{\bs i, \bs p}(W \widehat v) \widehat{B}_{\bs i, \bs p} \nonumber \\
		&= \ c \, \widehat B_{1,p_1} \sum_{i_2 = 1}^{n_2} w_{(1,i_2)}\widehat{B}_{i_2, p_2}  + \sum_{\bs i \in \bs I \setminus \bs I_{\pol}} \lambda_{\bs i, \bs p}(W \widehat v) \widehat{B}_{\bs i, \bs p} \nonumber \\
		& = c \, \widehat B_{1,p_1} W + \sum_{\bs i \in \bs I \setminus \bs I_{\pol}} \lambda_{\bs i, \bs p}(W \widehat v) \widehat{B}_{\bs i, \bs p} .
		\label{eq: representation in basis}
	\end{align}
	As $\widehat B_{1,p_1}(0) =1$ and $\widehat{B}_{\bs i, \bs p} (0, \cdot) = 0$ for all $\bs i \in \bs I \setminus \bs I_{\pol}$, it is
	\begin{align*}
		\left(\Pi^{\pol}_{\bs p , \bs \Xi}  (W( v \circ \bs F))\right) (0, \cdot) 
		=c \, \widehat B_{1,p_1}(0) W(0,\cdot) +  \sum_{\bs i \in \bs I \setminus \bs I_{\pol}} \lambda^{\pol}_{\bs i, \bs p}(W \widehat v) \widehat{B}_{\bs i, \bs p} (0, \cdot)
		=c \, W(0,\cdot) .
	\end{align*}
	Hence, condition \eqref{eq: condition on modified quasi-interpolant} holds and $\Pi_{V_h^{\pol}} v$ is $C^0$-continuous in $\bs P$. Moreover, by using representation \eqref{eq: representation in basis}, relation \eqref{eq: relation between NURBS and B-splines} and the basis \eqref{eq: basis of modified space as push-forwards}, we obtain
	\begin{align*}
		\Pi_{V_h^{\pol}} v &=
		\left(c \, \widehat B_{1,p_1} + \sum_{\bs i \in \bs I \setminus \bs I_{\pol}} \lambda_{\bs i, \bs p}(W \widehat v) \frac{\widehat{N}_{\bs i, \bs p}}{w_{\bs i}} \right)  \circ \bs F^{-1}
		=  c \, N_{\pol, \bs p} + \sum_{\bs i \in \bs I \setminus \bs I_{\pol}} \frac{\lambda_{\bs i, \bs p}(W \widehat v)}{w_{\bs i}} N_{\bs i, \bs p} 
	\end{align*}
	and thus $\Pi_{V_h^{\pol}} v \in V_h^{\pol}$. 
	
	Next, we show that $\Pi_{V_h^{\pol}}$ is indeed a projector. Therefore, let $v \in V_h^{\pol}$, that is, with \eqref{eq: modified basis function as push forward} and \eqref{eq: basis of modified space as push-forwards},
	\begin{align*}
		v  = c \, N_{\pol,\bs p} + \sum_{\bs i \in \bs I \setminus \bs I_{\pol}} \alpha_{\bs i} N_{\bs i, \bs p} 
		=  c \sum_{\bs i_{\pol} \in \bs I_{\pol}} N_{\bs i_{\pol}, \bs p} + \sum_{\bs i \in \bs I \setminus \bs I_{\pol}} \alpha_{\bs i} N_{\bs i, \bs p},
	\end{align*}
	where $c = v(\bs P)$ and $\alpha_{\bs i} \in \R,\bs i \in \bs I \setminus \bs I_{\pol}$. Then, it is $v = \widehat{v} \circ {\bs F}^{-1}$ with 
	\begin{align*}
		\widehat{v} = c \, \widehat{N}_{\pol, \bs p} + \sum_{\bs i \in \bs I \setminus \bs I_{\pol}}\alpha_{\bs i} \widehat{N}_{\bs i, \bs p} = c \sum_{\bs i_{\pol} \in \bs I_{\pol}} \widehat{N}_{\bs i_{\pol}, \bs p} + \sum_{\bs i \in \bs I \setminus \bs I_{\pol}}\alpha_{\bs i} \widehat{N}_{\bs i, \bs p} \in N_{\bs p}(\Xi),
	\end{align*}
	see \eqref{eq: modified parametric NURBS basis function}. Further, we have $\widehat{v} = \widehat{u}/W$ with 
	\[\widehat{u} = c \sum_{\bs i \in \bs I_{\pol}} w_{\bs i_{\pol}} \widehat{B}_{\bs i_{\pol}, \bs p} + \sum_{\bs i \in \bs I \setminus \bs I_{\pol}}\alpha_{\bs i}w_{\bs i} \widehat{B}_{\bs i, \bs p} \in S_{\bs p}(\Xi) .\]
	Since $\Pi^{\pol}_{\bs p , \bs \Xi}$ is a projector onto $S_{\bs p}(\Xi)$, it follows 
	\begin{align*}
		\Pi_{V_h^{\pol}} v = \frac{\Pi^{\pol}_{\bs p , \bs \Xi} (W( v \circ \bs F)) }{W} 
		= \frac{\Pi^{\pol}_{\bs p , \bs \Xi} (W \widehat v) }{W} 
		= \frac{\Pi^{\pol}_{\bs p , \bs \Xi} (W \widehat{u}/W) }{W} 
		= \frac{\Pi^{\pol}_{\bs p , \bs \Xi} (\widehat{u}) }{W} 
		= \frac{ \widehat{u}}{W} \circ \bs F^{-1} =  \widehat{v} \circ \bs F^{-1} = v .
	\end{align*}
	and the proof is complete.
\end{proof}

\begin{remark}
	Note that all the constructions of Section \ref{subsec: projectors in physical domain} can be extended to a non-isoparametric approach by simply setting $W = 1$.
\end{remark}

\subsection{Polar Sobolev spaces on the parametric domain}
\label{subsec: polar Sobolev spaces}
In standard isogeometric approximation theory, error estimates are typically derived on the parametric domain and then transferred to the physical domain, where the Jacobian of the geometry mapping $\bs F$, as well as its inverse, is assumed to be bounded  \cite{BazilevsBeiraoDaVeigaCottrellHughesSangalli2006,BeiraodaVeigaChoSangalli2012,BeiraodaVeigaBuffaSangalliVazquez2014}. This condition, however, is not satisfied by polar parameterizations, see equation \eqref{eq: Jacobian of polar parameterization collapses}. Therefore, the corresponding transformations must be treated more carefully, as demonstrated in this section. In addition, we introduce a novel class of Sobolev spaces on the parametric domain to characterize pull-backs of functions on polar domains with corners.

\subsubsection{Transformation of norms}
First, we define weighted Sobolev spaces on the parametric domain, which differ from the ones in Definition \ref{def: weighted Sobolev norm physical domain} in the sense that the distance of a point to the edge $\widehat{\Gamma}_1$ is weighted and not the distance to $\bs 0$. That is, for $Q \subset \widehat \Omega$, we set
\begin{align*}
	\widehat H^{s}_{\beta}(Q) := \{\widehat v \in \mathcal{D}^\prime(Q) : \norm{\widehat v}_{\widehat H^{s}_{\beta}(Q)} < \infty \}
	\quad \text{ and } \quad
	\widehat V^{s}_{\beta}(Q) := \{\widehat v \in \mathcal{D}^\prime(Q) : \norm{\widehat v}_{\widehat V^{s}_{\beta}(Q)} < \infty \} ,
\end{align*}
with the norms
\begin{align*}
	\norm{\widehat v}_{\widehat H^{s}_{\beta}(Q)} ^2
	&:= \sum_{\abs{\bs \alpha} \leq s} \norm{r^{\beta} \widehat D^{\bs \alpha}\widehat v}^2_{ \widehat L^2(Q)}
	:= \sum_{\abs{\bs \alpha} \leq s} \int_{Q} \abs{r^{\beta} \widehat D^{\bs \alpha} \widehat v(r,\varphi)}^2 \dr \dphi  , \\
	\norm{\widehat v}_{\widehat V^{s}_{\beta}(Q)} ^2
	&:= \sum_{\abs{\bs \alpha} \leq s} \norm{r^{\beta - s + \abs{\alpha}} \widehat D^{\bs \alpha}\widehat v}^2_{ \widehat L^2(Q)}
	:= \sum_{\abs{\bs \alpha} \leq s} \int_{Q} \abs{r^{\beta- s + \abs{\alpha}} \widehat D^{\bs \alpha} \widehat v(r,\varphi)}^2 \dr \dphi 
\end{align*}
and seminorms
\begin{align*}
	\abs{\widehat v}^2_{\widehat H^{s}_{\beta}(Q)} := \abs{\widehat v}_{\widehat V^{s}_{\beta}(Q)}^2 := 
	\sum_{\abs{\bs \alpha} = s} \norm{r^{\beta} \widehat D^{\bs \alpha} \widehat v}^2_{\widehat L^2(Q)}.
\end{align*}
The Sobolev spaces on $Q$ without weight are denoted by $\widehat H^{s}(Q) = \widehat H^{s}_{0}(Q)$ and $\widehat V^{s}(Q) = \widehat V^{s}_{0}(Q)$ and we write $\widehat L^{2}_{\beta} (Q) = \widehat H^{0}_{\beta}(Q) = \widehat V^{0}_{\beta}(Q)$ and $\widehat L^2(Q) = \widehat H^{0}(Q) = \widehat V^{0}(Q)$. Moreover, similar to \eqref{eq: embedding V and H spaces physical domain}, for all $s \in \N$ and $\beta \in \R$, we have the embedding
\begin{align}
	\label{eq: embedding V and H spaces parametric domain}
	\widehat{V}^{s}_{\beta}(Q) \hookrightarrow \widehat{H}^{s}_{\beta}(Q) . 
\end{align}

Polar parameterizations create the effect that Sobolev norms on the physical domain are converted to different norms on the parametric domain. Let $\widehat{v} = v \circ \bs F$ be the pull-back of a sufficiently regular function $v:\Omega \to \R$.  Moreover, let $K \in \mathcal{M}$ be an element of the B\'ezier mesh, and $Q = \bs F^{-1}(K) \in \widehat{\mathcal{M}}$ its corresponding parametric element. Taking into account the Jacobian \eqref{eq: Jacobian of parameterization} and the weighted spaces defined above, we obtain for $v \in L^2(K)$ by a simple integral transformation the relation
\begin{align}
	\norm{v}_{L^2(K)}^2
	= \int_{Q} \abs{\widehat{v}} \abs{\det(J_{\bs F})} \Dd \bs \zeta 
	\sim \int_{Q} \abs{\widehat{v}} ^2 r \dr \dphi 
	= \norm{\widehat{v}}_{\widehat{L}^2_{1/2}(Q)}^2 ,
	\label{eq: transformation of L2-norm}
\end{align}
where the notation $\sim$ here stands for the equivalence of isogeometric and polar coordinates as introduced in \eqref{eq: equivalence relation isogeometric and polar coordinates}. By transforming the first derivatives of a function on the physical domain to polar coordinates, $\partial_x v \sim \cos \varphi \, \partial_r \widehat v - \frac{1}{r} \sin \varphi \, \partial_\varphi \widehat v$ and $\partial_y v \sim \sin \varphi \, \partial_r \widehat v + \frac{1}{r} \cos \varphi \, \partial_\varphi \widehat v$, we further compute for $v \in H^1(K)$ that
\begin{align}
	\abs{v}^2_{H^1(K)}
	= \int_{K} \abs{\partial_x v}^2 + \abs{\partial_y v} ^2 \Dd \bs x 
	\sim \int_{Q} \left( \abs{\partial_r \widehat{v}}^2 + \abs{\frac{1}{r} \partial_\varphi \widehat{v}} ^2 \right) r \dr \dphi 
	= \norm{\partial_r \widehat{v}}_{\widehat L^2_{1/2}(Q)}^2 + \norm{r^{-1}\partial_\varphi \widehat{v}}_{\widehat L^2_{1/2}(Q)}^2 . 
	\label{eq: transformation of H1-seminorm}
\end{align}
To generalize this pattern, we introduce \textit{weighted polar Sobolev spaces} of order $s \in \N_0$ with weight $\beta \in \R$ on any subset $Q \subset \widehat \Omega$ of the parametric domain,
\begin{align*}
	\widehat{H}^s_{\pol,\beta}(Q) := \left\{ \widehat{v} \in \mathcal{D}'(Q) \colon \norm{\widehat v}_{\widehat{H}^{s}_{\pol,\beta}(Q)} < \infty \right\} 
	\quad \text { and } \quad 
	\widehat{V}^s_{\pol,\beta}(Q) := \left\{ \widehat{v} \in \mathcal{D}'(Q) \colon \norm{\widehat v}_{\widehat{V}^{s}_{\pol,\beta}(Q)} < \infty \right\} ,
\end{align*}
with the norms
\begin{align*}
	\norm{\widehat v}_{\widehat H^{s}_{\pol,\beta}(Q)} ^2
	&:= \sum_{\abs{\bs \alpha} \leq s} \norm{r^{\beta - \alpha_2} \widehat D^{\bs \alpha}\widehat v}^2_{ \widehat L^2_{1/2}(Q)}
	= \sum_{\abs{\bs \alpha} \leq s} \int_{Q} \abs{r^{\beta - \alpha_2} \widehat D^{\bs \alpha} \widehat v(r,\varphi)}^2 r \dr \dphi  , \\
	\norm{\widehat v}_{\widehat V^{s}_{\pol,\beta}(Q)} ^2
	&:= \sum_{\abs{\bs \alpha} \leq s} \norm{r^{\beta - s + \alpha_1} \widehat D^{\bs \alpha}\widehat v}^2_{ \widehat L^2_{1/2}(Q)}
	= \sum_{\abs{\bs \alpha} \leq s} \int_{Q} \abs{r^{\beta- s + \alpha_1} \widehat D^{\bs \alpha} \widehat v(r,\varphi)}^2 r \dr \dphi 
\end{align*}
and seminorms
\begin{align*}
	\abs{\widehat v}^2_{\widehat H^{s}_{\pol,\beta}(Q)} 
	:= \abs{\widehat v}_{\widehat V^{s}_{\pol,\beta}(Q)}^2 
	:= \sum_{\abs{\bs \alpha} = s} \norm{r^{\beta-\alpha_2} \widehat D^{\bs \alpha} \widehat v}^2_{\widehat L^2_{1/2}(Q)}.
\end{align*}
Polar Sobolev spaces without weight will be denoted by $\widehat{H}^s_{\pol}(Q) = \widehat{H}^s_{\pol,0}(Q)$ and $\widehat{V}^s_{\pol}(Q) = \widehat{V}^s_{\pol,0}(Q)$. Furthermore, it holds $\widehat{H}^{0}_{\pol}(Q)= \widehat{V}^{0}_{\pol}(Q) =\widehat{L}^2_{1/2}(Q)$ and $\widehat{H}^{0}_{\pol,\beta}(Q)= \widehat{V}^{0}_{\pol,\beta}(Q) =\widehat{L}^2_{1/2 + \beta}(Q)$.

With this definition, equations \eqref{eq: transformation of L2-norm} and \eqref{eq: transformation of H1-seminorm} can be rewritten as
\begin{align*}
	\norm{v}_{H^0(K)}^2
	\sim \norm{\widehat{v}}_{\widehat{H}^0_{\pol}(Q)}^2 
	= \norm{\widehat{v}}_{\widehat{V}^0_{\pol}(Q)}^2 
	\quad \text{ and } \quad
	\abs{v}^2_{H^1(K)}
	\sim \abs{\widehat{v}}_{\widehat{H}^1_{\pol}(Q)}^2 
	= \abs{\widehat{v}}_{\widehat{V}^1_{\pol}(Q)}^2 .
\end{align*}
We emphasize that the definition of the polar Sobolev spaces intrinsically includes the weight $\frac12$, although this is not explicitly indicated in the notation. Moreover, since the $V$-spaces scale naturally, the following more general transformation can be derived through explicit computation, 
\begin{align}
	\label{eq: transformation of norm on single element adjacent to polar point}
	\norm{v}_{V^{s}_{\beta}(K)} \sim \norm{\widehat{v}}_{\widehat{V}^{s}_{\pol,\beta}(Q)}
	\quad \text{ for all } v \in V^{s}_{\beta}(K) ,
\end{align}
with $s \in \N_0$ and $\beta \in \R$. Additionally, it is easy to see that
\begin{align}
	\label{eq: polar spaces and 1/2-spaces}
	\norm{\widehat{v}}_{\widehat{H}^s_{1/2+\beta}(Q)} 
	\leq C \norm{\widehat{v}}_{\widehat{H}_{\pol,\beta}^s(Q)} 
	\leq C \norm{\widehat{v}}_{\widehat{V}_{\pol,\beta}^s(Q)} 
	\quad \text{ for all } \widehat v \in V_{\pol,\beta}^s(Q) ,
\end{align}
recall also \eqref{eq: embedding V and H spaces parametric domain}, and thus the embeddings $\widehat{V}_{\pol,\beta}^s(Q) 
\subset \widehat{H}_{\pol,\beta}^s(Q) \subset \widehat{H}^s_{1/2+\beta}(Q)$ hold.


\subsubsection{Bent spaces}
\label{subsec: bent spaces}
In general, the Sobolev regularity of a pull-back function is affected not only by the polar singularity, but also by the reduced continuity of the isogeometric parameterization across element interfaces. In particular, polar parameterizations as defined in Section \ref{def: polar parameterization} are only $C^0$-continuous along the mesh lines that correspond to repeated knot values in the knot vector \eqref{eq: knot vector Xi_2}. Therefore, bent Sobolev spaces, as introduced in \cite{BazilevsBeiraoDaVeigaCottrellHughesSangalli2006}, must be taken into account.

Throughout this section, let $E \subset \Omega$ be a union of elements $K \in \mathcal{M}$ and let $\widehat E = \bs F^{-1}(E)\subset \widehat{\Omega}$ be the corresponding union of parametric elements. Besides, for two adjacent parametric elements $Q_1, Q_2 \in \widehat{\mathcal{M}}$, let $m_{Q_1,Q_2}$ denote the maximal order of continuous derivatives of $\bs F$ across their common edge $\partial Q_1 \cap \partial Q_2$. The pull-back $\widehat v = v \circ \bs F$ of a function $v \in H^s_{\beta}(E)$ or $v \in V^s_{\beta}(E)$ is not necessarily an element of $\widehat{H}^s_{\beta}(\widehat{E})$ or $\widehat{V}^s_{\beta}(\widehat{E})$, respectively. Instead, we define \textit{weighted bent Sobolev spaces} $\mathcal{X}^s_{\beta}(\widehat{E}) = \mathcal{H}^s_{\beta}(\widehat{E})$ and $\mathcal{X}^s_{\beta}(\widehat{E}) = \mathcal{V}^s_{\beta}(\widehat{E})$ of order $s \in \N_0$ with weight $\beta \in \R$ by 
\begin{align*}
	\mathcal{X}^s_{\beta}(\widehat{E}) := \left \{
	\widehat v \in \mathcal{D}'(\widehat E) \colon \begin{rcases}
		\begin{dcases} \widehat v_{|Q} \in \widehat X^s_{\beta}(Q) & \text{for all } Q \in \widehat{E} \text{ and } \\  \widehat \nabla^q(\widehat v_{|Q_{1}}) = \widehat \nabla^q (\widehat v_{|Q_2}) & \text{on } \partial Q_1 \cap \partial Q_2 \text{ for } q=0,1,\dots, \min\{m_{Q_1, Q_2},s-1\} \\ & \text{for all }  Q_1, Q_2 \text{ with } \partial Q_1 \cap \partial Q_2 \neq \emptyset \end{dcases} \end{rcases} \right\} ,
\end{align*}
where we insert $\widehat X^s_{\beta}(Q) = \widehat H^s_{\beta}(Q)$ and $\widehat X^s_{\beta}(Q) = \widehat V^s_{\beta}(Q)$, respectively, and $\widehat \nabla^q$ denotes the $q$-th order partial derivative operator, with $\widehat\nabla^0 v = v$. This definition yields well-defined Hilbert spaces endowed with the broken norm and seminorm
\begin{align*}
	\norm{\widehat v}^2_{\mathcal{X}^{s}_{\beta}(\widehat{E})} 
	:=\sum_{Q \in \widehat E} \norm{\widehat v}^2_{\widehat{X}^{s}_{\beta}(Q)} , \quad 
	\abs{\widehat v}^2_{\mathcal{X}^{s}_{\beta}(\widehat{E})} 
	:=\sum_{Q \in \widehat E} \abs{\widehat v}^2_{\widehat{X}^{s}_{\beta}(Q)} .
\end{align*}

The bent Sobolev spaces without weight are denoted by $\mathcal{H}^{s}(\widehat{E}) = \mathcal{H}^{s}_{0}(\widehat{E})$ and $\mathcal{V}^{s}(\widehat{E}) = \mathcal{V}^{s}_{0}(\widehat{E})$ and we set $\mathcal{L}^2_{\beta}(\widehat{E}) = \mathcal{H}^{0}_{\beta}(\widehat{E}) = \mathcal{V}^{0}_{\beta}(\widehat{E})$ and $\mathcal{L}^2(\widehat{E}) = \mathcal{H}^{0}(\widehat{E}) = \mathcal{V}^{0}(\widehat{E})$.


In the next step, we combine the concepts of polar and bent function spaces and introduce 
\textit{weighted bent polar Sobolev spaces} $\mathcal{X}^s_{\beta}(\widehat{E}) = \mathcal{H}^s_{\beta}(\widehat{E})$ and $\mathcal{X}^s_{\beta}(\widehat{E}) = \mathcal{V}^s_{\beta}(\widehat{E})$ of order $s \in \N_0$ with weight $\beta \in \R$, 
\begin{align*}
	\mathcal{X}^s_{\pol,\beta}(\widehat{E}) := \left \{
	\widehat v \in \mathcal{D}'(\widehat E) \colon
	\begin{rcases} \begin{dcases}
			\widehat v_{|Q} \in \widehat{X}^s_{\pol,\beta}(Q) & \text{for all } Q \in \widehat{E} \text{ and } \\
			\widehat \nabla^q( \widehat v_{|Q_{1}}) = \widehat \nabla^q (\widehat v_{|Q_2}) 
			& \text{on } \partial Q_1 \cap \partial Q_2 \text{ for } q=0,1,\dots, \min\{m_{Q_1, Q_2}, s-1 \} \\
			& \text{for all }  Q_1, Q_2 \text{ with } \partial Q_1 \cap \partial Q_2 \neq \emptyset
	\end{dcases} \end{rcases} \right\},
\end{align*}
where we insert $\widehat X^s_{\pol,\beta}(Q) = \widehat H^s_{\pol,\beta}(Q)$ and $\widehat X^s_{\pol,\beta}(Q) = \widehat V^s_{\pol,\beta}(Q)$, respectively, endowed with the broken norm and seminorm
\begin{align*}
	\norm{\widehat v}^2_{\mathcal{X}^{s}_{\pol,\beta}(\widehat{E})} 
	:=\sum_{Q \in \widehat E} \norm{\widehat v}^2_{\widehat{X}^{s}_{\pol,\beta}(Q)} , \quad 
	\abs{\widehat v}^2_{\mathcal{X}^{s}_{\pol,\beta}(\widehat{E})} 
	:=\sum_{Q \in \widehat E} \abs{\widehat v}^2_{\widehat{X}^{s}_{\pol,\beta}(Q)} .
\end{align*}
Bent polar Sobolev spaces without weight are denoted by $\mathcal{H}^s_{\pol}(\widehat{E}) = \mathcal{H}^s_{\pol,0}(\widehat{E})$ and $\mathcal{V}^s_{\pol}(\widehat{E}) = \mathcal{V}^s_{\pol,0}(\widehat{E})$ and we set $\mathcal{L}^2_{\pol,\beta}(\widehat{E}) = \mathcal{H}^{0}_{\pol,\beta}(\widehat{E}) = \mathcal{V}^{0}_{\pol,\beta}(\widehat{E})$ and $\mathcal{L}^2_{\pol}(\widehat{E})=  \mathcal{H}^{0}_{\pol}(\widehat{E}) = \mathcal{V}^{0}_{\pol}(\widehat{E})$.

From \eqref{eq: transformation of norm on single element adjacent to polar point}, we immediately obtain the relation
\begin{align}
	\label{eq: transformation polar bent Sobolev spaces}
	\norm{v}_{V^{s}_{\beta}(E)} \sim \norm{\widehat{v}}_{\mathcal{V}^{s}_{\pol,\beta}(\widehat{E})} 
	\quad \text{ for all } v \in V^{s}_{\beta}(E),
\end{align} 
with $s \in \N_0$ and $\beta \in \R$, and \eqref{eq: polar spaces and 1/2-spaces} yields the estimate
\begin{align}
	\label{eq: bent polar spaces and 1/2-spaces}
	\norm{\widehat{v}}_{\mathcal{H}^s_{1/2+\beta}(\widehat{E})} 
	\leq C \norm{\widehat{v}}_{\mathcal{H}_{\pol,\beta}^s(\widehat{E})} 
	\leq C \norm{\widehat{v}}_{\mathcal{V}_{\pol,\beta}^s(\widehat{E})}  
	\quad \text{ for all } \widehat v \in \mathcal{V}_{\pol,\beta}^s(\widehat{E}).
\end{align}
Thus, the embeddings $\mathcal{V}_{\pol,\beta}^s(\widehat{E}) \subset \mathcal{H}_{\pol,\beta}^s(\widehat{E}) \subset \mathcal{H}^s_{1/2+\beta}(\widehat{E})$ hold. 

Lastly, in this context, suitable projectors $\Pi: \mathcal{X}^{s}_{\pol,\beta} (\widehat{\Omega}) \to S_{(s-1,s-1)}(\bs \Xi)$ are required that map from weighted bent polar Sobolev spaces of order $s$ into the spline spaces of degree $s-1$ and satisfy the crucial property
\begin{align}
	\label{eq: projector for bent polar spaces}
	\widehat v - \Pi(\widehat v) \in \widehat{X}^{s}_{\pol,\beta}(\widehat{\Omega}) \quad \text{ for all } \widehat v \in \mathcal{X}^{s}_{\pol,\beta} (\widehat{\Omega}) .
\end{align}
The existence of such projectors can be shown by analogy to similar results from the literature \cite{BeiraodaVeigaBuffaSangalliVazquez2014}.

\subsubsection{Transformation of the projection error}
\label{subsec: transformation of projection error}
Finally, we perform a transformation of the projection error, which will be needed to prove Theorem \ref{theorem: main result}. Therefore, in addition to the computations above, the weight function has to be taken into account. For $q=0$, we obtain
\begin{align}
	\norm{v -  \Pi_{V_h^{\pol}} v }^2_{L^2(K)}
	&= \int_{K} \abs{v -  \Pi_{V_h^{\pol}} v} ^2 \Dd \bs x \nonumber \sim \int_{Q} \abs{v \circ \bs F -   \Pi_{V_h^{\pol}} v \circ \bs F} ^2 r \dr \dphi \nonumber \\
	&= \int_{Q} \abs{\frac{ W(v \circ \bs  F) - \Pi^{\pol}_{\bs p , \bs \Xi} (W( v \circ \bs F)) }{W}} ^2 r \dr \dphi \nonumber 
	\leq C \int_{Q} \abs{ W \widehat{v} - \Pi^{\pol}_{\bs p , \bs \Xi}( W \widehat{v})} ^2 r \dr \dphi \nonumber \\
	&= C \norm{W \widehat{v} - \Pi^{\pol}_{\bs p , \bs \Xi} (W \widehat{v})}_{\widehat L^2_{1/2}(Q)}^2 , \label{eq: transformation of projection error L^2}
\end{align}
where we use the boundedness of $W$ and its inverse $W^{-1}$, see \cite{BeiraodaVeigaBuffaSangalliVazquez2014}, and set $\widehat{v} =  v \circ \bs  F$. For $q=1$, we compute
\begin{align}
	\abs{v -  \Pi_{V_h^{\pol}} v }^2_{H^1(K)}
	&= \int_{K} \abs{\partial_x\left({v -  \Pi_{V_h^{\pol}} v}\right)}^2 + \abs{\partial_y \left({v -  \Pi_{V_h^{\pol}} v}\right)} ^2 \Dd \bs x \nonumber \\
	&\sim \int_{Q} \left(\abs{\partial_r\left({v \circ \bs F -  \Pi_{V_h^{\pol}} v \circ \bs F}\right)}^2 + \abs{\frac{1}{r} \partial_\varphi \left({v  \circ \bs F -  \Pi_{V_h^{\pol}} v \circ \bs F}\right)} ^2 \right) r \dr \dphi \nonumber \\
	&= \int_{Q} \left( \abs{\partial_r\left(\frac{ W \widehat{v} - \Pi^{\pol}_{\bs p, \bs \Xi}( W \widehat{v}) }{W}\right)}^2 + \abs{\frac{1}{r} \partial_\varphi \left(\frac{ W \widehat{v} - \Pi^{\pol}_{\bs p, \bs \Xi}( W \widehat{v} )}{W}\right)} ^2 \right) r \dr \dphi \nonumber \\
	& \leq C \bigg( \norm{\partial_r \left( W \widehat{v} - \Pi^{\pol}_{\bs p, \bs \Xi}( W \widehat{v}) \right)}_{\widehat L^2_{1/2}(Q)}^2 + \norm{W \widehat{v} -\Pi^{\pol}_{\bs p, \bs \Xi} (W \widehat{v}) }_{\widehat L^2_{1/2}(Q)}^2 \nonumber \\
	& \  \qquad + \norm{\partial_\varphi\left( W \widehat{v} -\Pi^{\pol}_{\bs p, \bs \Xi}( W \widehat{v}) \right)}_{\widehat L^2_{-1/2}(Q)}^2 + \norm{W \widehat{v} - \Pi^{\pol}_{\bs p, \bs \Xi} (W \widehat{v}) }_{\widehat L^2_{-1/2}(Q)}^2 \bigg) , \label{eq: transformation of projection error H^1}
\end{align}
where we use the Leibniz formula
\begin{align}
	\partial^{\bs \alpha} \left(\frac{g}{W}\right) = \sum_{\{\bs \beta: \bs \beta \leq \bs \alpha\}} \binom{\bs \alpha}{ \bs \beta} \left(\partial^{\bs \beta} g\right) \left(\partial^{\bs \alpha - \bs \beta} W^{-1}\right)
	\label{eq: Leibniz formula}
\end{align}
with $g=\widehat{v} - \Pi^{\pol}_{\bs p, \bs \Xi} \widehat{v}$ and Young's formula. Here, we remark again that the weight function $W$ and its inverse $W^{-1}$ as well as their derivatives are bounded \cite{BeiraodaVeigaBuffaSangalliVazquez2014}. Throughout the paper, the latter will be used multiple times, and the Leibniz formula will be applied in combination with the weight function without further explication.

Finally, let $E \subset \Omega$ be given by a union of elements $K \in \mathcal{M}$ and let $\widehat E = \bs F^{-1}(E)\subset \widehat{\Omega}$. Then, we have analogous results to \eqref{eq: transformation of projection error L^2} and \eqref{eq: transformation of projection error H^1}, but with bent Sobolev spaces,
\begin{align*}
	\norm{v -  \Pi_{V_h^{\pol}} v }^2_{L^2(E)}
	\leq \norm{W \widehat{v} - \Pi^{\pol}_{\bs p , \bs \Xi} (W \widehat{v})}_{\mathcal{L}^2_{1/2}(\widehat{E})}^2
\end{align*}
and 
\begin{align*}
	\abs{v -  \Pi_{V_h^{\pol}} v }^2_{H^1(E)}
	& \leq C \bigg( \norm{\partial_r \left( W \widehat{v} - \Pi^{\pol}_{\bs p, \bs \Xi}( W \widehat{v}) \right)}_{\mathcal L^2_{1/2}(\widehat{E})}^2 + \norm{W \widehat{v} -\Pi^{\pol}_{\bs p, \bs \Xi} (W \widehat{v}) }_{\mathcal L^2_{1/2}(\widehat{E})}^2 \nonumber \\
	& \  \qquad + \norm{\partial_\varphi\left( W \widehat{v} -\Pi^{\pol}_{\bs p, \bs \Xi}( W \widehat{v}) \right)}_{\mathcal L^2_{-1/2}(\widehat{E})}^2 + \norm{W \widehat{v} - \Pi^{\pol}_{\bs p, \bs \Xi} (W \widehat{v}) }_{\mathcal L^2_{-1/2}(\widehat{E})}^2 \bigg).
\end{align*}
Throughout the paper, the transformation of norms between the physical and parametric domains will sometimes be referred to as a change of variables.

\section{Proof of Theorem \ref{theorem: main result}}
\label{section: Error estimates}
In this section, we derive error estimates for the proposed modified projector with respect to graded $h$-refinement on polar domains with corners. To this end, the error is analyzed separately in the region surrounding the polar point and in the remainder of the domain, which is not affected by the corner singularity.

\subsection{Splitting of the polar domain}
\label{subsec: splitting of the model domain}
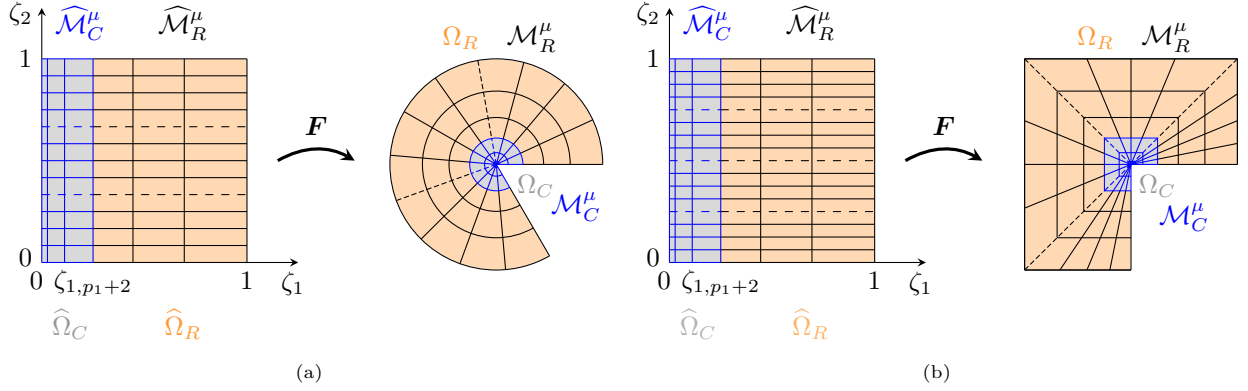
\begin{figure}
	\begin{center}
		\begin{subfigure}{0.495\textwidth}
			\begin{center}
				\input{splitting_circ_par.tikz}
				\hspace{1mm}
				\input{splitting_circ_phys.tikz}
				\caption{}
				\label{fig: partition of domain circular sector}
			\end{center}
		\end{subfigure}
		\begin{subfigure}{0.495\textwidth}
			\begin{center}
				\input{splitting_L_par.tikz}
				\hspace{2mm}
				\input{splitting_L_phys.tikz}
				\caption{}	
				\label{fig: partition of domain L-shape}
			\end{center}
		\end{subfigure}
		\caption{Splitting of the parametric and physical domain and the corresponding meshes. (a): Circular sector (b): L-shaped domain.}
		\label{fig: partition of domain}
	\end{center}
\end{figure}

To begin, we explain in more detail how the model domain and the underlying meshes can be split suitably. Let $\widehat{\mathcal{M}}^\mu$ and $\mathcal{M}^\mu$ be graded parametric and physical meshes, respectively, resulting from the knot vectors $\Xi_1^{h_1,\mu}$ and $\Xi_2^{h_2}$ for some grading parameter $\mu \in (0,1]$ and polynomial degrees $\bs p =(p_1,p_2)$, recall Section \ref{subsec: construction of graded mesh}. To describe a region surrounding the corner $\bs P$, we consider the physical B\'ezier elements $K_{\bs j} = \bs F(Q_{\bs j}) \in \mathcal{M}^\mu$ whose support extensions satisfy
\begin{align}
	\label{eq: characterization singular physical elements}
	\bs P \in \overline{\widetilde{K}_{\bs j}} = \overline{\bs F(\widetilde{Q}_{\bs j})} = \bs F(\overline{\widetilde{Q}_{\bs j}}) ,
\end{align}
that is, the closure of the support extension $\widetilde{K}_{\bs j}$ contains the singularity. The corresponding parametric mesh elements ${Q}_{\bs j}$ are characterized by the condition
\begin{align}
	\label{eq: characterization singular parametric elements}
	\Gamma_1 \cap \overline{\widetilde{Q}_{\bs j}} 
	= (\{0\} \times [0,1]) \cap \left( \overline{\widetilde{I}_{1,j_1}} \times \overline{\widetilde{I}_{2,j_2}} \right)
	=  \left(\{0\} \cap \overline{\widetilde{I}_{1,j_1}}\right) \times \overline{\widetilde{I}_{2,j_2}} \neq \emptyset ,
\end{align}
i.e., the singular edge $\Gamma_1$ touches the closure of the support extension $\widetilde{Q}_{\bs j}$, which is equivalent to $0 \in  \overline{\widetilde{I}_{1,j_1}}$, recall definition \eqref{eq: bivariate support extension}.
The refined knot vector in $\zeta_1$-direction can be expressed in terms of the breakpoints \eqref{eq: graded vector of breakpoints} by
\begin{align}
	\label{eq: knot vector p-refined}
	\quad \Xi_1^{h_1,\mu} = \{\xi_1, \xi_2, \dots, \xi_{n_1+p_1+1}\} = \left\{\underbrace{\zeta_{1,1},\dots,\zeta_{1,1}}_{p_1+1 \text{ times}}, \zeta_{2,1}, \zeta_{3,1}, \dots, \zeta_{N_1-2,1}, \zeta_{N_1-1,1}, \underbrace{\zeta_{N_1,1},\dots,\zeta_{N_1,1}}_{p_1+1 \text{ times}} \right\} .
\end{align}
Thus, it holds $n_1+p_1+1 = N_1 + 2p_1$ and
\begin{align*}
	{I}_{1,j_1}
	=\left(\zeta_{1,j_1}, \zeta_{1,j_1+1}\right) 
	=\left(\xi_{1,j_1+p_1}, \xi_{1,j_1+p_1+1}\right) , \quad j_1 = 1,2,\dots, N_1-1 ,
\end{align*}
and with definition \eqref{eq: univariate support extension}, we obtain
\begin{align}
	\label{eq: support extension graded direction}
	\widetilde{I}_{1,j_1} 
	= \left(\xi_{1,j_1+p_1-p_1}, \xi_{1,j_1+p_1+p_1+1}\right) 
	= \left(\xi_{1,j_1}, \xi_{1,j_1+2 p_1+1}\right) , \quad j_1 = 1,2,\dots, N_1-1 .
\end{align}
By construction \eqref{eq: knot vector p-refined}, it holds $\xi_{1,j_1} = 0$ and hence $0 \in \overline{\widetilde{I}_{1,j_1}}$ if and only if $j_1\leq p_1+1$. Therefore, we divide the graded parametric and physical B\'ezier mesh $\widehat{\mathcal{M}}^\mu$ and \ $\mathcal{M}^\mu$ into two submeshes,
\begin{align*}
	\widehat{\mathcal{M}}^\mu = \widehat{\mathcal{M}}^\mu_C \cup \widehat{\mathcal{M}}^\mu_R  \quad \text{ and } \quad
	\mathcal{M}^\mu = \mathcal{M}^\mu_C \cup \mathcal{M}^\mu_R  ,
\end{align*} 
respectively, with
\begin{alignat*}{2}
	\widehat{\mathcal{M}}^\mu_C &:= \{Q_{\bs j}\}_{j_1\in \{1,2, \dots, p_1+1\}, j_2\in \{1,2, \dots, N_2\}} , \quad
	\widehat{\mathcal{M}}^\mu_R &&:= \{Q_{\bs j}\}_{j_1\in \{p_1+2,p_1+3, \dots, N_1\}, j_2\in \{1,2, \dots, N_2\}} ,\\
	\mathcal{M}^\mu_C &:= \{K_{\bs j}\}_{j_1\in \{1,2, \dots, p_1+1\}, j_2\in \{1,2, \dots, N_2\}} , \quad
	\mathcal{M}^\mu_R &&:= \{K_{\bs j}\}_{j_1\in \{p_1+2,p_1+3, \dots, N_1\}, j_2\in \{1,2, \dots, N_2\}} .
\end{alignat*}
Induced by these partitions, we also divide the parametric domain into two subdomains, 
\begin{align}
	\widehat{\Omega} = \widehat{\Omega}_C \cup \widehat{\Omega}_R ,
	\text{ with }
	\widehat{\Omega}_C:= (0, \zeta_{1,p_1+2}] \times (0,1) \text{ and }
	\widehat{\Omega}_R:= [\zeta_{1,p_1+2},1) \times (0,1) .
	\label{eq: definition of parametric strips}
\end{align}
In this way, we obtain the desired partition of the physical domain,
\begin{align*}
	\Omega = \Omega_C \cup \Omega_R,
\end{align*}
into a small polar domain surrounding the corner,
\begin{align*}
	\Omega_C := \bs F(\widehat{\Omega}_C) =  \bs G(\widehat{\Omega}_C) = \left\{r R(\omega\varphi) (\cos (\omega\varphi) , \sin(\omega\varphi) )^T \in \R^2: 0 < r \leq \zeta_{1,p_1+2}, 0 < \varphi < \omega\right\} ,
\end{align*}
and the remaining ring-type part
\begin{align*}
	\Omega_R := \bs F(\widehat{\Omega}_R) = \left\{r R(\omega\varphi) (\cos (\omega\varphi) , \sin(\omega\varphi) )^T \in \R^2: \zeta_{1,p_1+2} \leq r < 1, 0 < \varphi < \omega\right\} ,
\end{align*}
where we used the reference parameterization in polar coordinates \eqref{eq: reference transformation}. In Figure \ref{fig: partition of domain}, we provide an illustration of the partitioned domains and meshes for a circular sector and an L-shaped domain. Note that the splitting depends on the mesh size $h_1$ and the polynomial degree $p_1$ of the first univariate direction, which is not pointed out explicitly in the notation. Lastly, we define the overlapping domains consisting of the union of the support extensions of all elements in $\mathcal{M}^\mu_C$ and $\mathcal{M}^\mu_R$, respectively,
\begin{align}
	\label{eq: overlapping splitting}
	\widetilde{\widehat{\Omega}}_{C}:=\bigcup_{Q_{\bs j} \in \widehat{\mathcal{M}}_C^\mu} \widetilde{Q}_{\bs j} , \quad
	\widetilde{\widehat{\Omega}}_{R}:=\bigcup_{Q_{\bs j} \in \widehat{\mathcal{M}}_R^\mu} \widetilde{Q}_{\bs j} , \quad
	\widetilde{\Omega}_{C}:=\bigcup_{K_{\bs j} \in \mathcal{M}_C^\mu} \widetilde{K}_{\bs j} , \quad
	\widetilde{\Omega}_{R}:=\bigcup_{K_{\bs j} \in \mathcal{M}_R^\mu} \widetilde{K}_{\bs j} .
\end{align}

To demonstrate the purpose of the constructed domain splitting, we show as a first simple consequence that the modified projector $\Pi_{V_h^{\pol}}$ coincides with the standard projector $\Pi_{V_h}$ on the ring-type part $\Omega_R$ of the domain, where the projection is not affected by the singularity.

\begin{lemma}
	\label{lemma: coinciding projectors}
	For all $v \in C(\overline{\Omega})$ and $\widehat{w} \in C([0,1]^2)$, it holds
	\begin{align*}
		\Pi_{V_h^{\pol}} v \Big| _{\Omega_R} =	\Pi_{V_h} v \Big| _{\Omega_R} 
		\quad \text{ and } \quad
		\Pi^{\pol}_{\bs p , \bs \Xi} (\widehat{w}) \Big| _{\widehat{\Omega}_R}  
		= \Pi_{\bs p , \bs \Xi} (\widehat{w}) \Big| _{\widehat{\Omega}_R} .
	\end{align*}	
\end{lemma}
\begin{proof}
	By construction of $\widehat{\Omega}_C$, the supports of the singular basis functions \eqref{eq: singular basis functions} satisfy 
	\begin{align*}
		\text{supp}(\widehat{B}_{\bs i, \bs p}) \subset \widehat{\Omega}_C \quad \text{ for all } \bs i_{\pol} \in \bs I_{\pol} = \{(1,i_2):1\leq i_2 \leq n_2\} 
	\end{align*}
	The modified quasi-interpolant \eqref{eq: modified quasi-interpolant} is thus given on $\widehat{\Omega}_R$ by
	\begin{align*}
		\Pi^{\pol}_{\bs p , \bs \Xi} (\widehat{w}) \Big|_{\widehat{\Omega}_R}
		= \sum_{\bs i \in \bs I \setminus \bs I_{\pol}} \lambda^{\pol}_{\bs i, \bs p}(\widehat{w}) \widehat{B}_{\bs i,\bs p} \Big| _{\widehat{\Omega}_R}
		= \sum_{\bs i \in \bs I} \lambda_{\bs i, \bs p}(\widehat{w}) \widehat{B}_{\bs i,\bs p} \Big| _{\widehat{\Omega}_R} 
		= \Pi_{\bs p , \bs \Xi} (\widehat{w}) \Big|_{\widehat{\Omega}_R} .
	\end{align*}
	The assertion for the projectors $\Pi_{V_h^{\pol}}$ and $\Pi_{V_h}$ then follows by using their definitions \eqref{eq: modified projector physical domain} and \eqref{eq: standard projector physical domain}, respectively, and setting $\widehat{w} = W(v \circ F)$ for $v \in C(\overline{\Omega})$.
\end{proof}

\subsection{Projection error estimates}
\label{subsec: projection error estimates}
Next, we proceed with the projection error estimates for the two regions of the polar domain.

\subsubsection{Estimates on $\Omega_R$}
We begin by deriving estimates on $\Omega_R$, where the distance to the singularity is strictly positive. In this setting, a similar strategy to that used in standard anisotropic approximation theory can be applied \cite{BeiraodaVeigaChoSangalli2012,BeiraodaVeigaBuffaSangalliVazquez2014}, in which element-wise estimates on the parametric domain are mapped to the physical domain via a change of variables. However, in contrast to the literature, the transformed quantities depend on the distance between each element and the singularity, and thus, by construction of $\Omega_R$, also on the mesh size.

The following element-wise error estimate on the parametric domain from the literature serves as the starting point for our analysis.
\begin{lemma}
	\label{lemma: standard projection error parametric domain}
	Let the integers $0\leq r_1 \leq s_1 \leq p_1 + 1$ and $ 0 \leq r_2 \leq s_2 \leq p_2 +1$. Then, there exists a constant $C$ depending only on $\bs p, \theta$ such that for all elements $Q_{\bs j} \in \widehat{\mathcal{M}}^\mu$, we have
	\begin{align*}
		\norm{\widehat{D}^{(r_1,r_2)} (\widehat v - \Pi_{\bs p,\bs \Xi} \widehat v)}_{\widehat{L}^2(Q_{\bs j})}
		\leq C \left(\widetilde{h}_{1,j_1}^{s_1-r_1} \norm{\widehat{D}^{(s_1,r_2)} \widehat v }_{\mathcal{L}^2(\widetilde{Q}_{\bs j})} + \widetilde{h}_{2,j_2}^{s_2-r_2} \norm{\widehat{D}^{(r_1,s_2)} \widehat v }_{\mathcal{L}^2(\widetilde{Q}_{\bs j})} \right)
	\end{align*}
	for all functions $\widehat v: \widehat{\Omega} \to \R$ with $\widehat{D}^{(s_1,r_2)} \widehat v$, $\widehat{D}^{(r_1,s_2)} \widehat v \in \mathcal{L}^2(\widetilde{\Omega}_{R})$, where $\widetilde{\Omega}_{R}$ is defined as in \eqref{eq: overlapping splitting}.
\end{lemma}

\begin{proof}
	The result is proven in \cite[Proposition 4.1]{BeiraodaVeigaChoSangalli2012}.
\end{proof}

Following the framework in \cite{BeiraodaVeigaChoSangalli2012,BeiraodaVeigaBuffaSangalliVazquez2014}, we define the coordinate system naturally induced by $\bs F$, referred to as the $\bs F$-coordinate system, whose tangent base vectors $\bs g_1$ and $\bs g_2$ can be defined by
\begin{align*}
	\bs g_i = \bs g_i(\bs x) = \frac{\partial \bs F}{\partial \zeta_i}(\bs F^{-1}(\bs x)) , \quad i =1,2 . 
\end{align*}
We also consider the derivatives of functions $v: \Omega \to \R$ with respect to the $\bs F$-coordinates,
\begin{align*}
	&\frac{\partial v}{\partial \bs g_i}(\bs x) 
	= \nabla v(\bs x) \cdot \bs g_i(\bs x)
	= \lim\limits_{t\to 0} \frac{v(\bs x+t\bs g_i(\bs x)) - v(\bs x)}{t} , \quad i =1,2 , \\
	&\frac{\partial^{\alpha_i} v}{\partial \bs g_i^{\alpha_i}} 
	= \frac{\partial}{\partial \bs g_i} \left( \frac{\partial^{\alpha_i-1}v}{\partial \bs g_i^{\alpha_i-1}} \right) 
	= \frac{\partial}{\partial \bs g_i} \left( \cdots \left( \frac{\partial}{\partial \bs g_i} \left( \frac{\partial v}{\partial \bs g_i} \right) \right) \right) , \quad \alpha_i \in \N , i =1,2 , \\
	& D_{\bs F}^{\bs \alpha} v = \frac{\partial^{\alpha_1}}{\partial \bs g_1^{\alpha_1}} \frac{\partial^{\alpha_2} v}{\partial \bs g_2^{\alpha_2}} , \quad \bs \alpha=(\alpha_1,\alpha_2) \in \N^2 .
\end{align*}
In our setting, the tangential basis of the $\bs F$-coordinate system for all $\bs x \neq (0,0)^T$ is given by
\begin{align*}
	\bs g_1(\bs x) \sim \begin{pmatrix} \cos \varphi (\bs x) \\ \sin \varphi (\bs x) \end{pmatrix}
	\quad \text { and } \quad 
	\bs g_2(\bs x) \sim \begin{pmatrix} -r \sin \varphi (\bs x) \\ r \cos \varphi (\bs x) \end{pmatrix}, 
\end{align*}
where the equivalence between polar and isogeometric coordinates is exploited, recall \eqref{eq: equivalence relation isogeometric and polar coordinates}. We exemplarily compute the derivatives $D_{\bs F}^{(1,0)} v$ and $\mathrm{D}_{\bs F}^{(0,1)} v$ in more detail. By applying standard calculus, it follows
\begin{align}
	\mathrm{D}_{\bs F}^{(1,0)} v (\bs x)
	&= \frac{\partial v}{\partial \bs g_1}(\bs x)
	= \nabla v(\bs x) \cdot \bs g_1(\bs x)
	\sim \begin{pmatrix} \frac{\partial v}{\partial x} (\bs x) \\ \frac{\partial v}{\partial y} (\bs x) \end{pmatrix} \cdot \begin{pmatrix} \cos \varphi(\bs x) \\ \sin \varphi(\bs x) \end{pmatrix}
	=  \cos \varphi(\bs x) \frac{\partial v}{\partial x} (\bs x) + \sin \varphi(\bs x) \frac{\partial v}{\partial y}(\bs x) \label{eq: derivative wrt r} \\
	&= \frac{\partial \widehat v}{\partial r} (r(\bs x), \varphi(\bs x)) 
	= \frac{\partial \widehat v}{\partial r} (\bs F^{-1}(\bs x)) \nonumber \\
	\mathrm{D}_{\bs F}^{(0,1)} v (\bs x)
	&= \frac{\partial v}{\partial \bs g_2} (\bs x)
	= \nabla v(\bs x) \cdot \bs g_2(\bs x)
	\sim -r(\bs x) \sin \varphi(\bs x) \frac{\partial v}{\partial x}(\bs x) + r(\bs x) \cos \varphi(\bs x) \frac{\partial v}{\partial y}(\bs x)  \label{eq: derivative wrt phi}\\
	&= \frac{\partial \widehat v}{\partial \varphi} (r(\bs x), \varphi(\bs x)) 
	= \frac{\partial \widehat v}{\partial \varphi} (\bs F^{-1}(\bs x)) \nonumber
\end{align}
In fact, as shown in \cite[Proposition 5.1]{BeiraodaVeigaChoSangalli2012}, we have
\begin{align}
	\mathrm{D}_{\bs F}^{\bs \alpha} v = (\widehat{D}^{\bs \alpha} \widehat v) \circ \bs F^{-1}.
	\label{eq: relation parametric derivative and F-derivative}
\end{align}
Let $E \subset \Omega$ be a union of elements $K \in \mathcal{M}^{\mu}$ and let $v \in H^s(E)$. Then, due to the reduced regularity of the parameterization $\bs F$ along mesh lines corresponding to repeated knot vectors, recall Section \ref{subsec: bent spaces}, the derivatives \eqref{eq: relation parametric derivative and F-derivative} are not necessarily in $L^2(E)$, but surely in $L^2(K)$ for all $K \in E$. Therefore, in the spirit of bent Sobolev spaces, we introduce the norms
\begin{align*}
	\norm{\mathrm{D}_{\bs F}^{\bs \alpha} v}_{\mathcal{L}^2(E)}^2 = \sum_{K \in E} \norm{\mathrm{D}_{\bs F}^{\bs \alpha} v}_{L^2(K)}.
\end{align*}
Finally, due to \eqref{eq: derivative wrt r}, \eqref{eq: derivative wrt phi} and the splitting of our domain, we have for every element $K \in \mathcal{M}_R^{\mu}$ 
\begin{align*}
	\norm{\mathrm{D}_{\bs F}^{(1,0)} v}_{\mathcal{L}^2(\widetilde{K})} \leq C \left( \norm{\frac{\partial v}{\partial x}}_{L^2(\widetilde{K})} + \norm{\frac{\partial v}{\partial y}}_{L^2(\widetilde{K})} \right) , \\
	\norm{r^{-1}\mathrm{D}_{\bs F}^{(0,1)} v}_{\mathcal{L}^2(\widetilde{K})} \leq C \left( \norm{\frac{\partial v}{\partial x}}_{L^2(\widetilde{K})} + \norm{\frac{\partial v}{\partial y}}_{L^2(\widetilde{K})} \right),
\end{align*}
and, more generally, it holds
\begin{align}
	\norm{r^{-\alpha_2}\mathrm{D}_{\bs F}^{\bs \alpha} v}_{\mathcal{L}^2(\widetilde{K})} \leq C 
	\sum_{\abs{\bs \beta} \leq \abs{\bs \alpha}} \norm{D^{\bs \beta} v}_{L^2(\widetilde{K})}.
	\label{eq: relation physical derivative and F-derivative}
\end{align}

%
%
%

With these things at hand, we can develop an estimate on the full ring-type domain part $\Omega_{R}$. For simplicity, we set $p = \min\{p_1,p_2\}$ from now on.

\begin{Theorem}
	\label{theorem: projection error ring domain}
	Let $q \in \{0,1\}$ and $s,s_0 \in \N$ with $2 \leq s \leq s_0 \leq p +1$. Further, let $v \in V^{s}_{\beta}(\widetilde{\Omega}_R)$ for all $\beta >s-1-\nu$, recall \eqref{eq: weighted regularity of singular part V spaces}. If the mesh grading parameter satisfies condition \eqref{eq: mesh grading parameter condition}, then there exists a constant $C$ depending only on $p,\bs \theta, \bs F, W$ such that
	\begin{align*}
		\abs{v - \Pi_{V_h^{\pol}} v}_{H^q(\Omega_R)} 
		\leq C h^{s-q} \norm{v}_{V^{s}_{\beta}(\widetilde{\Omega}_R)}   ,
	\end{align*}
	with $\widetilde{\Omega}_R$ as defined in \eqref{eq: overlapping splitting}.
\end{Theorem}

\begin{proof}
	Let $K_{\bs j} \in \mathcal{M}_R^\mu$ and $Q_{\bs j} = \bs F ^{-1}(K_{\bs j}) \in \widehat{\mathcal{M}}_R^\mu$. By construction of the splitting defined in Section \ref{subsec: splitting of the model domain}, the element $K_{\bs j}$ is not adjacent to the polar point, $\bs P \not \in \overline{K_{\bs j}}$, and thus the Jacobian of $\bs F$ is bounded on $\overline{Q_{\bs j}}$. However, the bounds depend on the distance of $K_{\bs j}$ to the polar point, which is equivalent to the distance of $Q_{\bs j}$ to the singular edge $\widehat{\Gamma_1}$. In more detail, let $Q_{\bs j}=(\zeta_{1,j_1} , \zeta_{1,j_1 + 1}) \times (\zeta_{2,j_2} , \zeta_{2,j_2 + 1}) \in \mathcal{M}^{\mu}_R$. Then, due to \eqref{eq: Jacobian of parameterization}, we have 
	\begin{align*}
		\zeta_{1,j_1} \leq \det(J_{\bs F}) \leq \zeta_{1,j_1 + 1} \quad \text{ and } \quad 
		\zeta_{1,j_1 + 1}^{-1} \leq \det(J_{\bs F^{-1}}) \leq \zeta_{1,j_1}^{-1} \quad \text{ in } Q_{\bs j}.
	\end{align*}
	As the considered graded meshes are locally quasi-uniform, see Lemma \ref{lemma: local quasi-uniformity of the graded mesh}, it holds $\zeta_{1,j_1 + 1} / \zeta_{1,j_1} \sim 1$.
	In addition, there are even constants $C_1,C_2 \geq 0$ such that
	\begin{align}
		C_1 \, \zeta_{1,j_1} \leq \det(J_{\bs F}) \leq C_2 \, \zeta_{1,j_1} \quad \text{ and } \quad 
		C_2^{-1} \, \zeta_{1,j_1}^{-1} \leq \det(J_{\bs F^{-1}}) \leq C_1^{-1} \, \zeta_{1,j_1}^{-1} \quad \text{ in } \widetilde{Q}_{\bs j}.
		\label{eq: bounded determinant on ring elements}
	\end{align}
	With \eqref{eq: transformation of projection error L^2} and Lemma \ref{lemma: coinciding projectors}, we thus obtain
	\begin{align}
		\norm{v - \Pi_{V_h}^{\pol} v}_{L^2(K_{\bs j})}
		\leq C \norm{W \widehat{v} - \Pi^{\pol}_{\bs p , \bs \Xi} (W \widehat{v})}_{\widehat L^2_{1/2}(Q_{\bs j})}
		\leq C \zeta_{1,j_1}^{1/2} \norm{W \widehat{v} - \Pi_{\bs p , \bs \Xi} (W \widehat{v})}_{\widehat L^2(Q_{\bs j})}
		\label{eq: proof projection error estimate}
	\end{align}
	and with \eqref{eq: transformation of projection error H^1}, it follows
	\begin{align}
		\abs{v - \Pi_{V_h} v}^2_{H^1(K_{\bs j})}
		&\leq C \bigg( \zeta_{1,j_1}^{1/2} \norm{\partial_r \left( W \widehat{v} - \Pi_{\bs p, \bs \Xi}( W \widehat{v}) \right)}_{\widehat L^2(Q_{\bs j})}^2 + \zeta_{1,j_1}^{1/2} \norm{W \widehat{v} -\Pi_{\bs p, \bs \Xi} (W \widehat{v}) }_{\widehat L^2(Q_{\bs j})}^2 \label{eq: proof of ring domain terms positive weight}\\
		& \  \qquad + \zeta_{1,j_1}^{-1/2} \norm{\partial_\varphi\left( W \widehat{v} -\Pi_{\bs p, \bs \Xi}( W \widehat{v}) \right)}_{\widehat L^2(Q_{\bs j})}^2 + \zeta_{1,j_1}^{-1/2} \norm{W \widehat{v} - \Pi_{\bs p, \bs \Xi} (W \widehat{v}) }_{\widehat L^2(Q_{\bs j})}^2 \bigg) \label{eq: proof of ring domain terms negative weight}. 
	\end{align}
	
	First, we estimate the terms in line \eqref{eq: proof of ring domain terms positive weight}. Therefore, we use Lemma \ref{lemma: standard projection error parametric domain} with $s_1=s$, $s_2=s-q$, $r_1=q$ and $r_2=0$ and apply the Leibniz formula, where we note that the weight function $W$ and its inverse $W^{-1}$ as well as their derivatives are bounded \cite{BeiraodaVeigaBuffaSangalliVazquez2014}. 
	In more detail, we recall that the weight function has the form \eqref{eq: weight function calculated} and consequently, all derivatives $\widehat{D}^{\bs \alpha}W$ for $\bs \alpha \geq (1,0)$ vanish, where we use the notation $\bs \alpha \geq \bs \gamma$ if and only if $\alpha_1 \geq \gamma_1$ and $\alpha_2 \geq \gamma_2$.
	It follows for $q \in \{0,1\}$ that
	\begin{align*}
		\norm{\partial_r^q \left( W \widehat{v} - \Pi_{\bs p, \bs \Xi}( W \widehat{v}) \right)}_{\widehat L^2(Q_{\bs j})}^2
		&\leq C \left(\widetilde{h}_{1,j_1}^{s-q} \norm{\widehat{D}^{(s,0)}(W\widehat v) }_{\mathcal{L}^2(\widetilde{Q}_{\bs j})} + \widetilde{h}_{2,j_2}^{s-q} \norm{\widehat{D}^{(q,s-q)} (W\widehat v) }_{\mathcal{L}^2(\widetilde{Q}_{\bs j})} \right)\\
		&\leq C \left(\widetilde{h}_{1,j_1}^{s-q} \sum_{\substack{\bs \alpha \leq (s,0) \\ \bs \alpha \geq (s,0)}} \norm{\widehat{D}^{\bs \alpha}\widehat v }_{\mathcal{L}^2(\widetilde{Q}_{\bs j})} + \widetilde{h}_{2,j_2}^{s-q} \sum_{\substack{\bs \alpha \leq (q,s-q) \\ \bs \alpha \geq (q,0)}}  \norm{ \widehat{D}^{\bs \alpha}\widehat v}_{\mathcal{L}^2(\widetilde{Q}_{\bs j})} \right) \\ 
		&= \left(\widetilde{h}_{1,j_1}^{s-q}  \norm{\widehat{D}^{(s,0)}\widehat v }_{\mathcal{L}^2(\widetilde{Q}_{\bs j})} + \widetilde{h}_{2,j_2}^{s-q} \sum_{\alpha_2 = 0}^{s-q}  \norm{ \widehat{D}^{(q,\alpha_2)}\widehat v}_{\mathcal{L}^2(\widetilde{Q}_{\bs j})} \right).
	\end{align*}
	Now, we perform a change of variables back to the physical domain using \eqref{eq: relation parametric derivative and F-derivative}, \eqref{eq: bounded determinant on ring elements} and \eqref{eq: relation physical derivative and F-derivative} and insert the mesh properties \eqref{eq: edge lengths parametric support extension} with $h = \max\{h_1,h_2\}$,
	\begin{align}
		& \quad \zeta_{1,j_1}^{1/2} \norm{\partial_r^q \left( W \widehat{v} - \Pi_{\bs p, \bs \Xi}( W \widehat{v}) \right)}_{\widehat L^2(Q_{\bs j})}^2  \nonumber \\
		&\leq C\left(\frac{\zeta_{1,j_1}}{\zeta_{1,j_1}}\right)^{1/2} \left(\widetilde{h}_{1,j_1}^{s-q} \norm{D^{(s,0)}_{\bs F}v }_{\mathcal{L}^2(\widetilde{K}_{\bs j})} + \widetilde{h}_{2,j_2}^{s-q} \sum_{\alpha_2 = 0}^{s-q} \norm{r^{\alpha_2}r^{-\alpha_2} D^{(q,\alpha_2)}_{\bs F}v}_{\mathcal{L}^2(\widetilde{K}_{\bs j})} \right) \nonumber\\
		&\leq C \left(\widetilde{h}_{1,j_1}^{s-q} \sum_{\abs{\bs \gamma} \leq s} \norm{D^{\bs \gamma}v }_{\mathcal{L}^2(\widetilde{K}_{\bs j})} + \widetilde{h}_{2,j_2}^{s-q} \sum_{\alpha_2 = 0}^{s-q} \sum_{\abs{\gamma}\leq q+\alpha_2}  \norm{r^{\alpha_2} D^{\bs \gamma}v}_{\mathcal{L}^2(\widetilde{K}_{\bs j})} \right) \nonumber\\
		&\leq C\, h^{s-q} \left( \left(\zeta_{1,j_1}\right)^{(1-\mu)(s-q)} \sum_{\abs{\bs \gamma} \leq s} \norm{D^{\bs \gamma}v }_{\mathcal{L}^2(\widetilde{K}_{\bs j})} + \sum_{\alpha_2 = 0}^{s-q} \sum_{\abs{\bs \gamma}\leq q+\alpha_2} \norm{r^{\alpha_2} D^{\bs \gamma}v}_{\mathcal{L}^2(\widetilde{K}_{\bs j})} \right).
		\label{eq: proof ring domain transform back to physical domain}
	\end{align}
	Due to the local quasi-uniformity of the mesh, recall Lemma \ref{lemma: local quasi-uniformity of the graded mesh}, and the construction of polar parameterizations, we have
	\begin{align}
		\zeta_{1,j_1} \sim \zeta_1 \sim  r (\bs x) \quad \text{ for all } \bs \zeta = (\zeta_1,\zeta_2) \in \widetilde{Q}_{\bs j}
		\text{ and } \bs x \in \widetilde{K}_{\bs j}.
		\label{eq: proof of ring domain pull zeta into norm}
	\end{align}
	Moreover, we set $\beta =(1 - \mu)(s-q)$. Then, as $\mu \in (0,1]$ satisfies condition \eqref{eq: mesh grading parameter condition}, it holds
	\begin{align}
		\beta \leq s -q \quad \text{ and } \quad \beta > \left (1 - \frac{\nu-q+1}{s-q} \right)(s-q) = s-q -(\nu-q+1) = s-1-\nu. \label{eq: beta property}
	\end{align}
	By combining \eqref{eq: proof ring domain transform back to physical domain}, \eqref{eq: proof of ring domain pull zeta into norm} and \eqref{eq: beta property}, we obtain
	\begin{align*}
		&\quad \zeta_{1,j_1}^{1/2} \norm{\partial_r^q \left( W \widehat{v} - \Pi_{\bs p, \bs \Xi}( W \widehat{v}) \right)}_{\widehat L^2(Q_{\bs j})}^2 \\ 
		&\leq C \, h^{s-q} \left(\sum_{\abs{\bs \gamma} \leq s} \norm{r^{\beta} D^{\bs \gamma}v}_{L^2(\widetilde{K}_{\bs j})} + \sum_{\alpha_2 = 0}^{s-q} \sum_{\abs{\bs \gamma}\leq q+\alpha_2} \norm{r^{\beta - (s-q) + \alpha_2} D^{\bs \gamma}v}_{L^2(\widetilde{K}_{\bs j})} \right) \\
		&\leq C \, h^{s-q} \left(\sum_{\abs{\bs \gamma} \leq s} \norm{r^{\beta} D^{\bs \gamma}v}_{L^2(\widetilde{K}_{\bs j})} + \sum_{\alpha_2 = 0}^{s-q} \sum_{\abs{\bs \gamma}\leq q+\alpha_2} \norm{r^{\beta - s + \abs{\bs \gamma}} D^{\bs \gamma}v}_{L^2(\widetilde{K}_{\bs j})} \right) \\
		&\leq C \, h^{s-q} \left(\sum_{\abs{\bs \gamma} \leq s} \norm{r^{\beta} D^{\bs \gamma}v}_{L^2(\widetilde{K}_{\bs j})} + \sum_{\abs{\bs \gamma}\leq s} \norm{r^{\beta - s + \abs{\bs \gamma}} D^{\bs \gamma}v}_{L^2(\widetilde{K}_{\bs j})} \right) \\
		&\leq C \, h^{s-q} \left(\norm{v}_{H^s_{\beta}(\widetilde{K}_{\bs j})} + \norm{v}_{V^{s}_{\beta}(\widetilde{K}_{\bs j})} \right) \leq C \, h^{s-q} \norm{v}_{V^{s}_{\beta}(\widetilde{K}_{\bs j})}.
	\end{align*}
	Next, we estimate the first term in line \eqref{eq: proof of ring domain terms negative weight} by applying Lemma \ref{lemma: standard projection error parametric domain} for $s_1=s-1$, $s_2=s$, $r_1=0$ and $r_2=1$ as well as the Leibniz formula,
	\begin{align*}
		\norm{\partial_\varphi \left( W \widehat{v} -\Pi_{\bs p, \bs \Xi}( W \widehat{v}) \right)}_{\widehat L^2(Q)}^2  
		&\leq C \left(\widetilde{h}_{1,j_1}^{s-1} \norm{\widehat{D}^{(s-1,1)}(W\widehat v) }_{\mathcal{L}^2(\widetilde{Q}_{\bs j})} + \widetilde{h}_{2,j_2}^{s-1} \norm{ \widehat{D}^{(0,s)}( W\widehat v) }_{\mathcal{L}^2(\widetilde{Q}_{\bs j})} \right) \\
		&\leq C \left( \widetilde{h}_{1,j_1}^{s-1} \sum_{\alpha_2 =0}^{1} \norm{\widehat{D}^{(s-1,\alpha_2)}\widehat v}_{\mathcal{L}^2(\widetilde{Q}_{\bs j})} + \widetilde{h}_{2,j_2}^{s-1} \sum_{\alpha_2=0}^{s}  \norm{ \widehat{D}^{(0, \alpha_2)}\widehat v}_{\mathcal{L}^2(\widetilde{Q}_{\bs j})} \right) .
	\end{align*}		
	Then, we perform a change of variables, see \eqref{eq: bounded determinant on ring elements} and \eqref{eq: relation physical derivative and F-derivative}, and exploit the mesh properties \eqref{eq: edge lengths parametric support extension} and \eqref{eq: proof of ring domain pull zeta into norm},
	\begin{align*}
		& \quad \zeta_{1,j_1}^{-1/2} \norm{\partial_\varphi \left( W \widehat{v} -\Pi_{\bs p, \bs \Xi}( W \widehat{v}) \right)}_{\widehat L^2(Q)}^2   \\
		&\leq C \, (\zeta_{1,j_1} \zeta_{1,j_1})^{-1/2} \left( \widetilde{h}_{1,j_1}^{s-1} \sum_{\alpha_2 =0}^{1} \norm{r^{\alpha_2} r^{-\alpha_2} D^{(s-1,\alpha_2)}_{\bs F}v }_{\mathcal{L}^2(\widetilde{K}_{\bs j})} + \widetilde{h}_{2,j_2}^{s-1}   \sum_{\alpha_2=0}^{s} \norm{r^{\alpha_2} r^{-\alpha_2} D^{(0,\alpha_2)}_{\bs F}v}_{\mathcal{L}^2(\widetilde{K}_{\bs j})} \right) \\
		&\leq C \, \zeta_{1,j_1}^{-1} \left( \widetilde{h}_{1,j_1}^{s-1} \sum_{\alpha_2 =0}^{1} \sum_{\abs{\bs \gamma} \leq s-1+\alpha_2} \norm{r^{\alpha_2} D^{\bs \gamma}v }_{\mathcal{L}^2(\widetilde{K}_{\bs j})} + \widetilde{h}_{2,j_2}^{s-1} \sum_{\alpha_2=0}^{s} \sum_{\abs{\bs \gamma} \leq \alpha_2} \norm{r^{\alpha_2} D^{\bs \gamma}v}_{\mathcal{L}^2(\widetilde{K}_{\bs j})} \right) \\
		&\leq C \, h^{s-1} \left( \left(\zeta_{1,j_1}\right)^{(1-\mu)(s-1)} \sum_{\alpha_2 =0}^{1} \sum_{\abs{\bs \gamma} \leq s+\alpha_2-1}\norm{r^{\alpha_2 -1} D^{\bs \gamma}v}_{L^2(\widetilde{K}_{\bs j})} + \sum_{\alpha_2=0}^{s} \sum_{\abs{\bs \gamma} \leq \alpha_2}  \norm{r^{\alpha_2-1} D^{\bs \gamma}v}_{L^2(\widetilde{K}_{\bs j})} \right) .
	\end{align*}
	We set again $\beta =(1 - \mu)(s-q)$ for $q=1$ and recall that \eqref{eq: beta property} is satisfied. It follows 
	\begin{align*}
		& \quad \zeta_{1,j_1}^{-1/2} \norm{\partial_\varphi \left( W \widehat{v} -\Pi_{\bs p, \bs \Xi}( W \widehat{v}) \right)}_{\widehat L^2(Q)}^2 \\
		&\leq C \, h^{s-1} \left(\sum_{\alpha_2 =0}^{1} \sum_{\abs{\bs \gamma} \leq s+\alpha_2-1} \norm{r^{\beta + \alpha_2 -1} D^{\bs \gamma}v}_{L^2(\widetilde{K}_{\bs j})} + \sum_{\alpha_2=0}^{s} \sum_{\abs{\bs \gamma} \leq \alpha_2} \norm{r^{ \beta - (s-1) + \alpha_2-1 } D^{\bs \gamma}v}_{L^2(\widetilde{K}_{\bs j})} \right) \\
		&\leq C \, h^{s-1} \left(\sum_{\alpha_2 =0}^{1} \sum_{\abs{\bs \gamma} \leq s+\alpha_2-1}  \norm{r^{\beta - s + \abs{\bs \gamma}} D^{\bs \gamma}v}_{L^2(\widetilde{K}_{\bs j})} + \sum_{\alpha_2=0}^{s} \sum_{\abs{\bs \gamma} \leq \alpha_2}  \norm{r^{\beta - s + \abs{\bs \gamma}} D^{\bs \gamma}v}_{L^2(\widetilde{K}_{\bs j})} \right) \\
		&\leq C \, h^{s-1} \left(\sum_{\abs{\bs \gamma}\leq s}  \norm{r^{\beta - s + \abs{\bs \gamma}} D^{\bs \gamma}v}_{L^2(\widetilde{K}_{\bs j})} + \sum_{\abs{\bs \gamma} \leq s}  \norm{r^{\beta - s + \abs{\bs \gamma}} D^{\bs \gamma}v}_{L^2(\widetilde{K}_{\bs j})} \right) 
		\leq C \, h^{s-1}\norm{v}_{V^{s}_{\beta}(\widetilde{K}_{\bs j})} .
	\end{align*}
	The second term in line \eqref{eq: proof of ring domain terms negative weight} can be computed similarly,
	\begin{align*}
		\quad \zeta_{1,j_1}^{-1/2} \norm{\left( W \widehat{v} -\Pi_{\bs p, \bs \Xi}( W \widehat{v}) \right)}_{\widehat L^2(Q)}^2
		&\leq C \left( \widetilde{h}_{1,j_1}^{s} \norm{\widehat{D}^{(s,0)}\widehat v}_{\mathcal{L}^2(\widetilde{Q}_{\bs j})} + \widetilde{h}_{2,j_2}^{s} \sum_{\alpha_2=0}^{s}  \norm{ \widehat{D}^{(0, \alpha_2)}\widehat v}_{\mathcal{L}^2(\widetilde{Q}_{\bs j})} \right) \\
		&\leq C \, \zeta_{1,j_1}^{-1} \left( \widetilde{h}_{1,j_1}^{s} \sum_{\abs{\bs \gamma} \leq s} \norm{D^{\bs \gamma}v }_{\mathcal{L}^2(\widetilde{K}_{\bs j})} + \widetilde{h}_{2,j_2}^{s} \sum_{\alpha_2=0}^{s} \sum_{\abs{\bs \gamma} \leq \alpha_2} \norm{r^{\alpha_2} D^{\bs \gamma}v}_{\mathcal{L}^2(\widetilde{K}_{\bs j})} \right) \\
		&\leq C \, h^{s-1} \left(\sum_{\abs{\bs \gamma} \leq s} \norm{r^{\beta} D^{\bs \gamma}v}_{L^2(\widetilde{K}_{\bs j})} + \sum_{\alpha_2=0}^{s} \sum_{\abs{\bs \gamma} \leq \alpha_2} \norm{r^{ \beta - s + \alpha_2 } D^{\bs \gamma}v}_{L^2(\widetilde{K}_{\bs j})} \right) \\
		&\leq C \, h^{s-1} \left(\norm{v}_{H^s_{\beta}(\widetilde{K}_{\bs j})} + \sum_{\alpha_2=0}^{s} \sum_{\abs{\bs \gamma} \leq \alpha_2}  \norm{r^{\beta - s + \abs{\bs \gamma}} D^{\bs \gamma}v}_{L^2(\widetilde{K}_{\bs j})} \right) \\
		&\leq C \, h^{s-1} \norm{v}_{V^{s}_{\beta}(\widetilde{K}_{\bs j})} .
	\end{align*}
	Finally, we sum over all elements and remark that the overlap of support extensions depends only on the polynomial degree $p$ and is thus bounded with respect to $h$. We obtain
	\begin{align*}
		\abs{v - \Pi_{V_h^{\pol}} v}^2_{H^q(\Omega_R)} 
		= \sum_{K \in \mathcal{M}_R^{\mu}} \abs{v - \Pi_{V_h}^{\pol} v}^2_{H^q(K)}
		\leq C \, h^{2(s-q)} \sum_{K \in \mathcal{M}_R^{\mu}} \norm{v}^2_{V^{s}_{\beta}(\widetilde{K}_{\bs j})}
		\leq C h^{2(s-q)}  \norm{v}^2_{V^{s}_{\beta}(\widetilde{\Omega}_R)}
	\end{align*}
	for $q\in\{0,1\}$ and the assertion follows.
\end{proof}

\subsubsection{Estimates on $\Omega_C$}
Next, we analyze the projection error on the small polar domain $\Omega_{C}$ that is surrounding the polar point. Therefore, we adopt the proof strategy from \cite{ApelSaendigSolovyev2002}, where a similar error estimate has been shown for interpolation with linear finite elements. The radial size of the corresponding parametric strip $\widehat{\Omega}_{C}$ will be denoted by $h_{r, \widehat{\Omega}_C}=\zeta_{1,p_1+2}$ such that $\widehat{\Omega}_C = \left(0, h_{r, \widehat{\Omega}_C}\right] \times (0,1)$, recall \eqref{eq: definition of parametric strips}.

\begin{Theorem}
	\label{theorem: local projection error polar stripe}
	Let $q \in \{0,1\}$, $s \in \N$ with $2 \leq s \leq p +1$ and $v \in V^{s}(\widetilde{\Omega}_{C})$. Then, there exists a positive constant $C$ depending only on $p,\bs \theta, \bs F, W$ such that
	\begin{equation*}
		\abs{v - \Pi_{V_h^{\pol}} v}_{H^q(\Omega_{C})} \leq C h_{r, \widehat{\Omega}_C}^{s-q} \norm{v}_{V^{s}(\widetilde{\Omega}_C)} , 
	\end{equation*}
	with $\widetilde{\Omega}_C$ as defined in \eqref{eq: overlapping splitting}.
\end{Theorem}

\begin{proof}
	As in the proof of Theorem \ref{theorem: projection error ring domain}, we transform the projection error to the polar coordinate system as illustrated in Section \ref{subsec: transformation of projection error}. For $q=0$, we thus have
	\begin{align*}
		\norm{v -  \Pi_{V_h^{\pol}} v }^2_{L^2(\Omega_C)}
		\leq C \norm{W \widehat{v} - \Pi^{\pol}_{\bs p , \bs \Xi} (W \widehat{v})}_{ \mathcal{L}^2_{1/2}(\widehat{\Omega}_C)}^2 ,
	\end{align*}
	and for $q=1$, it follows
	\begin{align}
		\abs{v -  \Pi_{V_h^{\pol}} v }^2_{H^1(\Omega_C)}
		&\leq C \bigg( \norm{\partial_r \left( W \widehat{v} - \Pi^{\pol}_{\bs p, \bs \Xi}( W \widehat{v}) \right)}_{\mathcal L^2_{1/2}(\widehat{\Omega}_C)}^2 + \norm{W \widehat{v} -\Pi^{\pol}_{\bs p, \bs \Xi} (W \widehat{v}) }_{\mathcal L^2_{1/2}(\widehat{\Omega}_C)}^2 \nonumber \\
		& \ + \norm{\partial_\varphi\left( W \widehat{v} -\Pi^{\pol}_{\bs p, \bs \Xi}( W \widehat{v}) \right)}_{\mathcal L^2_{-1/2}(\widehat{\Omega}_C)}^2 + \norm{W \widehat{v} - \Pi^{\pol}_{\bs p, \bs \Xi} (W \widehat{v}) }_{\mathcal L^2_{-1/2}(\widehat{\Omega}_C)}^2 \bigg) .
		\label{eq: error components H^1 projection}
	\end{align}
	However, in contrast to Theorem \ref{theorem: projection error ring domain}, where the projection error on $\Omega_R$ is considered, the Jacobian of the inverse of the polar parameterization $\bs F$ is not bounded on $\Omega_{C}$, see \eqref{eq: Jacobian of polar parameterization collapses}. In the following Lemmata \ref{lemma: estimate of r-derivative}, \ref{lemma: estimate of phi-derivative} and \ref{lemma: estimate of L_{-1/2}-norm}, we will thus show that each of the weighted error terms is bounded by the bent polar Sobolev norm of $v$ introduced in Section \ref{subsec: bent spaces}. With equation \eqref{eq: bent polar spaces and 1/2-spaces} and the transformation  \eqref{eq: transformation polar bent Sobolev spaces} back to the physical domain, we obtain the assertion,
	\begin{align*}
		\abs{v -  \Pi_{V_h^{\pol}} v }_{H^q(\Omega_C)}
		\leq C h_{r, \widehat{\Omega}_C}^{s-q} \norm{\widehat{v}}_{\mathcal{H}_{\pol}^{s}\left(\widetilde{\widehat{\Omega}}_C\right)} 
		\leq C h_{r, \widehat{\Omega}_C}^{s-q} \norm{\widehat{v}}_{\mathcal{V}_{\pol}^{s}\left(\widetilde{\widehat{\Omega}}_C\right)} 
		\leq C h_{r, \widehat{\Omega}_C}^{s-q} \norm{v}_{V^{s}(\widetilde \Omega_C)} ,
	\end{align*}	
	where $\widetilde{\widehat{\Omega}}_C$ and $\widetilde \Omega_C$ are defined as in \eqref{eq: overlapping splitting}.
\end{proof}

\begin{figure}[t]
	\begin{center}
		\input{ref_elem_ext.tikz}
		\input{ref_elem_ext_ori.tikz}
	\end{center}	
	
	\vspace{-18mm}
	
	\begin{center}
		\hspace{1mm}
		\input{ref_elem_ext_phys.tikz}
		\hspace{5mm}
		\input{ref_elem_ext_ori_phys.tikz}
		\caption{First row: Illustration of the parametric scaling transformation \eqref{eq: scaling transformation}, with the extended parametric strip $\widetilde{\widehat{\Omega}}_C$ and the corresponding reference domain $\widetilde{\mr{\Omega}}_C$ being hatched in green. Second row: Illustration of the physical scaling transformation \eqref{eq: scaling transformation physical reference domain}, with the extended small polar domain $\widetilde{\Omega}_C$ and the corresponding reference domain $\widetilde{\breve{\Omega}}_C$ being hatched in green.}
		\label{fig: illustration of support extension of reference domain}
	\end{center}
\end{figure}
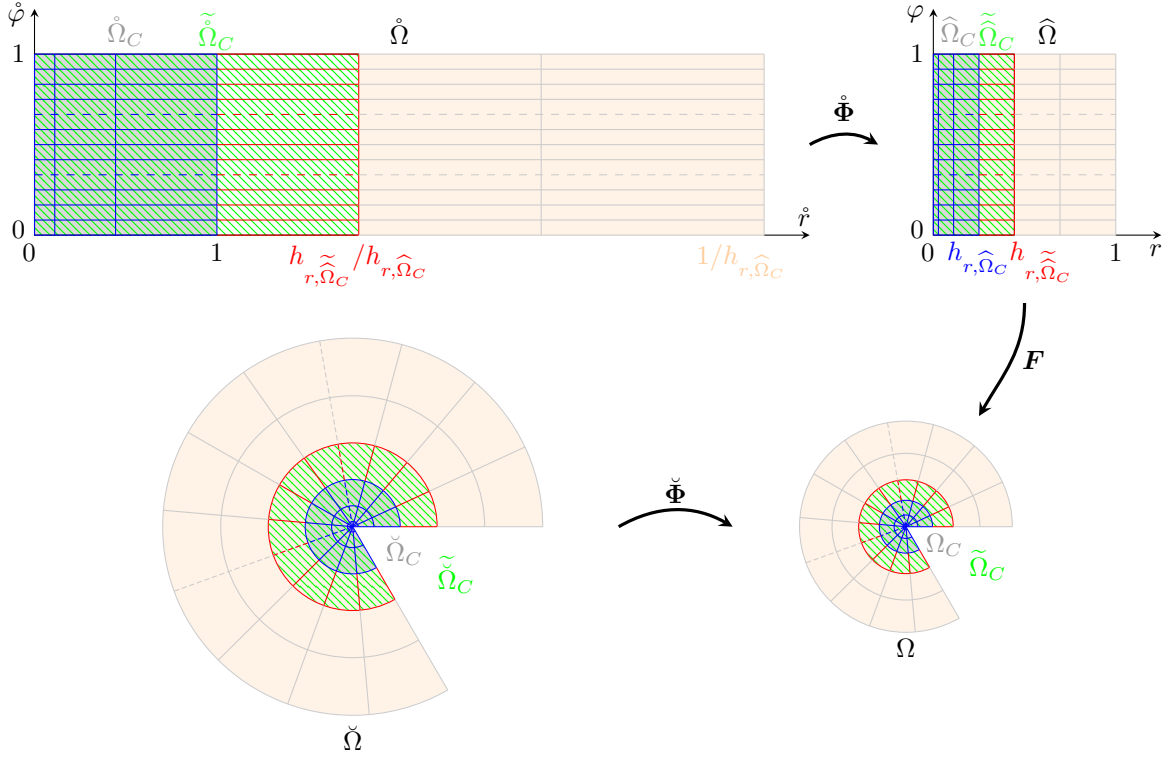

An essential ingredient for the following estimates is the bijective scaling transformation between the parametric domain $\widehat{\Omega}$ and the reference domain $\mr\Omega := (0,1/h_{r,\widehat{\Omega}_C}) \times (0,1)$, given by
\begin{align}
	\mr{\bs \Phi} : \mr\Omega \to \widehat{\Omega}, \quad \mr{\bs \Phi}(\mr r,\mr \varphi) 
	= \begin{pmatrix} h_{r,\widehat{\Omega}_C} \mr{r} \\ \mr \varphi \end{pmatrix} 
	= \begin{pmatrix} r \\ \varphi \end{pmatrix}.
	\label{eq: scaling transformation}
\end{align}
In this way, the size of the reference domain corresponding to the parametric strip $\widehat{\Omega}_C$, namely
\begin{align*}
	\mr\Omega_C := \mr{\bs \Phi}^{-1} \left(\widehat{\Omega}_C\right) = (0,1)^2,
\end{align*}
is independent of the global mesh size $h$, and a typical scaling argument can be applied. The transformation and reference domains are illustrated visually in Figure \ref{fig: illustration of reference domain}. Any quantity that is transformed to the reference domain will be marked accordingly with a circle on top, e.g., $\widehat{v} \in H_{\pol}^s(\widehat\Omega)$ becomes $\mr{v} \in H_{\pol}^s(\mr\Omega)$, and will be called the corresponding reference quantity or scaled quantity. Moreover, we define the reference domain of the extended parametric strip $\widetilde{\widehat{\Omega}}_{C}$, recall \eqref{eq: overlapping splitting}, by
\begin{align}
	\widetilde{\mr{\Omega}}_C := \mr{\bs \Phi}^{-1} \left( \widetilde{\widehat{\Omega}}_{C} \right)
	= \left(0, h_{r,\widetilde{\widehat{\Omega}}_{C}} \Big/ h_{r, \widehat{\Omega}_C}\right)\times (0,1) ,
	\label{eq: reference of extended parametric strip}
\end{align}
which is displayed in the top row of Figure \ref{fig: illustration of support extension of reference domain}. Its size is also independent of the mesh size, since the relation $h_{r,\widetilde{\widehat{\Omega}}_{C}} / h_{r, \widehat{\Omega}_C} \sim 1$ holds due to the local quasi-uniformity of the mesh, recall Lemma \ref{lemma: local quasi-uniformity of the graded mesh}. Finally, we introduce the corresponding physical reference domains, which are obtained using the transformation
\begin{align}
	\breve{\bs \Phi}: \breve{\Omega} \to \Omega , \quad \breve{\bs \Phi}(\breve x, \breve y) 
	= \begin{pmatrix} h_{r,\widehat{\Omega}_C} \breve x \\ h_{r,\widehat{\Omega}_C} \breve y  \end{pmatrix} 
	= \begin{pmatrix}  x \\  y \end{pmatrix} , 
	\label{eq: scaling transformation physical reference domain}
\end{align}
with $\breve{\Omega} := \breve{\bs \Phi}^{-1} (\Omega)$, as illustrated in the bottom row of Figure \ref{fig: illustration of support extension of reference domain}. Similarly to the parametric reference configuration, the size of $\breve{\Omega}_C := \breve{\bs \Phi}^{-1} \left({\Omega}_C\right)$ and $\widetilde{\breve\Omega}_C := \breve{\bs \Phi}^{-1} \left(\widetilde{\Omega}_C\right)$ is independent of $h$. Any scaled quantity will be marked accordingly with a breve on top.

\begin{figure}
	\begin{center}
		\input{ref_elem.tikz}
		\input{ref_elem_ori.tikz}
		\hspace{1mm}
		\caption{Illustration of the scaling transformation $\mr {\bs \Phi}$ with the reference domains $\mr{\Omega}$ and $\mr{\Omega}_C$ and the corresponding parametric domains $\widehat{\Omega}$ and $\widehat{\Omega}_C$, respectively.}
		\label{fig: illustration of reference domain}
	\end{center}
\end{figure}
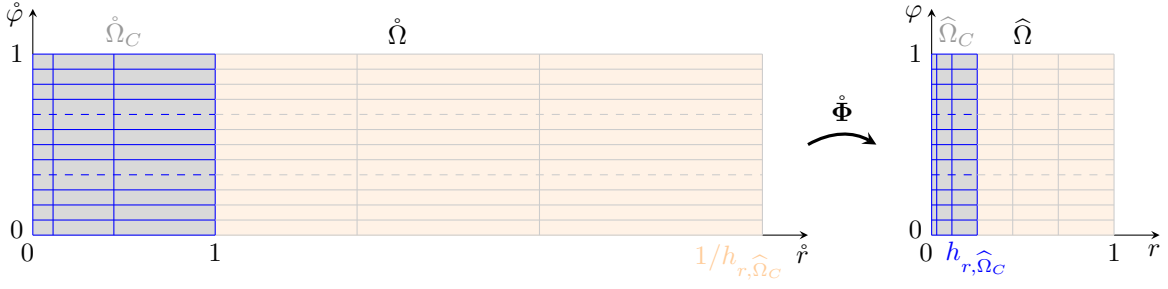

We start with a stability result for the scaled version of the quasi-interpolant defined in \eqref{eq: modified quasi-interpolant}, which will be used more often in the remainder of this section.
\begin{lemma}
	\label{lemma: stability of quasi-interpolant}
	Let either $\mr v \in  \widehat H^{s}_{1/2}\left(\widetilde{\mr{\Omega}}_C\right)$ or $ \mr v \in \mathcal H^{s}_{1/2}\left(\widetilde{\mr{\Omega}}_C\right)$ for $s \geq 2$. Then, it holds
	\begin{align*}
		\norm{\partial^q_{\mr r} \left(\mr{\Pi}^{\pol}_{\bs p, \bs \Xi} \, \mr{v}\right)}_{\widehat L^2_{1/2}(\mr\Omega_C)} ^2
		\leq C \norm{\mr{v}}_{\widehat H^{s}_{1/2}\left(\widetilde{\mr{\Omega}}_C\right)}^2
		\quad \text { or } \quad
		\norm{\partial^q_{\mr r} \left(\mr{\Pi}^{\pol}_{\bs p, \bs \Xi} \, \mr{v}\right)}_{\mathcal L^2_{1/2}(\mr\Omega_C)} ^2
		\leq C \norm{\mr{v}}_{\mathcal H^{s}_{1/2}\left(\widetilde{\mr{\Omega}}_C\right)}^2 ,
	\end{align*}
	respectively, for $q \in \{0,1\}$.
\end{lemma}

\begin{proof}
	Let $q \in \{0,1\}$ and $\mr v \in  \widehat H^{s}_{1/2}\left(\widetilde{\mr{\Omega}}_C\right)$. By definition of $\mr{\Pi}^{\pol}_{\bs p, \bs \Xi}$, recall \eqref{eq: modified quasi-interpolant}, and the stability property of the dual functionals \cite[Theorem 4.41]{Schumaker2007}, we compute
	\begin{align*}
		\norm{\partial^q_{\mr r} \left(\mr{\Pi}^{\pol}_{\bs p, \bs \Xi} \, \mr{v}\right)}_{\widehat L^2_{1/2}(\mr\Omega_C)} ^2
		&= \int_{\mr \Omega_C} \abs{\partial^q_{\mr r} \left(\sum_{\bs i \in \bs I} {\mr \lambda}^{\pol}_{\bs i, \bs p}(\mr{v}) \mr{B}_{\bs i,\bs p}\right)}^2 \mr{r} \Dd \mr r \Dd \mr \varphi
		\leq \int_{\mr \Omega_C} \left( \sum_{\bs i \in \bs I} \abs{{\mr \lambda}^{\pol}_{\bs i, \bs p}(\mr{v}) \, \partial^q_{\mr r}  \mr{B}_{\bs i,\bs p}} \right)^2 \mr{r} \Dd \mr r \Dd \mr \varphi \\
		& \leq \int_{\mr \Omega_C} \left(\sum_{\bs i \in \bs I} C \norm{\mr v}_{L^{\infty}\left(\widetilde{\mr \Omega}_C\right)} \abs{\partial^q_{\mr r} \mr{B}_{\bs i,\bs p}} \right)^2 \mr{r} \Dd \mr r \Dd \mr \varphi \\
		&\leq C \norm{\mr v}_{C^{0}\left(\overline{\widetilde{\mr \Omega}_C}\right)}^2  \int_{\mr \Omega_C} \left(\sum_{\bs i \in \bs I} \abs{\partial^q_{\mr r}  \mr{B}_{\bs i,\bs p}} \right)^2 \mr r \Dd \mr r \Dd \mr \varphi \\
		&\leq C \norm{\mr{v}}_{C^0\left(\overline{\widetilde{\mr{\Omega}}_C}\right)}^2  \int_{\mr \Omega_C} 1 \Dd \mr r \Dd \mr \varphi 
		\leq C \norm{\mr{v}}_{C^0\left(\overline{\widetilde{\mr{\Omega}}_C}\right)}^2 
		\leq C \norm{\mr{v}}_{\widehat H^{s}_{1/2}\left(\widetilde{\mr{\Omega}}_C\right)}^2 ,
	\end{align*}
	where, in the last step, we used the embedding $\widehat H^{s}_{1/2}\left(\widetilde{\mr{\Omega}}_C\right) \hookrightarrow C^0\left(\overline{\widetilde{\mr{\Omega}}_C}\right)$, see \cite[Theorem 4.7]{MercierRaugel1982}. To prove the assertion for bent Sobolev norms, we simply apply the estimate above on every mesh element and sum over all elements.
\end{proof}

Now, we estimate the first two terms of \eqref{eq: error components H^1 projection}, that is, the projection error and its $r$-derivative in the $\mathcal L^2_{1/2}(\widehat{\Omega}_C)$-norm.
\begin{lemma}
	Let the assumptions of Theorem \ref{theorem: local projection error polar stripe} hold. Then, it follows for $q\in \{0,1\}$
	\begin{align*}
		\norm{\partial^q_r\left( W \widehat{v} - \Pi^{\pol}_{\bs p, \bs \Xi}(W\widehat{v})\right)}_{\mathcal L^2_{1/2}(\widehat{\Omega}_C)} 
		\leq C h_{r, \widehat{\Omega}_C}^{s-q} \norm{\widehat{v}}_{\mathcal{H}^{s}_{\pol}\left(\widetilde{\widehat{\Omega}}_{C}\right)} .
	\end{align*}
	\label{lemma: estimate of r-derivative}
\end{lemma}

\begin{proof}
	First, we apply a typical scaling argument using the reference configuration introduced above and, to simplify notation, we set $\mr v = \mr W \mr v$. It follows
	\begin{align}
		\norm{\partial^q_r\left( W \widehat{v} - \Pi^{\pol}_{\bs p, \bs \Xi}(W\widehat{v})\right)}_{\mathcal L^2_{1/2}(\widehat{\Omega}_C)}
		&= \left(\int_{\widehat{\Omega}_C}  \left(\partial^q_r \left(W \widehat v -  \Pi^{\pol}_{\bs p, \bs \Xi}(W \widehat v) \right) \right)^2 r \dr \dphi \right)^{1/2} \nonumber \\
		&= \left( \int_{\mr\Omega_C} \left( \frac{\partial^q_{\mr r} \left( \mr W \mr{v} - \mr{\Pi}^{\pol}_{\bs p, \bs \Xi} (\mr W \mr{v})\right)}{h_{r,\widehat{\Omega}_C}^q} \right) ^2  h_{r, \widehat{\Omega}_C} \, \mr{r} \, h_{r, \widehat{\Omega}_C} \, \Dd \mr{r} \Dd \mr{\varphi} \right)^{1/2} \nonumber \\
		&= h_{r, \widehat{\Omega}_C} ^{1-q}  \norm{\partial^q_{\mr r} \left(\mr{v} -\mr{\Pi}^{\pol}_{\bs p, \bs \Xi} \mr{v}\right)}_{\mathcal L^2_{1/2}(\mr\Omega_C)}. \label{eq: scaling argument r-derivative}
	\end{align}
	In the rest of the paper, such estimates will be used more often. The underlying steps will not be carried out in detail anymore, but we will refer to it as scaling argument.
	
	Next, let $\mathcal{P}_{s-1}$ be the space of piecewise polynomials up to degree $s-1$ on $\mr{\Omega}_C$ and $\mr{\Pi}$ be the reference projector from the bent Sobolev space to the spline space of degree $s-1$ that satisfies property \eqref{eq: projector for bent polar spaces}. Due to Lemma \ref{lemma: stability of quasi-interpolant}, we obtain
	\begin{align*}
		\norm{\partial^q_{\mr r} \left(\mr{v} -\mr{\Pi}^{\pol}_{\bs p, \bs \Xi} \mr{v}\right)}_{\mathcal L^2_{1/2}(\mr\Omega_C)}
		& = \inf_{\mr{w} \in \mathcal{P}_{s-1}} \norm{\partial^q_{\mr r} \left(\mr{v} - \mr{\Pi}(\mr v) - \mr{w} - \mr{\Pi}^{\pol}_{\bs p, \bs \Xi}(\mr{v} - \mr{\Pi}(\mr v) - \mr{w})\right)}_{\mathcal L^2_{1/2}(\mr\Omega_C)} \\
		& = \inf_{\mr{w} \in \mathcal{P}_{s-1}} \left( \norm{\partial^q_{\mr r} \left(\mr{v} - \mr{\Pi}(\mr v) - \mr{w}\right)}_{\mathcal L^2_{1/2}(\mr\Omega_C)} + \norm{\partial^q_{\mr r} \left(\mr{\Pi}^{\pol}_{\bs p, \bs \Xi} (\mr{v} - \mr{\Pi}(\mr v) - \mr{w})\right)}_{\mathcal L^2_{1/2}(\mr\Omega_C)}\right) \\
		& \leq C \inf_{\mr{w} \in \mathcal{P}_{s-1}} \norm{\mr{v} - \mr{\Pi}(\mr v) - \mr{w}}_{\mathcal H^{s}_{1/2}\left(\widetilde{\mr{\Omega}}_C\right)} 
		= \inf_{\mr{w} \in \mathcal{P}_{s-1}} \norm{\mr{v} - \mr{\Pi}(\mr v) - \mr{w}}_{\widehat H^{s}_{1/2}\left(\widetilde{\mr{\Omega}}_C\right)} .
	\end{align*}
	Then, we apply the weighted Deny-Lions-type argument \cite[Theorem 4.6]{MercierRaugel1982},
	\begin{align}
		\label{eq: deny-lions argument higher order}
		\inf_{\mr{w}\in \mathcal{P}_{s-1}} \norm{\mr{z} - \mr{w}}_{\widehat H^{s}_{1/2}\left(\widetilde{\mr{\Omega}}_C\right)}
		\leq C \abs{\mr{z}}_{\widehat H^{s}_{1/2}\left(\widetilde{\mr{\Omega}}_C\right)} 
	\end{align}
	for $\mr{z} = \mr{v} - \mr{\Pi}(\mr v) \in \widehat H^{s}_{1/2}\left(\widetilde{\mr{\Omega}}_C\right)$. Finally, we combine all the steps and use again the definition of $\mr{\Pi}$, a scaling argument, the embedding \eqref{eq: bent polar spaces and 1/2-spaces} and the Leibniz formula \eqref{eq: Leibniz formula} to compute
	\begin{align*}
		\norm{\partial^q_r\left( W\widehat{v} - \Pi^{\pol}_{\bs p, \bs \Xi}(W\widehat{v})\right)}_{\widehat L^2_{1/2}(\widehat{\Omega}_C)}
		&\leq C  h_{r, \widehat{\Omega}_C}^{1-q} \inf_{\mr{w} \in \mathcal{P}_{s-1}} \norm{\mr{v} - \mr{\Pi}(\mr v) - \mr{w}}_{\widehat H^{s}_{1/2}\left(\widetilde{\mr{\Omega}}_C\right)}
		\leq C h_{r, \widehat{\Omega}_C}^{1-q} \abs{\mr{v} - \mr{\Pi}(\mr v)}_{\widehat H^{s}_{1/2}\left(\widetilde{\mr{\Omega}}_C\right)}  \\
		& = C h_{r, \widehat{\Omega}_C}^{1-q} \abs{\mr{v}}_{\mathcal{H}^{s}_{1/2}\left(\widetilde{\mr{\Omega}}_C\right)}
		\leq C h_{r, \widehat{\Omega}_C}^{s-q} \abs{W\widehat{v}}_{\mathcal{H}^{s}_{1/2}\left(\widetilde{\widehat{\Omega}}_{C}\right)}
		\leq C h_{r, \widehat{\Omega}_C}^{s-q} \norm{\widehat{v}}_{\mathcal{H}^{s}_{\pol}\left(\widetilde{\widehat{\Omega}}_{C}\right)} ,
	\end{align*}
	and the demonstration is complete.
\end{proof} 

In the proof of Lemma \ref{lemma: estimate of r-derivative}, the projector $\mr{\Pi}$ that satisfies property \eqref{eq: projector for bent polar spaces} was only inserted for the interplay between bent Sobolev norms and their non-bent counterparts, recall Section \ref{subsec: bent spaces}. This adds a further technicality to the proof, which is not in the scope of this paper, and has already been investigated in detail in standard approximation literature \cite{BeiraodaVeigaBuffaSangalliVazquez2014}. For the sake of simplicity, we will not elaborate this argument anymore in the rest of our paper and, w.l.o.g., we will assume that pull-backs of functions are sufficiently smooth across mesh lines. The generalization can always be done by inserting the projector $\mr{\Pi}$ at the correct places, as shown in the demonstration of Lemma \ref{lemma: estimate of r-derivative}.

To continue the proof of Theorem \ref{theorem: local projection error polar stripe}, we estimate the third term in \eqref{eq: error components H^1 projection} in the next Lemma.

\begin{lemma}
	Let the assumptions of Theorem \ref{theorem: local projection error polar stripe} hold. Then, it follows
	\begin{align*}
		\norm{\partial_\varphi\left(W\widehat{v} - \Pi^{\pol}_{\bs p, \bs \Xi} (W\widehat{v})\right)}_{\widehat{L}^2_{-1/2}(\widehat{\Omega}_C)} 
		\leq C h_{r, \widehat{\Omega}_C}^{s-1} \norm{\widehat{v}}_{\mathcal{H}^{s}_{\pol}\left(\widetilde{\widehat{\Omega}}_{C}\right)} .
	\end{align*}
	\label{lemma: estimate of phi-derivative}
\end{lemma}

\begin{proof} 
	Let $c \in \R$ such that $\widehat{v}(0,\cdot)= v \circ \bs F(0,\cdot) = v(\bs P) = c$. We define the auxiliary function $\widehat{v}_0 = \widehat v - c \widehat{B}_{1,p_1}$, which satisfies $(W\widehat{v}_0)(0,\cdot) = (W  \widehat v) (0,\cdot) - c W \widehat{B}_{1,p_1}(0) = 0$. Moreover, we have
	\begin{align*}
		\Pi^{\pol}_{\bs p, \bs \Xi} (W \widehat{v}_0)
		=\Pi^{\pol}_{\bs p, \bs \Xi} (W \widehat{v}) - \Pi^{\pol}_{\bs p, \bs \Xi} (c W B_{1,p_1}) 
		=\Pi^{\pol}_{\bs p, \bs \Xi} (W \widehat{v}) - c W B_{1,p_1}
	\end{align*}
	and it follows $W \widehat{v}_0 - \Pi^{\pol}_{\bs p, \bs \Xi} (W\widehat{v}_0) = W\widehat{v} - \Pi^{\pol}_{\bs p, \bs \Xi} (W\widehat{v})$ and $\Pi^{\pol}_{\bs p, \bs \Xi}( W \widehat{v}_0)(0,\cdot) = 0$. Using a scaling argument and the notation $\mr v_0 = \mr W \mr v_0$, we compute
	\begin{align*}
		\norm{\partial_\varphi\left(W\widehat{v} - \Pi^{\pol}_{\bs p, \bs \Xi} (W\widehat{v})\right)}_{\widehat L^2_{-1/2}(\widehat{\Omega}_C)} 
		& = \norm{\partial_\varphi\left( W \widehat{v}_0 - \Pi^{\pol}_{\bs p, \bs \Xi} (W\widehat{v}_0)\right)}_{\widehat L^2_{-1/2}(\widehat{\Omega}_C)} \nonumber \\
		&\leq \norm{\partial_\varphi (W\widehat{v}_0)}_{\widehat L^2_{-1/2}(\widehat{\Omega}_C)} + \norm{\partial_\varphi\left( \Pi^{\pol}_{\bs p, \bs \Xi} (W\widehat{v}_0) \right)}_{\widehat L^2_{-1/2}(\widehat{\Omega}_C)} \\ 
		&\leq \norm{\partial_\varphi \mr{v}_0}_{\widehat L^2_{-1/2}(\mr{\Omega}_C)} + \norm{\partial_{\mr \varphi}\left( \mr{\Pi}^{\pol}_{\bs p, \bs \Xi} \mr{v}_0\right)}_{\widehat L^2_{-1/2}(\mr{\Omega}_C)},
	\end{align*}
	where all terms are well-defined despite the negative weight of the norms as the corresponding functions vanish for $r=0$. In Lemma \ref{lemma: derivative of projector}, we show that the $\varphi$-derivative of our projector can be written as another projector, and in Lemma \ref{lemma: estimate of derivative of projector}, we estimate the corresponding norm. Hence, we obtain 
	\begin{align*}
		\norm{\partial_\varphi\left(W\widehat{v} - \Pi^{\pol}_{\bs p, \bs \Xi} (W\widehat{v})\right)}_{\widehat L^2_{-1/2}(\widehat{\Omega}_C)}  
		& \leq C \norm{\partial_{\mr \varphi} \mr{v}_0}_{\widehat L^2_{-1/2}\left(\widetilde{\mr{\Omega}}_C\right)} \leq C \norm{\partial_{\mr \varphi}\mr{v}_0}_{\widehat L^{2}_{-1}\left(\widetilde{\mr{\Omega}}_C\right)}.
	\end{align*}
	It is shown in \cite[Corollary 4.1]{MercierRaugel1982} that, under the condition $\partial_{\mr \varphi}\mr{v}_0 (0, \cdot)=0$, the inequality
	\begin{align*}
		\norm{\partial_{\mr \varphi}\mr{v}_0}_{\widehat L^{2}_{-1}\left(\widetilde{\mr{\Omega}}_C\right)}
		\leq C \abs{\partial_{\mr \varphi}\mr{v}_0}_{\widehat H^{1}\left(\widetilde{\mr{\Omega}}_C\right)} 
	\end{align*}
	holds. By a scaling argument, it follows
	\begin{align*}
		\abs{\partial_{\mr \varphi}\mr{v}_0}_{\widehat H^{1}\left(\widetilde{\mr{\Omega}}_C\right)} 
		&\leq C \left( \norm{\partial_{\mr r \mr \varphi} \mr{v}_0 }_{\widehat L^2\left(\widetilde{\mr{\Omega}}_C\right)} + \norm{\partial_{\mr \varphi \mr \varphi} \mr{v}_0 }_{\widehat L^2\left(\widetilde{\mr{\Omega}}_C\right)} \right) \nonumber \\
		&\leq C \left( \norm{\partial_{\mr r \mr \varphi} \mr{v}_0 }_{\widehat L^2_{-1/2}\left(\widetilde{\mr{\Omega}}_C\right)} + \norm{\partial_{\mr \varphi \mr \varphi} \mr{v}_0 }_{\widehat L^2_{-3/2}\left(\widetilde{\mr{\Omega}}_C\right)} \right) \\
		&\leq C h_{r, \widehat{\Omega}_C} \left( \norm{\partial_{r \varphi} (W\widehat{v}_0) }_{\widehat L^2_{-1/2}\left(\widetilde{\widehat{\Omega}}_{C}\right)} + \norm{\partial_{\varphi \varphi} (W\widehat{v}_0) }_{\widehat L^2_{-3/2}\left(\widetilde{\widehat{\Omega}}_{C}\right)} \right) \nonumber \\
		& \leq C h_{r, \widehat{\Omega}_C} \norm{\widehat{v}}_{\mathcal{H}^{2}_{\pol}\left(\widetilde{\widehat{\Omega}}_{C}\right)} , \nonumber
	\end{align*}
	which yields the assertion for $s=2$. For higher regularity $s \geq 3$, the proof needs to be adapted. We rewrite the error in terms of the projector $\mr{\Pi}^{\pol,\partial_{\varphi}}_{\bs p, \bs \Xi}$ derived in Lemma \ref{lemma: derivative of projector}, and compute
	\begin{align}
		\norm{\partial_\varphi\left(W\widehat{v} - \Pi^{\pol}_{\bs p, \bs \Xi}( W \widehat{v})\right)}_{\widehat L^2_{-1/2}(\widehat{\Omega}_C)}
		&= \norm{\partial_{\mr \varphi}\left( \mr{v}_0 - \mr{\Pi}^{\pol}_{\bs p, \bs \Xi} \mr{v}_0\right)}_{\widehat L^2_{-1/2}(\mr \Omega_C)}
		= \norm{\partial_{\mr \varphi} \mr{v}_0 - \mr{\Pi}^{\pol,\partial_\varphi}_{(p_1,p_2-1), \bs \Xi} (\partial_{\varphi}\mr{v}_0)}_{\widehat L^2_{-1/2}(\mr \Omega_C)} \nonumber \\
		& \leq C \norm{\partial_{\mr \varphi} \mr{v}_0 - \mr{\Pi}^{\pol,\partial_\varphi}_{(p_1,p_2-1), \bs \Xi} (\partial_{\varphi}\mr{v}_0)}_{\widehat L^2_{-1}(\mr \Omega_C)} \nonumber \\
		& \leq C \abs{\partial_{\mr \varphi} \mr{v}_0 - \mr{\Pi}^{\pol,\partial_\varphi}_{(p_1,p_2-1), \bs \Xi} (\partial_{\varphi}\mr{v}_0)}_{\widehat H^1(\mr \Omega_C)} \nonumber \\
		&= C \inf_{\mr{w} \in \mathcal{P}_{s-2}}\abs{\partial_{\mr \varphi} \mr{v}_0 - \mr{w} + \mr{\Pi}^{\pol,\partial_\varphi}_{(p_1,p_2-1), \bs \Xi} \left(\partial_{\varphi}\mr{v}_0 - \mr{w} \right)}_{\widehat H^1(\mr \Omega_C)} \nonumber \\
		&\leq C \inf_{\mr{w} \in \mathcal{P}_{s-2}} \left(\abs{\partial_{\mr \varphi} \mr{v}_0 - \mr{w}}_{\widehat H^1(\mr \Omega_C)} + \abs{\mr{\Pi}^{\pol,\partial_\varphi}_{(p_1,p_2-1), \bs \Xi}  \left(\partial_{\varphi}\mr{v}_0 - \mr{w} \right)}_{\widehat H^1(\mr \Omega_C)} \right) \nonumber \\
		&\leq C \inf_{\mr{w} \in \mathcal{P}_{s-2}} \norm{\partial_{\mr \varphi} \mr{v}_0 - \mr{w}}_{\widehat H^1\left(\widetilde{\mr{\Omega}}_C\right)} . 	\label{eq: estimate of phi-derivative with weights 1}
	\end{align}
	In the last step, we used the $H^1$-stability of the new projector $\mr{\Pi}^{\pol,\partial_{\varphi}}_{\bs p, \bs \Xi}$ which can be shown as in \cite{BeiraodaVeigaBuffaSangalliVazquez2014}. Finally, by applying the embedding $\widehat H^{s-1}_{1/2}(\mr \Omega_C) \hookrightarrow \widehat H^1(\mr \Omega_C)$ from \cite[Remark 4.1]{MercierRaugel1982}, the Deny-Lions type estimate \eqref{eq: deny-lions argument higher order}, a scaling argument and the relation \eqref{eq: bent polar spaces and 1/2-spaces}, we obtain
	\begin{align*}
		\norm{\partial_\varphi\left(W\widehat{v} - \Pi^{\pol}_{\bs p, \bs \Xi} (W\widehat{v})\right)}_{\widehat L^2_{-1/2}(\widehat{\Omega}_C)}
		&\leq C \inf_{\mr{w} \in \mathcal{P}_{s-2}} \norm{\partial_{\mr \varphi} \mr{v}_0 - \mr{w}}_{\widehat H^{s-1}_{1/2}\left(\widetilde{\mr{\Omega}}_C\right)} 
		\leq C \abs{\partial_{\mr \varphi} \mr{v}_0}_{\widehat H^{s-1}_{1/2}\left(\widetilde{\mr{\Omega}}_C\right)} \nonumber \\
		&\leq C \abs{\mr{v}_0}_{\widehat H^{s}_{1/2}\left(\widetilde{\mr{\Omega}}_C\right)} 
		\leq C h_{r, \widehat{\Omega}_C}^{s-1} \abs{W\widehat{v}_0}_{\widehat {H}^{s}_{1/2}\left(\widetilde{\widehat{\Omega}}_{C}\right)} \nonumber \\
		& \leq C h_{r, \widehat{\Omega}_C}^{s-1} \norm{\widehat{v}}_{\mathcal{H}^{s}_{1/2}\left(\widetilde{\widehat{\Omega}}_{C}\right)}
		\leq C h_{r, \widehat{\Omega}_C}^{s-1} \norm{\widehat{v}}_{\mathcal{H}^{s}_{\pol}\left(\widetilde{\widehat{\Omega}}_{C}\right)} ,
	\end{align*}
	and the assertion follows.
\end{proof}	

Finally, we estimate the last term of the sum \eqref{eq: error components H^1 projection}.
\begin{lemma}
	Let the assumptions of Theorem \ref{theorem: local projection error polar stripe} hold. Then, it follows
	\begin{align*}
		\norm{W\widehat{v} - \Pi^{\pol}_{\bs p, \bs \Xi}(W\widehat{v})}_{\widehat L^2_{-1/2}(\widehat \Omega_C)}
		\leq C h_{r, \widehat{\Omega}_C}^{s-1} \norm{\widehat{v}}_{\mathcal{H}^{s}_{\pol}\left(\widetilde{\widehat{\Omega}}_{C}\right)} .
	\end{align*}
	\label{lemma: estimate of L_{-1/2}-norm}
\end{lemma}

\begin{proof}
	As in the proof of Lemma \ref{lemma: estimate of phi-derivative}, we use again the auxiliary functions $\widehat{v}_0 = \widehat v - c \widehat{B}_{1,p_1}$ and set $\mr v_0 = \mr W \mr v_0$. By a scaling argument, it follows
	\begin{align*}
		\norm{W\widehat{v} - \Pi^{\pol}_{\bs p, \bs \Xi}(W\widehat{v})}_{\widehat L^2_{-1/2}(\widehat \Omega_C)}
		=\norm{W\widehat{v}_0 - \Pi^{\pol}_{\bs p, \bs \Xi}(W\widehat{v}_0)}_{\widehat L^2_{-1/2}(\widehat \Omega_C)}
		=\norm{\mr{v}_0 - \mr{\Pi}^{\pol}_{\bs p, \bs \Xi} \mr{v}_0}_{\widehat L^2_{-1/2}(\mr \Omega_C)} .
	\end{align*}
	Then, due to the $H^1$-stability of $\mr{\Pi}^{\pol}_{\bs p, \bs \Xi}$, we can show similar to \eqref{eq: estimate of phi-derivative with weights 1} that
	\begin{align*}
		\norm{\mr{v}_0 - \mr{\Pi}^{\pol}_{\bs p, \bs \Xi} \mr{v}_0}_{\widehat L^2_{-1/2}(\mr \Omega_C)} 
		\leq  C \inf_{\mr{w} \in \mathcal{P}_{s-1}} \norm{\mr{v}_0 - \mr w}_{\widehat H^1\left(\widetilde{\mr\Omega}_C\right)} .
	\end{align*} 
	With the embedding $\widehat H^{s}_{1/2}(\mr \Omega_C) \hookrightarrow \widehat H^1(\mr \Omega_C)$ for $s\geq 2$, see \cite[Remark 4.1]{MercierRaugel1982}, the Deny-Lions type estimate \eqref{eq: deny-lions argument higher order}, another scaling argument and the relation \eqref{eq: bent polar spaces and 1/2-spaces}, it follows
	\begin{align*}
		\norm{\mr{v}_0 - \mr{\Pi}^{\pol}_{\bs p, \bs \Xi} \mr{v}_0}_{\widehat L^2_{-1/2}(\mr \Omega_C)} 
		&\leq C \inf_{\mr{w} \in \mathcal{P}_{s-1}} \norm{\mr{v}_0 - \mr w}_{\widehat H^s_{1/2}\left(\widetilde{\mr\Omega}_C\right)} 
		\leq C \abs{\mr{v}_0}_{\widehat H^{s}_{1/2}\left(\widetilde{\mr{\Omega}}_C\right)} \\
		&\leq C h_{r, \widehat{\Omega}_C}^{s-1} \abs{W\widehat{v}_0}_{\widehat {H}^{s}_{1/2}\left(\widetilde{\widehat{\Omega}}_{C}\right)}
		\leq C h_{r, \widehat{\Omega}_C}^{s-1} \norm{\widehat{v}}_{\mathcal{H}^{s}_{\pol}\left(\widetilde{\widehat{\Omega}}_{C}\right)}
		\nonumber
	\end{align*}
	and the proof is complete.
\end{proof}

Recall that by proving Lemmata \ref{lemma: estimate of r-derivative}, \ref{lemma: estimate of phi-derivative} and \ref{lemma: estimate of L_{-1/2}-norm}, the proof of Theorem \ref{theorem: local projection error polar stripe} is complete. For now, we have assumed that $\widehat{v} \in \mathcal{H}^{s}_{\pol}\left(\widetilde{\widehat{\Omega}}_{C}\right)$ and $\norm{v}_{V^{s}(\widetilde \Omega_C)}< \infty$. However, as seen in Section \ref{subsec: model problems and regularity properties}, this is generally not the case for solutions of PDEs on polar domains with corners. Therefore, we show a corresponding result to Theorem \ref{theorem: local projection error polar stripe} for functions with reduced regularity by moving to a weighted function space setting.

\begin{Theorem}
	\label{theorem: local projection error polar stripe weighted}
	Let $q \in \{0,1\}$ and $s,s_0 \in \N$ with $2 \leq s \leq s_0 \leq p +1$. Further, let $v \in V^{s}_{\beta}(\widetilde{\Omega}_C)$ for all $s-1>\beta >s-1-\nu$. Then, there is a positive constant $C$, depending only on $p,\bs \theta, \bs F, W$ such that
	\begin{equation*}
		\abs{v - \Pi_{V_h^{\pol}} v}_{H^q(\Omega_{C})} \leq C h_{r, \widehat{\Omega}_C}^{s - q - \beta} \norm{v}_{V^{s}_{\beta}(\widetilde{\Omega}_C)} , 
	\end{equation*}
	with $\widetilde{\Omega}_C$ as defined in \eqref{eq: overlapping splitting}.
\end{Theorem}

\begin{proof}
	Just as in the proof of Theorem \ref{theorem: local projection error polar stripe}, we estimate the four error terms of the sum \eqref{eq: error components H^1 projection}. First, we adapt the stability result stated in Lemma \ref{lemma: stability of quasi-interpolant}, which is needed to estimate the first two terms. To that end, we use the embedding $H^{s}_{\beta}\left(\widetilde{\breve{\Omega}}_C\right) \hookrightarrow H^{s-\beta}\left(\widetilde{\breve{\Omega}}_C\right)$ from \cite{Rossmann1992} on the physical scaled reference domain defined in \eqref{eq: scaling transformation physical reference domain}, and combine it with the usual Sobolev embedding $H^{s-\beta}\left(\widetilde{\breve{\Omega}}_C\right) \hookrightarrow C^0\left(\overline{\widetilde{\breve{\Omega}}_C}\right)$ for $\beta<s-1$, see for instance \cite{AdamsFournier2003}. Then, with the proof of Lemma \ref{lemma: stability of quasi-interpolant}, embedding \eqref{eq: embedding V and H spaces physical domain} and transformation \eqref{eq: transformation of norm on single element adjacent to polar point}, we obtain 
	\begin{align}
		\label{eq: stability of projector weighted}
		\norm{\partial^q_{\mr r} \left(\mr{\Pi}^{\pol}_{\bs p, \bs \Xi} \, \mr{v}\right)}_{\widehat L^2_{1/2}(\mr\Omega_C)}
		&\leq C \norm{\mr{v}}_{C^0\left(\overline{\widetilde{\mr{\Omega}}_C}\right)}
		= C \norm{\breve{v}}_{C^0\left(\overline{\widetilde{\breve{\Omega}}_C}\right)}
		\leq C \norm{\breve{v}}_{H^{s-\beta}\left(\widetilde{\breve{\Omega}}_C\right)} \nonumber\\
		&\leq C \norm{\breve{v}}_{H^{s}_{\beta}\left(\widetilde{\breve{\Omega}}_C\right)}
		\leq C \norm{\breve{v}}_{V^{s}_{\beta}\left(\widetilde{\breve{\Omega}}_C\right)}
		\leq C \norm{\mr{v}}_{\widehat{V}^{s}_{\pol,\beta}\left(\widetilde{\mr{\Omega}}_{C}\right)} .
	\end{align}
	Moreover, we use transformation \eqref{eq: transformation of norm on single element adjacent to polar point} and the embeddings $V^{s-\beta}\left(\widetilde{\breve{\Omega}}_C\right) \hookrightarrow V^{1}\left(\widetilde{\breve{\Omega}}_C\right)$ from \cite[Lemma 2.29]{PfeffererPhDThesis} and $V^{s}_{\beta}\left(\widetilde{\breve{\Omega}}_C\right) \hookrightarrow V^{s-\beta}\left(\widetilde{\breve{\Omega}}_C\right)$ from \cite{Rossmann1992}, where $s-\beta>1$ is required, to show
	\begin{align}
		\label{eq: embedding weighted}
		\norm{ \mr{v}}_{\widehat V^{1}_{\pol}\left(\widetilde{\mr{\Omega}}_{C}\right)}
		\leq C \norm{\breve{v}}_{V^{1}\left(\widetilde{\breve{\Omega}}_C\right)}	
		\leq C \norm{\breve{v}}_{V^{s-\beta}\left(\widetilde{\breve{\Omega}}_C\right)}
		\leq C \norm{\breve{v}}_{V^s_{\beta}\left(\widetilde{\breve{\Omega}}_C\right)}
		\leq C \norm{\mr{v}}_{\widehat{V}^{s}_{\pol,\beta}\left(\widetilde{\mr{\Omega}}_{C}\right)}.
	\end{align}
	By the scaling argument \eqref{eq: scaling argument r-derivative} and the estimates \eqref{eq: stability of projector weighted}, \eqref{eq: polar spaces and 1/2-spaces} and \eqref{eq: embedding weighted}, it follows
	\begin{align*}
		\norm{\partial^q_r\left( W\widehat{v} - \Pi^{\pol}_{\bs p, \bs \Xi}(W\widehat{v})\right)}_{\widehat L^2_{1/2}(\widehat{\Omega}_C)}
		&\leq C h_{r, \widehat{\Omega}_C} ^{1-q}  \norm{\partial^q_{\mr r} \left(\mr{v} -\mr{\Pi}^{\pol}_{\bs p, \bs \Xi} \mr{v}\right)}_{\widehat L^2_{1/2}(\mr\Omega_C)} \\
		&  \leq C h_{r, \widehat{\Omega}_C} ^{1-q} \left(\norm{\mr{v}}_{\widehat H^{1}_{1/2}(\mr\Omega_C)} + \norm{\partial^q_{\mr r} \left( \mr{\Pi}^{\pol}_{\bs p, \bs \Xi}\mr{v} \right)}_{\widehat L^2_{1/2}(\mr\Omega_C)} \right)  \\
		&  \leq C h_{r, \widehat{\Omega}_C} ^{1-q} \left(\norm{\mr{v}}_{\widehat V^{1}_{\pol}\left(\widetilde{\mr{\Omega}}_{C}\right)} + \norm{ \mr{v}}_{\widehat{V}^{s}_{\pol,\beta}\left(\widetilde{\mr{\Omega}}_{C}\right)} \right)  \\
		&\leq C h_{r, \widehat{\Omega}_C} ^{1-q} \norm{ \mr{v}}_{\widehat{V}^{s}_{\pol,\beta}\left(\widetilde{\mr{\Omega}}_{C}\right)} 
		\leq C h_{r, \widehat{\Omega}_C}^{s-q-\beta} \norm{\widehat{v}}_{\mathcal{V}^{s}_{\pol,\beta}\left(\widetilde{\widehat{\Omega}}_{C}\right)} .
	\end{align*}
	Note that the $V$-spaces scale perfectly, which was exploited in the last step in combination with the Leibniz rule. Next, we adapt the proof of Lemma \ref{lemma: estimate of phi-derivative}. Let $\widehat v_0$ be the auxiliary function defined there such that $\widehat v_0(0,\cdot) = 0$ and let $\mr v_0 = \mr W \mr v_0$ be the corresponding scaling, recall \eqref{eq: scaling transformation}.
	By combining Lemmata \ref{lemma: derivative of projector} and \ref{lemma: estimate of derivative of projector}, the definition of polar weighted Sobolev norms and \eqref{eq: embedding weighted}, we obtain
	\begin{align*}
		\norm{\partial_\varphi\left(W\widehat{v} - \Pi^{\pol}_{\bs p, \bs \Xi}( W \widehat{v})\right)}_{\widehat 	L^2_{-1/2}(\widehat{\Omega}_C)}
		&= \norm{\partial_{\mr \varphi}\left( \mr{v}_0 - \mr{\Pi}^{\pol}_{\bs p, \bs \Xi} \mr{v}_0\right)}_{\widehat 	L^2_{-1/2}(\mr \Omega_C)}
		= \norm{\partial_{\mr \varphi} \mr{v}_0 - \mr{\Pi}^{\pol,\partial_\varphi}_{(p_1,p_2-1), \bs \Xi}	(\partial_{\varphi}\mr{v}_0)}_{\widehat L^2_{-1/2}(\mr \Omega_C)} \\
		& \leq C  \norm{ \partial_{\mr \varphi} \mr{v}_0}_{\widehat L^2_{-1/2}\left(\widetilde{\mr{\Omega}}_{C}\right)} 
		\leq C \norm{\mr v_0}_{\widehat V^1_{\pol}\left(\widetilde{\mr{\Omega}}_C\right)} 
		\leq C  \norm{\mr{v}_0}_{\widehat{V}^{s}_{\pol,\beta}\left(\widetilde{\mr{\Omega}}_{C}\right)} \\
		& \leq C h_{r, \widehat{\Omega}_C}^{s-1-\beta} 	\norm{\widehat{v}}_{\mathcal{V}^{s}_{\pol,\beta}\left(\widetilde{\widehat{\Omega}}_{C}\right)} ,
	\end{align*}
	where, in the last step, we used again the perfect scaling of $V$-spaces.
	
	Finally, the proof of Lemma \ref{lemma: estimate of L_{-1/2}-norm} is modified analogously. Similar to Lemma \ref{lemma: estimate of derivative of projector}, stability of the projector $\mr{\Pi}^{\pol}_{\bs p, \bs \Xi}$ in the $L^2_{-1/2}(\mr \Omega_C)$-norm can be shown, and it follows
	\begin{align*}
		\norm{\mr{v}_0 - \mr{\Pi}^{\pol}_{\bs p, \bs \Xi} \mr{v}_0}_{\widehat L^2_{-1/2}(\mr \Omega_C)} 
		\leq C \norm{\mr{v}_0}_{\widehat L^2_{-1/2}\left(\widetilde{\mr{\Omega}}_C\right)} 
		\leq C \norm{\mr{v}_0}_{\widehat V^1_{\pol}\left(\widetilde{\mr{\Omega}}_C\right)} 
		\leq C  \norm{\mr{v}_0}_{\widehat{V}^{s}_{\pol,\beta}\left(\widetilde{\mr{\Omega}}_{C}\right)}
		\leq C h_{r, \widehat{\Omega}_C}^{s-1-\beta} 	\norm{\widehat{v}}_{\mathcal{V}^{s}_{\pol,\beta}\left(\widetilde{\widehat{\Omega}}_{C}\right)} .
	\end{align*}
	
	Hence, we have estimated all four terms in \eqref{eq: error components H^1 projection}. By the transformation \eqref{eq: transformation polar bent Sobolev spaces} of the weighted polar Sobolev norms back to the physical domain, we obtain
	\begin{align*}
		\abs{v -  \Pi_{V_h^{\pol}} v }_{H^q(\Omega_C)}
		\leq C h_{r, \widehat{\Omega}_C}^{s-q-\beta} \norm{\widehat{v}}_{\mathcal{V}_{\pol,\beta}^{s}\left(\widetilde{\widehat{\Omega}}_C\right)} 
		\leq C h_{r, \widehat{\Omega}_C}^{s-q-\beta} \norm{v}_{V^{s}_{\beta}(\widetilde \Omega_C)} ,
	\end{align*}	
	and the demonstration is complete.
\end{proof}

\subsubsection{Estimates on the full domain $\Omega$}
Now, we are in the position to show an overall projection error estimate on polar domains with corners and prove optimal convergence for suitable grading parameters.

\begin{Theorem}
	\label{theorem: projection error estimate}
	Let $q \in \{0,1\}$ and $s,s_0 \in \N$ with $2 \leq s \leq s_0 \leq p +1$. Further, let $v \in V^{s}_{\beta}(\Omega)$ for all $s-1>\beta >s-1-\nu$. If the mesh grading parameter satisfies condition \eqref{eq: mesh grading parameter condition}, we obtain optimal convergence of the projection with respect to the mesh size $h$,
	\begin{align*}
		\norm{v - \Pi_{V_h^{\pol}}v}_{H^q(\Omega)} \leq C h^{s-q} \norm{v}_{V^{s}_{\beta}(\Omega)}.
	\end{align*}
\end{Theorem}

\begin{proof}
	By the splitting of the model domain defined in Section \ref{subsec: splitting of the model domain}, we obtain
	\begin{align*}
		\norm{v-\Pi_{V_h^{\pol}}v}_{H^q(\Omega)}^2
		= \norm{v-\Pi_{V_h^{\pol}}v}_{H^q(\Omega_C)}^2 + \norm{v-\Pi_{V_h^{\pol}}v}_{H^q(\Omega_R)}^2 .
	\end{align*}
	The second term has been estimated in Theorem \ref{theorem: projection error ring domain}. For the first term, we use Theorem \ref{theorem: local projection error polar stripe weighted}, the mesh properties \eqref{eq: edge lengths parametric support extension}, where $h_{r,\Omega_C} \sim \widetilde{h}_{1,j_1}$, and set $\beta = (1-\mu)(s-q)$ such that \eqref{eq: beta property} is satisfied. It follows
	\begin{align*}
		\norm{v-\Pi_{V_h^{\pol}}v}_{H^q(\Omega_C)}^2  
		&\leq C h_{r,\Omega_C}^{2(s - q- \beta)} \norm{v}_{V^{s}_{\beta}(\widetilde{\Omega}_C)} 
		= C h_{r,\Omega_C}^{2(s - q)(1-1+\mu)} \norm{v}_{V^{s}_{\beta}(\widetilde{\Omega}_C)} \\
		&\leq C h_1^{2(s - q)\mu/\mu} \norm{v}_{V^{s}_{\beta}(\widetilde{\Omega}_C)} 
		= C h_1^{2(s-q)} \norm{v}_{V^{s}_{\beta}(\widetilde{\Omega}_C)}^2 .
	\end{align*}
	In total, for comparable refinement in both directions, that is, $h_1 \sim h_2$ and $h=\max\{h_1,h_2\}$, we have
	\begin{align*}
		\norm{v-\Pi^{\pol}_{V_h}v}_{H^q(\Omega)}^2
		\leq C h^{2(s-q)} \left(\norm{v}_{V^{s}_{\beta}(\widetilde{\Omega}_C)}^2 + \norm{v}_{V^{s}_{\beta}(\widetilde{\Omega}_R)}^2 \right)
		\leq C  h^{2(s-q)} \norm{v}_{V^{s}_{\beta}(\Omega)}^2 ,
	\end{align*}
	where we use that $\widetilde{\Omega}_C$ and $\widetilde{\Omega}_R$ have bounded intersection with respect to $h$, and the assertion is proven.	
\end{proof}

The main results of our paper, in particular Theorem \ref{theorem: main result}, follow immediately from Theorem \ref{theorem: projection error estimate}. Finally, we point out the sub-optimal convergence of the projection on uniform meshes.

\begin{remark}[Convergence on uniform meshes]
	The proof of Theorem \ref{theorem: projection error estimate} shows that on uniform meshes, the convergence rate for the $H^q$-error of the projection is bounded by
	\begin{align*}
		s-q-\beta<s-q-(s-1-\nu) = \nu + 1 - q .
	\end{align*}
	The approximation error in the $H^1$-norm decreases with the same order as the projection error. In contrast, this is not the case in the $L^2$-norm. By using a standard duality argument, it can be proven that the order on uniform meshes is given by $2\nu$.
	\label{remark: convergence on uniform meshes}
\end{remark}

\section{Numerical results}
\label{section: Numerical results}
In this section, we confirm our theoretical findings with numerical results. All experiments are carried out using the computing package GeoPDEs 3.0 \cite{Vazquez2016,deFalcoRealiVazquez2011}. In particular, we apply the graded $h$-refinement method introduced in Section \ref{subsec: construction of graded mesh}, which can be implemented easily in GeoPDEs by adapting one line of the standard knot refining routine \textit{kntrefine} provided in the package, to  validate the optimal convergence orders that have been established theoretically in Section \ref{section: Error estimates}.

All integrals are computed using a Gauss-Legendre quadrature rule with $36$ Gauss points per mesh element, where $6$ quadrature nodes are employed in both univariate directions. 


\begin{example}[\textbf{Pacman domain}]
	In the first example, we solve the Poisson equation \eqref{eq: poisson equation} on the Pacman domain, which is defined as a circular sector with angle $\omega = \frac{5}{3}\pi$, i.e., $\Omega=\{(r \cos \varphi, r \sin \varphi) : r \in (0,1), \varphi \in (0, \omega)\}$, with homogeneous Dirichlet boundary conditions on the circular edge, $\Gamma_D=\{(\cos \varphi, \sin \varphi) : \varphi \in [0,\omega]\}$, and homogeneous Neumann boundary conditions on the re-entrant edges, $\Gamma_N = \{(r \cos \varphi, r \sin \varphi) : 0\leq r\leq 1, \varphi \in \{0,\omega\}\}$. With the right-hand side $ f: \Omega \to \R, \ f(x,y) = \widehat f \circ \bs F^{-1}(x,y) = \widehat f(r,\varphi) = -r^{\nu-1} \cos(\nu\varphi)(-2\nu-1)$, we obtain the manufactured solution 
	\begin{align*}
		u: \Omega \to \R, \quad u(x,y) = \widehat u(r,\varphi) = r^\nu \cos(\nu \varphi) (1-r) .
	\end{align*}
	In Figure \ref{fig: Poisson Pacman - exact solution contour}, we illustrate the exact solution using a contour plot. Figures \ref{fig: Poisson Pacman - numerical solution mesh contour} and \ref{fig: Poisson Pacman - numerical solution mesh surf} show the numerical solution obtained with a polynomial degree of $p=3$ on a coarse graded mesh. As explained in Section \ref{subsec: model problems and regularity properties}, the solution has a singularity of type $r^{\nu}$ with $\nu = \frac{\pi}{\omega} = \frac35$ and satisfies $u \in V^{s}_{\beta}(\Omega)$ for all $s\in \N$ and $s-1> \beta > s- 1 -\nu = s - \frac85$.
	
	The approximation errors in the $H^1(\Omega)$-norm and $L^2(\Omega)$-norm are plotted for decreasing mesh size $h$ in Figures \ref{fig: Poisson Pacman - H1-error} and \ref{fig: Poisson Pacman - L2-error}. For each NURBS degree $p$, we consider $C^{p-1}$-smooth NURBS on uniform and graded meshes, corresponding to grading parameters $\mu =1$ and $\mu = 0.9 \frac{\nu}{p}$, respectively, recall \eqref{eq: choice of grading parameter}. The experimental convergence rates on uniform meshes are observed to be $\nu$ and $2\nu$ independent of the polynomial degree, as expected from Remark \ref{remark: convergence on uniform meshes}. In contrast, graded meshes produce optimal rates of $p$ and $p+1$ in the $H^1(\Omega)$-norm and $L^2(\Omega)$-norm, respectively, in agreement with our main result stated in Corollary \ref{cor: approximation error}.
	\begin{figure}[t]
		\begin{subfigure}{0.33 \linewidth}
			\begin{center}
				\includegraphics[width= \linewidth,  trim=1cm 1cm 1cm 1cm, clip]{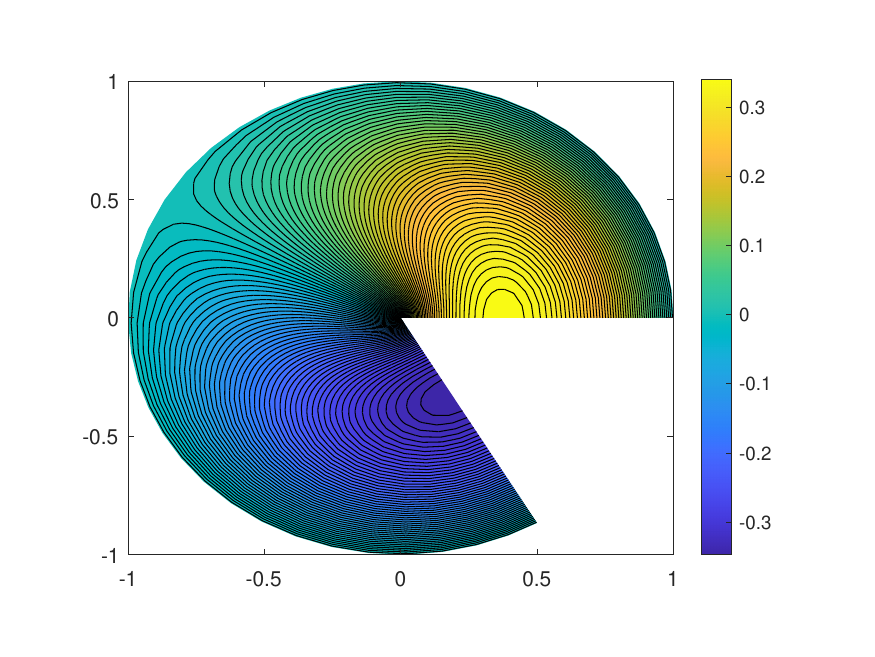}
			\end{center}	
			\caption{}	
			\label{fig: Poisson Pacman - exact solution contour}
		\end{subfigure}
		\hspace*{1mm}
		\begin{subfigure}{0.33 \linewidth}
			\includegraphics[width= \linewidth, trim=1cm 0.1cm 1cm 1cm, clip]{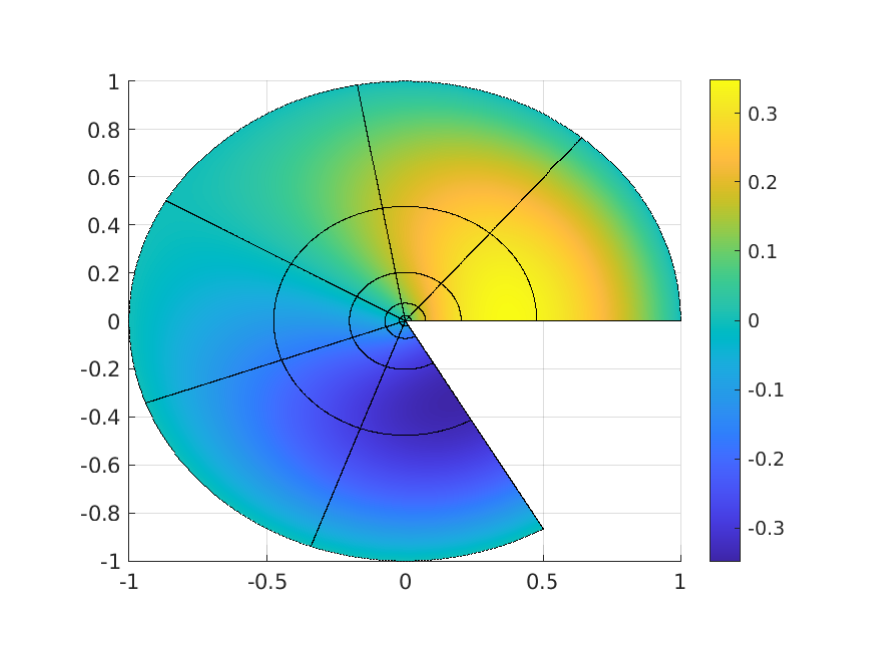}
			\caption{}	
			\label{fig: Poisson Pacman - numerical solution mesh contour}
		\end{subfigure}
		\begin{subfigure}{0.33 \linewidth}
			\includegraphics[width= \linewidth, trim=1cm 0cm 1cm 1.2cm, clip]{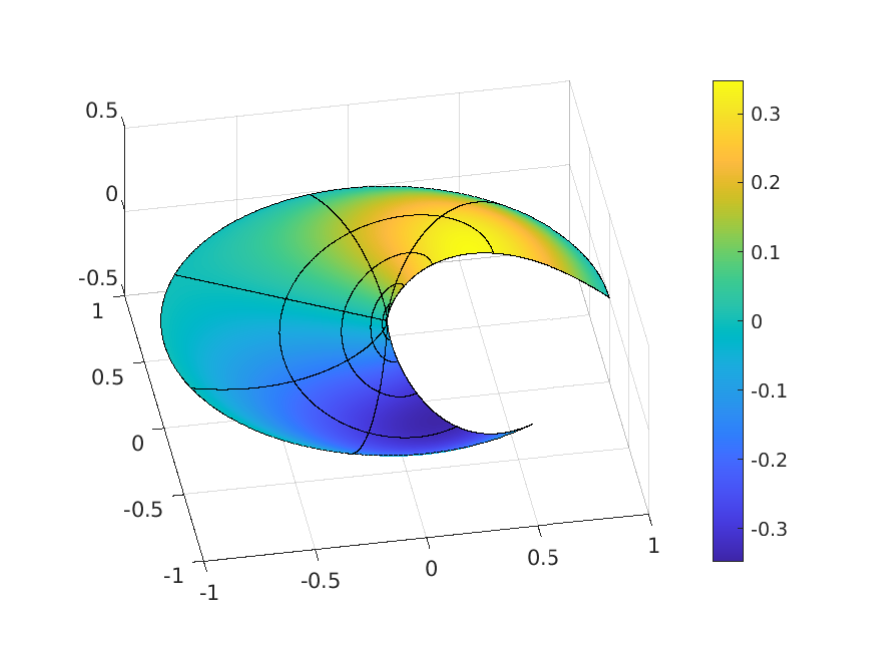}
			\caption{}	
			\label{fig: Poisson Pacman - numerical solution mesh surf}
		\end{subfigure}
		\caption{Solving the Poisson equation on the Pacman domain: (a): Contour plot of the exact solution. (b): Contour plot of an approximation with cubic NURBS on a coarse graded mesh. (c): Surface plot of an approximation with cubic NURBS on a coarse graded mesh.}
		\label{fig: Poisson Pacman - solution plots}
	\end{figure}
	
	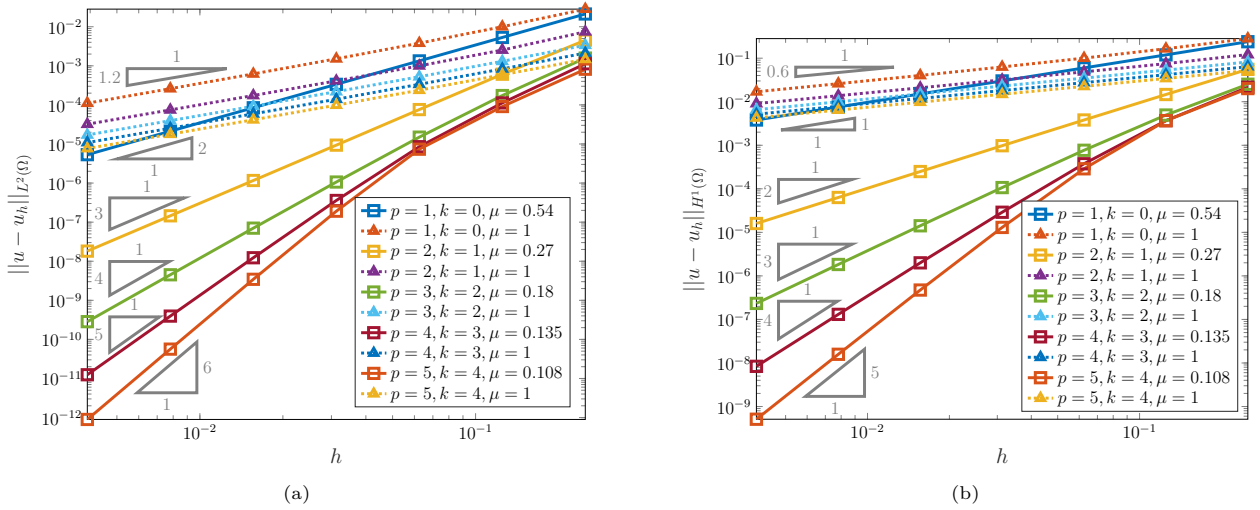
\begin{figure}
		\begin{subfigure}{0.47 \linewidth}
			\resizebox{\linewidth}{!}{\input{laplace_unitdisk300_mg_neumann_hom_H1.tikz}}
			\caption{}	
			\label{fig: Poisson Pacman - H1-error}
		\end{subfigure}%
		\hfill
		\begin{subfigure}{0.46 \linewidth}
			\resizebox{\linewidth}{!}{\input{laplace_unitdisk300_mg_neumann_hom_L2.tikz}}
			\caption{}	
			\label{fig: Poisson Pacman - L2-error}
		\end{subfigure}
		\caption{Solving the Poisson equation on the Pacman domain: (a): $H^1(\Omega)$-error of the approximation. (b): $L^2(\Omega)$-error of the approximation.}
		\label{fig: Poisson Pacman - errors}
	\end{figure}
\end{example}

\begin{example}[\textbf{L-shaped domain}]
	In the next numerical test, we solve the Poisson equation \eqref{eq: poisson equation} on an L-shaped domain, which has been described in Example \ref{ex: L-Shape}. We impose mixed boundary conditions at the corner, that is, homogeneous Neumann boundary conditions on one re-entrant edge, $\Gamma_N = \{(0,-r)^T : r \in[0,1]\}$, and homogeneous Dirichlet boundary conditions elsewhere, $\Gamma_D = \Gamma \setminus \Gamma_N$. The exact solution is illustrated in Figure \ref{fig: Poisson on L-shape with mixed bc - solution} using a surface plot, where some mesh lines are inserted for visualization purposes. As discussed in Section \ref{subsec: model problems and regularity properties}, the solution has a singularity of type $r^{\nu}$ with $\nu = \frac{\pi}{2\omega}=\frac13$, where $\omega=\frac32 \pi$ is the corner angle. With the right-hand side 
	\begin{align*}
		f: \Omega \to \R, \quad f(x,y) = \widehat f(r,\varphi) = r^{\nu} \left( (4 + 4\nu - 2r^2 - \nu r^2) \sin(\nu \varphi) 
		+ \nu r^2 \sin(\varphi (\nu - 4))\right),
	\end{align*}
	we obtain the manufactured solution 
	\begin{align*}
		u: \Omega \to \R, \quad u(x,y) = \widehat u(r,\varphi) = r^\nu \sin(\nu \varphi) (1-r^2 + r^4 \cos^2(\varphi) \sin^2(\varphi)) ,
	\end{align*}
	which satisfies $u \in V^{s}_{\beta}(\Omega)$ for all $s \in \N$ and $s-1> \beta > s- 1 -\nu = s - \frac43$. Note that, due to the relation $(1-x)(1+x)(1-y)(1+y) = 1-r^2 + r^4 \cos^2(\varphi) \sin^2(\varphi)$, it holds $u = 0$ on $\Gamma_{D}$.
	
	The approximation errors in the $H^1(\Omega)$-norm and $L^2(\Omega)$-norm are plotted for decreasing mesh size $h$ in Figures \ref{fig: Poisson on L-shape with mixed bc - H1-error} and \ref{fig: Poisson on L-shape with mixed bc - L2-error}. For each degree $p$, we consider $C^{p-1}$-smooth NURBS on uniform and graded meshes, corresponding to grading parameters $\mu =1$ and $\mu = 0.9 \frac{\nu}{p}$, respectively. Similar to the previous example, the experimental convergence rates on uniform meshes are observed to be $\nu$ and $2\nu$ independent of the polynomial degree, whereas graded meshes produce optimal convergence orders of $p$ and $p+1$ in the $H^1(\Omega)$-norm and $L^2(\Omega)$-norm, respectively.
	\begin{figure}[t]
		\begin{subfigure}{0.29 \linewidth}
			\begin{center}
				\includegraphics[width= \linewidth, trim=1cm 0.8cm 1cm 1.2cm, clip]{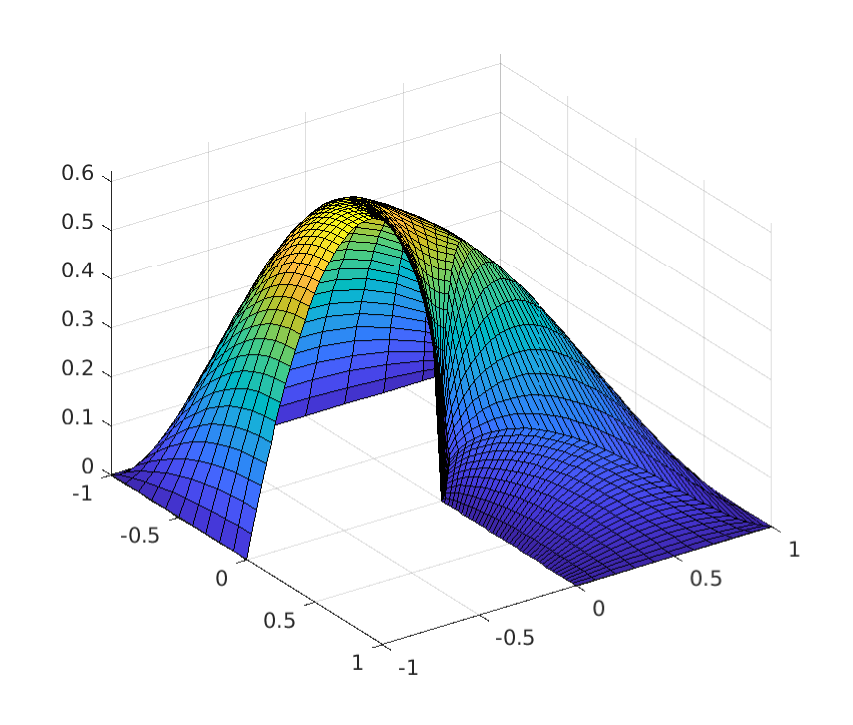}
			\end{center}	
			\caption{}	
			\label{fig: Poisson on L-shape with mixed bc - solution}
		\end{subfigure}
		\begin{subfigure}{0.35 \linewidth}
			\resizebox{\linewidth}{!}{\input{laplace_Lshaped_mixed_hom_p=1to5_n=6_H1_err_h_slopes.tikz}}
			\caption{}	
			\label{fig: Poisson on L-shape with mixed bc - H1-error}
		\end{subfigure}
		\hfill
		\begin{subfigure}{0.34 \linewidth}
			\resizebox{\linewidth}{!}{\input{laplace_Lshaped_mixed_hom_p=1to5_n=6_L2_err_h_slopes.tikz}}
			\caption{}	
			\label{fig: Poisson on L-shape with mixed bc - L2-error}
		\end{subfigure}
		\caption{Solving the Poisson equation on the L-shaped domain: (a): Plot of the exact solution. (b): $H^1(\Omega)$-error of the approximation. (c): $L^2(\Omega)$-error of the approximation.}
		\label{fig: Poisson on L-shape with mixed bc}
	\end{figure}
\end{example}

\begin{example}[\textbf{Transverse vibrations of a circular drum with crack}]
	The Laplace eigenvalue problem on circular sectors, for which exact solutions are known, has been studied in detail in our recent work \cite{ApelZilk2024}. The numerical results in Section 5.1 of this paper serve as a further verification of our theory. Corollary \ref{cor: approximation error} can be extended in a standard way to the setting of eigenvalue problems.
\end{example}

\section{Conclusion \& Outlook}
\label{section: Conclusion}
In this paper, we analyzed the approximation power of splines on polar domains with corners. To address the lack of regularity in polar parameterizations, which standard isogeometric approximation theory does not accommodate, we introduced polar Sobolev spaces on the parametric domain, derived corresponding error estimates, and mapped them to the polar domain. Furthermore, we proved optimal convergence rates for the approximation of singular solutions on graded meshes with suitable grading parameters. With that, we established a new local refinement approach in IGA, which, unlike comparable methods, does not require to break the tensor-product structure. Additionally, we demonstrated that our grading scheme is easy to implement and produces excellent numerical results.

This work serves as a foundation for further research in various directions. First, the approach may be extended to polar domains without corner, such as the unit disk, which has many practical applications, e.g., in mathematical biology \cite{LaguneroFellnerApelKempfZilk2025}. In addition, it may be generalized for other model problems and equations, and the treatment of inhomogeneous Dirichlet boundary conditions requires further investigation. In that context, for higher-order PDEs, such as the biharmonic equation, an extension of our approach to $C^k$-smooth polar splines with $k\geq 1$ \cite{ToshniwalSpeleersHiemstraHughes2017,SpeleersToshniwal2021} is required.

Moreover, the proposed graded refinement can be combined with a hierarchical scheme to avoid anisotropic elements near the corner \cite{ApelZilk2024}, and corresponding error estimates can be derived. In particular, further analysis is needed to understand why smooth solutions are not approximated optimally on non-graded hierarchical meshes. Other than that, developing error estimates for interpolation operators as an alternative to the projector used in this paper could potentially yield new insights. Finally, extensions to more complex domains and multipatch formulations should be explored.

\section*{Acknowledgments}
The authors thank Johannes Pfefferer and Max Winkler for sharing their expertise on the concepts of mesh grading for finite elements and weighted Sobolev spaces. Further, the authors express their gratitude to Giancarlo Sangalli and Thomas Takacs for motivating scientific discussions, which helped to place the results of this paper within the existing literature on isogeometric approximation theory, and to Rafael Vázquez for providing an open-source software package for IGA and for his prompt assistance with any related questions.

\appendix

\renewcommand{\thesection}{\Alph{section}} 
\makeatletter
\renewcommand\@seccntformat[1]{\appendixname\ \csname the#1\endcsname.\hspace{0.5em}}
\makeatother
 
\section{Properties of the modified projector}
The following technical computations, which are required for the proofs of Lemma \ref{lemma: estimate of phi-derivative} and Theorem \ref{theorem: local projection error polar stripe weighted}, are included in the appendix to improve the flow of the main text. In particular, we establish key properties of the modified projector \eqref{eq: modified projector physical domain}.

\begin{lemma}
	Let the assumptions of Theorem \ref{theorem: local projection error polar stripe} hold. Then, there is a projector
	\begin{align*}
		\mr{\Pi}^{\pol,\partial_\varphi}_{(p_1,p_2-1), \bs \Xi}: C(\mr \Omega) \to S_{(p_1,p_2-1)}( \bs \Xi),
	\end{align*}
	with $\mr \Omega$ as in \eqref{eq: scaling transformation}, such that
	\begin{align*}
		\partial_{\mr \varphi}\left( \mr{\Pi}^{\pol}_{\bs p, \bs \Xi} \mr{v}_0\right)
		= \mr{\Pi}^{\pol,\partial_\varphi}_{(p_1,p_2-1), \bs \Xi} (\partial_{\mr \varphi} \mr{v}_0 ) 
	\end{align*}
	for all $\mr v_0 \in  C(\mr \Omega)$ with $\mr v_0(0,\cdot) = 0$.
	\label{lemma: derivative of projector}
\end{lemma}

\begin{proof}
	By representation \eqref{eq: representation in basis} and the formula for differentiation of B-splines \cite[Section 1.2.4]{LycheManniSpeleers2018}, it follows for $\mr v_0 \in  C(\mr \Omega)$ with $\mr v_0(0,\cdot) = 0$ that
	\allowdisplaybreaks
	\begin{align*}
		&\partial_{\mr \varphi}\left( \mr{\Pi}^{\pol}_{\bs p, \bs \Xi} \mr{v}_0\right) (\mr r, \mr \varphi)
		= \partial_{\mr \varphi}\left(\sum_{\bs i \in \bs I \setminus \bs I_{\pol}} {\mr \lambda}_{\bs i, \bs p}(\mr{v}_0) \mr{B}_{\bs i,\bs p} (\mr r, \mr \varphi) \right) 
		= \sum_{\bs i \in \bs I \setminus \bs I_{\pol}} {\mr \lambda}_{\bs i, \bs p}(\mr{v}_0) \, \partial_{\mr \varphi} \mr{B}_{\bs i,\bs p} (\mr r, \mr \varphi)  \\
		&= \sum_{i_1=2}^{n_1} \sum_{i_2=1}^{n_2} (\mr \lambda_{i_1, p_1} \otimes \mr \lambda_{i_2, p_2})(\mr{v}_0) \, \mr{B}_{i_1,p_1}(\mr r) \, \partial_{\mr \varphi} \mr{B}_{i_2,p_2}(\mr \varphi)   \\
		&= \sum_{i_1=2}^{n_1} \sum_{i_2=1}^{n_2} (\mr \lambda_{i_1, p_1} \otimes \mr \lambda_{i_2, p_2})(\mr{v}_0) \, \mr{B}_{i_1,p_1}(\mr r) \, p_2 \left( \frac{\mr B_{i_2,p_2-1}(\mr \varphi)}{\mr \xi_{i_2+p_2,2} - \mr \xi_{i_2,2}} - \frac{\mr B_{i_2+1,p_2-1}(\mr \varphi)}{\mr \xi_{i_2+p_2+1,2} - \mr  \xi_{i_2+1,2}}\right) \\
		&= p_2 \sum_{i_1=2}^{n_1} \sum_{i_2=1}^{n_2} \mr \lambda_{i_1, p_1} \left(  \int_{\mr \xi_{i_2,2}}^{\mr \xi_{i_2+p_2+1,2}} \mr{v}_0(\cdot, t) \, \mathring{D}^{p_2+1} \mr \psi_{i_2}(t) \, \Dd t \right)  \mr{B}_{i_1,p_1}(\mr r) \left( \frac{\mr B_{i_2,p_2-1}(\mr \varphi)}{\mr \xi_{i_2+p_2,2} - \mr \xi_{i_2,2}} - \frac{\mr B_{i_2+1,p_2-1}(\mr \varphi)}{\mr \xi_{i_2+p_2+1,2} - \mr \xi_{i_2+1,2}}\right) \\
		&= p_2 \sum_{i_1=2}^{n_1} \mr \lambda_{i_1, p_1} \left(\sum_{i_2=1}^{n_2} \int_{\mr \xi_{i_2,2}}^{\mr \xi_{i_2+p_2+1,2}} \mr{v}_0(\cdot, t) \, \mathring{D}^{p_2+1} \mr \psi_{i_2}(t) \, \Dd t \left( \frac{\mr B_{i_2,p_2-1}(\mr \varphi)}{\mr \xi_{i_2+p_2,2} - \mr \xi_{i_2,2}} - \frac{\mr B_{i_2+1,p_2-1}(\mr \varphi)}{\mr \xi_{i_2+p_2+1,2} - \mr \xi_{i_2+1,2}}\right) \right)  \mr{B}_{i_1,p_1}(\mr r) .
	\end{align*}
	By reordering the sum that appears as an argument of the dual functional $\mr \lambda_{i_1,p_1}$, we obtain
	\begin{align*}
		&\sum_{i_2=1}^{n_2} \int_{\mr \xi_{i_2,2}}^{\mr \xi_{i_2+p_2+1,2}} \mr{v}_0(\cdot, t) \, \mathring{D}^{p_2+1} \mr \psi_{i_2}(t) \, \Dd t \left( \frac{\mr B_{i_2,p_2-1}(\mr \varphi)}{\mr \xi_{i_2+p_2,2} - \mr \xi_{i_2,2}} - \frac{\mr B_{i_2+1,p_2-1}(\mr \varphi)}{\mr \xi_{i_2+p_2+1,2} - \mr \xi_{i_2+1,2}} \right) \\
		&= \int_{\mr \xi_{1,2}}^{\mr \xi_{1+p_2+1,2}} \mr{v}_0(s,t) \, \mathring{D}^{p_2+1} \mr \psi_{1}(t) \, \dt \, \frac{\mr B_{1,p_2-1}(\mr \varphi)}{\mr \xi_{1+p_2,2} - \mr \xi_{1,2}} \\
		&\quad +\sum_{i_2=2}^{n_2} \left( \int_{\mr \xi_{i_2,2}}^{\mr \xi_{i_2+p_2+1,2}} \mr{v}_0(s,t) \, \mathring{D}^{p_2+1} \mr \psi_{i_2}(t) \, \dt - \int_{\mr \xi_{i_2-1,2}}^{\mr \xi_{i_2+p_2,2}} \mr{v}_0(s,t) \, \mathring{D}^{p_2+1} \mr \psi_{i_2-1}(t) \, \dt \right)  \frac{\mr B_{i_2,p_2-1}(\mr \varphi)}{\mr \xi_{i_2+p_2,2} - \mr \xi_{i_2,2}} \\
		&\quad - \int_{\mr \xi_{n_2,2}}^{\mr \xi_{n_2+p_2+1,2}} \mr{v}_0(s,t) \, \mathring{D}^{p_2+1} \mr \psi_{n_2}(t) \,\dt \, \frac{\mr B_{n_2+1,p_2-1}(\mr \varphi)}{\mr \xi_{n_2+p_2+1,2} - \mr \xi_{n_2+1,2}} \\
		&= \sum_{i_2=2}^{n_2} \left( \int_{\mr \xi_{i_2,2}}^{\mr \xi_{i_2+p_2+1,2}} \mr{v}_0(s,t) \, \mathring{D}^{p_2+1} \mr \psi_{i_2}(t) \, \dt - \int_{\mr \xi_{i_2-1,2}}^{\mr \xi_{i_2+p_2,2}} \mr{v}_0(s,t) \, \mathring{D}^{p_2+1} \mr \psi_{i_2-1}(t) \, \dt \right)  \frac{\mr B_{i_2,p_2-1}(\mr \varphi)}{\mr \xi_{i_2+p_2,2} - \mr \xi_{i_2,2}} .
	\end{align*}	
	The first and last term of the reordered sum vanish since $\Xi_2$ is an open knot vector and fractions with zero denominator are defined to have value zero. Next, we rewrite the term in brackets,
	\begin{align}
		& \int_{\mr \xi_{i_2,2}}^{\mr \xi_{i_2+p_2+1,2}} \mr{v}_0(s,t) \, \mathring{D}^{p_2+1} \mr \psi_{i_2}(t) \, \dt - \int_{\mr \xi_{i_2-1,2}}^{\mr \xi_{i_2+p_2,2}} \mr{v}_0(s,t) \, \mathring{D}^{p_2+1} \mr \psi_{i_2-1}(t) \, \dt \nonumber \\
		&= \int_{\mr \xi_{i_2,2}}^{\mr \xi_{i_2+p_2+1,2}} \mr{v}_0(s,t) \, \mathring{D}^{p_2+1} \mr \psi_{i_2}(t) \, \dt - \int_{\mr \xi_{i_2,2}}^{\mr \xi_{i_2+p_2+1,2}} \mr{v}_0(s,\alpha(t)) \, \mathring{D}^{p_2+1}( \mr \psi_{i_2-1} \circ \alpha(t)) \abs{\alpha'(t)} \, \dt \label{eq: new dual functional two parts} ,
	\end{align}
	where we use the transformation $\alpha(t) := t - h_{\varphi}$ based on the mesh size in angular direction $h_{\varphi}$. With definitions \eqref{eq: dual functional definition of Phi} and \eqref{eq: dual functional definition of G}, we further compute
	\begin{align*}
		\psi_{i_2-1}(t- h_{\varphi}) 
		&= G_{i_2-1}(t - h_{\varphi}) \Phi_{i_2-1} (t- h_{\varphi}) \\
		&= g\left(\frac{2 (t- h_{\varphi}) - \xi_{i_2-1,2} - \xi_{i_2+p_2,2}}{\xi_{i_2+p_2,2}- \xi_{i_2-1,2}}\right) \left(\frac{(t- h_{\varphi}-\xi_{i_2,2}) \cdots (t- h_{\varphi} - \xi_{i_2+p_2-1,2})}{p_2!}\right) \\
		&= g\left(\frac{2 t - (\xi_{i_2-1,2} + h_{\varphi}) - (\xi_{i_2+p_2,2}+h_{\varphi})}{\xi_{i_2+p_2+1,2} - h_{\varphi}- (\xi_{i_2,2} - h_{\varphi})}\right) \left(\frac{(t-(\xi_{i_2,2}+h_{\varphi})) \cdots (t - (\xi_{i_2+p_2-1,2}+h_{\varphi}))}{p_2!}\right) \\
		&= g\left(\frac{2 t - \xi_{i_2,2} - \xi_{i_2+p_2+1,2}}{\xi_{i_2+p_2+1,2}- \xi_{i_2,2}}\right) \left(\frac{(t-\xi_{i_2+1,2}) \cdots (t - \xi_{i_2+p_2,2}))}{p_2!} \right) \\
		& = G_{i_2}(t) \Phi_{i_2} (t) \\
		& = \psi_{i_2}(t) .
	\end{align*}	
	It follows that the expression \eqref{eq: new dual functional two parts} can be written in terms of a new dual functional,
	\begin{align*}
		\mu_{i_2,p_2-1} (\mr{w}) := \int_{\mr \xi_{i_2,2}}^{\mr \xi_{i_2+p_2+1,2}} \int_{\alpha(t)}^{t} \mr{w}(z) \Dd z \, \mathring{D}^{p_2+1} \mr \psi_{i_2}(t) \, \dt ,
	\end{align*}
	with $ \mr w(\cdot) = \partial_{\mr \varphi}\mr{v}_0(s,\cdot)$. By combining all results, we obtain
	\begin{align*}
		\partial_{\mr \varphi}\left( \mr{\Pi}^{\pol}_{\bs p, \bs \Xi} \mr{v}_0\right) 
		&= p_2 \sum_{i_1=2}^{n_1} \lambda_{i_1,p_1} \left(\sum_{i_2=2}^{n_2} \mu_{i_2,p_2-1} (\partial_{\mr \varphi}\mr{v}_0(s,\cdot)) \frac{\mr B_{i_2,p_2-1}(\mr \varphi)}{\mr \xi_{i_2+p_2,2} - \mr \xi_{i_2,2}}  \right) \mr{B}_{i_1,p_1}(\mr r) \\
		&= p_2 \sum_{i_1=2}^{n_1} \sum_{i_2=2}^{n_2}  \left(\lambda_{i_1,p_1} \otimes \mu_{i_2,p_2-1}\right) (\partial_{\mr \varphi}\mr{v}_0) \mr{B}_{i_1,p_1}(\mr r)  \frac{\mr B_{i_2,p_2-1}(\mr \varphi)}{\mr \xi_{i_2+p_2,2} - \mr \xi_{i_2,2}} \\
		&= \mr{\Pi}^{\pol,\partial_\varphi}_{(p_1,p_2-1), \bs \Xi} (\partial_{\mr \varphi} \mr{v}_0 )  ,
	\end{align*}
	where
	\begin{align}
		\mr{\Pi}^{\pol,\partial_\varphi}_{(p_1,p_2-1), \bs \Xi}(\mr w):= p_2 \sum_{i_1=2}^{n_1} \sum_{i_2=2}^{n_2}  \left(\lambda_{i_1,p_1} \otimes \mu_{i_2,p_2-1}\right) (\mr w) \mr{B}_{i_1,p_1}(\mr r)  \frac{\mr B_{i_2,p_2-1}(\mr \varphi)}{\mr \xi_{i_2+p_2,2} - \mr \xi_{i_2,2}} 
		\label{eq: new projector derivative of projector}
	\end{align}
	is the resulting new projector.
\end{proof}

\begin{lemma}
	Let the notation of Lemmata \ref{lemma: estimate of phi-derivative} and \ref{lemma: derivative of projector} hold. Then, there is a constant $C>0$ such that
	\begin{align*}
		\norm{\mr{\Pi}^{\pol,\partial_\varphi}_{(p_1,p_2-1), \bs \Xi} (\partial_{\mr \varphi} \mr{v}_0 )  }_{\widehat L^2_{-1/2}(\mr \Omega_C)}
		\leq C \norm{\partial_{\mr \varphi}\mr{v}_0}_{\widehat L^{2}_{-1/2}\left(\widetilde{\mr{\Omega}}_C\right)}
	\end{align*}
	\label{lemma: estimate of derivative of projector}
	for all $\mr v_0 \in  C(\mr \Omega)$ with $\mr v_0(0,\cdot) = 0$.
\end{lemma}

\begin{proof}
	By definition \eqref{eq: new projector derivative of projector} and the continuous Cauchy-Schwarz inequality, it follows
	\allowdisplaybreaks
	\begin{align*}
		& \norm{\mr{\Pi}^{\pol,\partial_\varphi}_{(p_1,p_2-1), \bs \Xi} (\partial_{\mr \varphi} \mr{v}_0 )  }_{\widehat L^2_{-1/2}(\mr \Omega_C)}^2 \\
		&\leq p_2^2 \sum_{i_1=2}^{n_1} \sum_{i_2=2}^{n_2} \abs{(\lambda_{i_1,p_1} \otimes \mu_{i_2,p_2-1}) (\partial_{\mr \varphi}\mr{v}_0)}^2 \int_{\mr \Omega_C} \abs{\mr r ^{-1/2} \, \mr{B}_{i_1,p_1}(\mr r)  \frac{\mr B_{i_2,p_2-1}(\mr \varphi)}{\mr \xi_{i_2+p_2+1,2} - \mr \xi_{i_2,2}}}^2 \Dd \mr r \Dd \mr \varphi \\
		&\leq C \sum_{i_1=2}^{n_1} \sum_{i_2=2}^{n_2} \abs{(\lambda_{i_1,p_1} \otimes \mu_{i_2,p_2-1})(\partial_{\mr \varphi}\mr{v}_0)}^2 C h_{\varphi}^{-2} \\
		&=  C h_{\varphi}^{-2}  \sum_{i_1=2}^{n_1}\sum_{i_2=2}^{n_2}\abs{ \int_{\mr \xi_{i_1,1}}^{\mr \xi_{i_1+p_1+1,1}} \int_{\mr \xi_{i_2,2}}^{\mr \xi_{i_2+p_2+1,2}} \int_{\alpha(t)}^{t} \partial_{\mr \varphi}\mr{v}_0(s,z) \Dd z \, \mathring{D}^{p_2+1} \mr \psi_{i_2}(t) \, \dt \, \mathring{D}^{p_1+1} \mr \psi_{i_1}(s) \Dd s }^2 \\
		& \leq C h_{\varphi}^{-2} \sum_{i_1=2}^{n_1}\sum_{i_2=2}^{n_2}\abs{ \int_{\mr \xi_{i_1,1}}^{\mr \xi_{i_1+p_1+1,1}} \int_{\mr \xi_{i_2,2}}^{\mr \xi_{i_2+p_2+1,2}} \int_{\mr \xi_{i_2-1,2}}^{\mr \xi_{i_2+p_2+1,2}} \abs{\partial_{\mr \varphi}\mr{v}_0(s,z)} \Dd z \, \mathring{D}^{p_2+1} \mr \psi_{i_2}(t) \, \dt \, \mathring{D}^{p_1+1} \mr \psi_{i_1}(s) \Dd s }^2 \\
		&= C h_{\varphi}^{-2} \left(\int_{\mr \xi_{i_2-1,2}}^{\mr \xi_{i_2+p_2+1,2}}  \mathring{D}^{p_2+1} \mr \psi_{i_2}(t) \, \dt\right)^2 \sum_{i_1=2}^{n_1}\sum_{i_2=2}^{n_2}\abs{ \int_{\mr \xi_{i_1,1}}^{\mr \xi_{i_1+p_1+1,1}} \int_{\mr \xi_{i_2-1,2}}^{\mr \xi_{i_2+p_2+1,2}} \abs{\partial_{\mr \varphi}\mr{v}_0(s,z)} \Dd z \, \mathring{D}^{p_1+1} \mr \psi_{i_1}(s) \Dd s }^2 \\
		&\leq C \sum_{i_1=2}^{n_1}\sum_{i_2=2}^{n_2}\abs{ \int_{\mr \xi_{i_1,1}}^{\mr \xi_{i_1+p_1+1,1}} \int_{\mr \xi_{i_2-1,2}}^{\mr \xi_{i_2+p_2+1,2}} \abs{\partial_{\mr \varphi}\mr{v}_0(s,z) \mathring{D}^{p_1+1} \mr \psi_{i_1}(s)} \Dd z \Dd s }^2 . 
	\end{align*}
	Then, we set $E_{(i_1,i_2)}:=(\mr \xi_{i_1,1},\mr \xi_{i_1+p_1+1,1}) \times (\mr \xi_{i_2-1,2},\mr \xi_{i_2+p_2+1,2}) \subset \widetilde{\mr{\Omega}}$ and simplify the term above using the discrete Cauchy-Schwarz inequality and the equivalence of finite-dimensional norms,
	\begin{align*}
		& C \sum_{i_1=2}^{n_1}\sum_{i_2=2}^{n_2}\abs{ \int_{\mr \xi_{i_1,1}}^{\mr \xi_{i_1+p_1+1,1}} \int_{\mr \xi_{i_2-1,2}}^{\mr \xi_{i_2+p_2+1,2}} \abs{\partial_{\mr \varphi}\mr{v}_0(s,z) \mathring{D}^{p_1+1} \mr \psi_{i_1}(s)} \Dd z \Dd s }^2 \\
		&\leq C \sum_{i_1=2}^{n_1}\sum_{i_2=2}^{n_2} \left(\int_{E_{(i_1,i_2)}} \abs{\partial_{\mr \varphi}\mr{v}_0(s,z)}^2 \Dd z \Dd s \right) \left(\int_{E_{(i_1,i_2)}}  \abs{\mathring{D}^{p_1+1} \mr \psi_{i_1}(s)}^2 \Dd z \Dd s \right) \\
		&\leq C \left(\sum_{i_1=2}^{n_1}\sum_{i_2=2}^{n_2} \left(\int_{E_{(i_1,i_2)}} \abs{\partial_{\mr \varphi}\mr{v}_0(s,z)}^2 \Dd z \Dd s \right)^{2} \right)^{1/2} \left(\sum_{i_1=2}^{n_1}\sum_{i_2=2}^{n_2} \left( \int_{E_{(i_1,i_2)}}  \abs{\mathring{D}^{p_1+1} \mr \psi_{i_1}(s)} \Dd z \Dd s \right)^{2} \right)^{1/2} \\
		&\leq C \left(\sum_{i_1=2}^{n_1}\sum_{i_2=2}^{n_2} \int_{E_{(i_1,i_2)}} \abs{\partial_{\mr \varphi}\mr{v}_0(s,z)}^2 \Dd z \Dd s \right) \left(\sum_{i_1=2}^{n_1}\sum_{i_2=2}^{n_2} \int_{E_{(i_1,i_2)}} \abs{\mathring{D}^{p_1+1} \mr \psi_{i_1}(s)} \Dd z \Dd s \right) \\
		&\leq C \sum_{i_1=2}^{n_1}\sum_{i_2=2}^{n_2} \norm{\partial_{\mr \varphi}\mr{v}_0}^2_{L^{2}\left( E_{(i_1,i_2)}\right)} \sum_{i_1=2}^{n_1}\sum_{i_2=2}^{n_2}  h_{\varphi} h_r  
		\leq C \norm{\partial_{\mr \varphi}\mr{v}_0}_{\widehat L^{2}\left(\widetilde{\mr{\Omega}}_C\right)}^2 \cdot 1
		\leq C \norm{\partial_{\mr \varphi}\mr{v}_0}_{\widehat L^{2}_{-1/2}\left(\widetilde{\mr{\Omega}}_C\right)}^2 ,
	\end{align*}
	and the demonstration is complete.
\end{proof} 

\section{Proof of local quasi-uniformity of graded meshes}
Finally, we show an auxiliary result that is needed to prove the local quasi-uniformity of the proposed graded meshes as stated in Lemma \ref{lemma: local quasi-uniformity of the graded mesh}.
\begin{lemma}
	\label{lemma: function monotonous}
	Let $\alpha \geq 1$. The sequence $(\Theta(j))_{j = 1}^\infty$ with 
	\begin{align*}
		\Theta(j) := \frac{j^\alpha - (j-1)^\alpha}{(j+1)^\alpha - j^\alpha}
	\end{align*}
	is monotonously increasing and satisfies $\lim_{j \to \infty} \Theta(j) = 1$.
\end{lemma}
\begin{proof}
	For $\alpha=1$, it holds $\Theta(j) = 1$ for all $j \in \N$ and the assertion follows directly. For $\alpha >1$, we consider the function 
	\begin{align*}
		\Theta : [1,\infty) \to \R, \quad \Theta(x) := \frac{x^\alpha - (x-1)^\alpha}{(x+1)^\alpha - x^\alpha} ,
	\end{align*}
	which is illustrated exemplarily for $\alpha \in \{2,3,\dots,6\}$ in Figure \ref{fig: quasi-uniformity function}, where the claimed properties can be anticipated visually. 
	\begin{figure}
		\centering
		\includegraphics[width=0.7\linewidth]{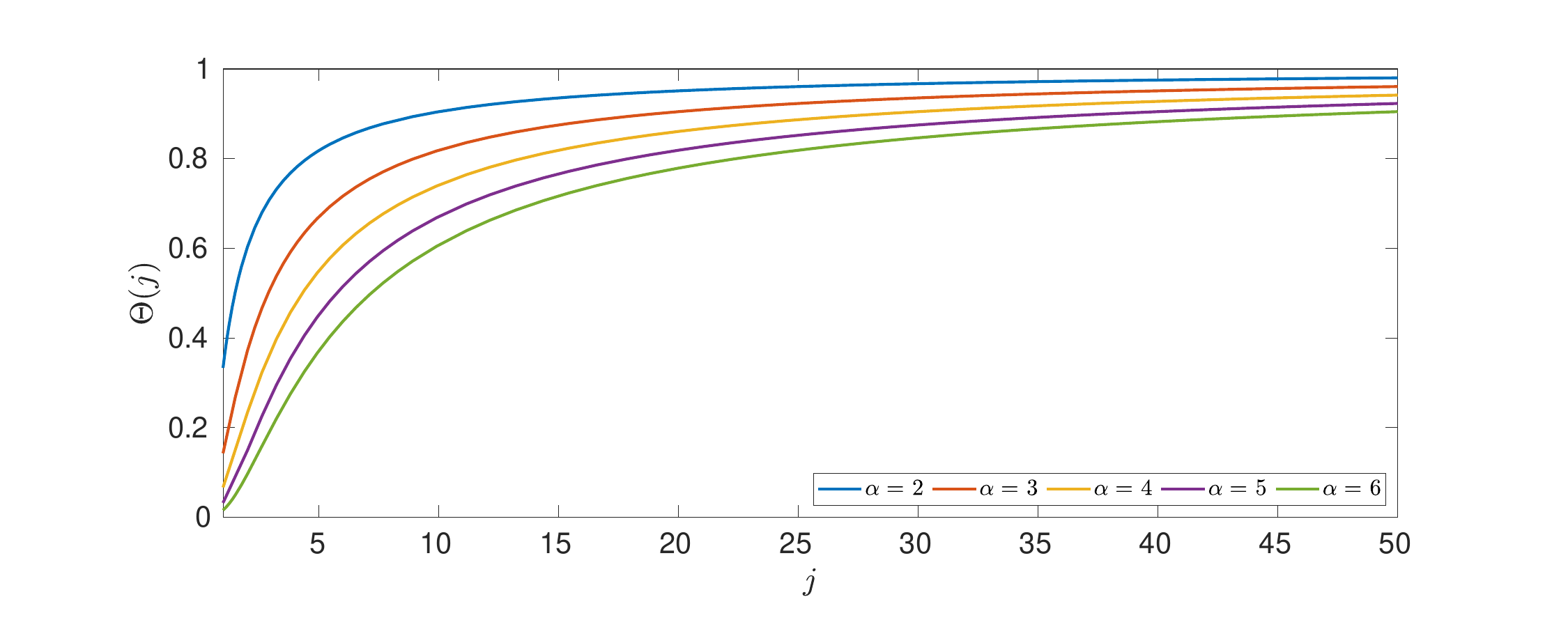}
		\caption{Illustration of the function $\Theta$ for $\alpha \in \{2,3,\dots,6\}$}
		\label{fig: quasi-uniformity function}
	\end{figure}
	In more detail, we define the functions $f:[1,\infty) \to \R, \ f(x):=x^{\alpha}$ and $g:[1,\infty) \to \R, \ g(x):=(x+1)^{\alpha}$. Then, by the generalized mean value theorem, for every $j \in \N$, there is a point $x \in (j-1,j)$ such that
	\begin{align*}
		\Theta(j) = \frac{f(j) - f(j-1)}{g(j)-g(j-1)} = \frac{f'(x)}{g'(x)}
		= \left( \frac{x}{x + 1} \right)^{\alpha - 1} =:\widetilde{\Theta}(x) .
	\end{align*}
	The derivative of the new function $\widetilde{\Theta}$ is given by
	\begin{align*}
		\widetilde{\Theta}'(x) 
		= \left(\alpha-1\right)\left( \frac{x}{x + 1} \right)^{\alpha - 2}  \frac{1}{(x + 1)^2}
		= \left(\alpha-1\right) x^{\alpha - 2} (x + 1)^{-\alpha},
	\end{align*}
	and, since $\alpha > 1$ and $x > j-1 \geq 0$, it always holds $\widetilde{\Theta}'(x) >0$. Thus, the function $\widetilde{\Theta}$ and consequently also the sequence $(\Theta(j))_{j = 1}^\infty$ are monotonously increasing. Moreover, we have 
	\begin{align*}
		\lim\limits_{j \to \infty} \Theta(j) 
		= \lim\limits_{x \to \infty} \widetilde{\Theta}(x) 
		= \lim\limits_{x \to \infty} \left( \frac{x}{x + 1} \right)^{\alpha - 1} = 1 
	\end{align*}
	and the proof is concluded.
\end{proof}

\section*{Declaration of generative AI and AI-assisted technologies in the writing process}
During the preparation of this work, the authors used ChatGPT 4o and DeepL Write to fix writing errors and to improve the readability and language in some parts of the manuscript. After using these tools, the authors reviewed and edited the content and take full responsibility for the content of the published article.



\biboptions{sort&compress}
\bibliographystyle{elsarticle-harv} 
\bibliography{Bibliography.bib}


%
%
%
%
%
\end{document}

%% file: param_circ.tikz
	\begin{tikzpicture}[baseline, decoration = {snake, pre length=3pt, post length=3pt,},line width=0.4mm, scale = 0.6]
	\fill[gray!30] (0,0) rectangle (4,4);		
	\draw[orange] (0,4) -- (4, 4);
	\draw[cyan] (0,0) -- (0, 4);
	\draw[olive] (0,0) -- (4, 0);
	\draw[blue] (4,0) -- (4, 4);
	\draw[dash pattern=on 2pt off 1pt, line width=0.05mm] (0,1.333) -- (4,1.333);
	\draw[dash pattern=on 2pt off 1pt, line width=0.05mm] (0,2.666) -- (4,2.666);
	
	\node[gray] at (2,2) {$\widehat{\Omega}$};
	\node[left,yshift=2] at(0,0) {\small $0$};
	\node[left] at(0, 4) {\small $1$};
	
	\node[left,cyan] at (0, 2) {\small $\widehat{\Gamma}_1$};
	\node[right,blue] at (4, 2) {\small $\widehat{\Gamma}_2$};
	\node[below,olive] at (2, 0) {\small $\widehat{\Gamma}_3$};
	\node[above,orange] at (2, 4) {\small $\widehat{\Gamma}_4$};
	
	\node at (0,-0.4) {\small $0$};
	\node at (4,-0.4) {\small $1$};
	
	\draw[-stealth] (4, 0) -- (4.8,0) node[right] {$\zeta_1$};
	\draw[-stealth] (0, 4) -- (0, 4.8) node[above] {$\zeta_2$};
	
	\draw [-stealth, very thick] (5.1,2) to [out=30,in=150] (6.9,2) node[above] at (6,2.25) {$\bs F$};
\end{tikzpicture}

%% file: phys_circ.tikz
\begin{tikzpicture}[scale=1.4, baseline=-12mm, line width=0.4mm]	
	\filldraw[fill = gray, fill opacity = 0.3, draw = white] (1,0) arc (0:300:1) -- (0,0) -- cycle;
	\draw[blue] (1,0) arc (0:300:1) node[pos = 0.5, rotate = 50, below] {$\Gamma_2$};	
	\draw[olive] (0,0.013) -- (1, 0.013) node[midway, below, xshift =5] {$\Gamma_3$};
	\draw[orange] (0,0) -- (300:1) node[midway, right, yshift=-3] {$\Gamma_4$};				
	\draw[-stealth] (0.1:0.3) arc (0:300:0.3) node[sloped, midway, below] {$\omega$};
	
	\draw[dash pattern=on 2pt off 1pt, line width=0.05mm] (100:1) -- (0, 0);
	\draw[dash pattern=on 2pt off 1pt, line width=0.05mm] (200:1) -- (0, 0);
	
	\fill[cyan] (0,0) circle (1pt);
	
	\node[cyan] at (330:0.3) {$\bs P$};		
	\node[gray] at (225:0.625) {$\Omega$};
\end{tikzpicture}

%% file: param_L.tikz
\begin{tikzpicture}[baseline, decoration = {snake, pre length=3pt, post length=3pt,},line width=0.4mm, scale = 0.6]
	\fill[gray!30] (0,0) rectangle (4,4);		
	\draw[orange] (0,4) -- (4, 4);
	\draw[cyan] (0,0) -- (0, 4);
	\draw[olive] (0,0) -- (4, 0);
	\draw[blue] (4,0) -- (4, 4);
	\draw[dash pattern=on 2pt off 1pt, line width=0.05mm] (0,1) -- (4,1);
	\draw[dash pattern=on 2pt off 1pt, line width=0.05mm] (0,2) -- (4,2);
	\draw[dash pattern=on 2pt off 1pt, line width=0.05mm] (0,3) -- (4,3);
	
	\node[gray] at (2,2.5) {$\widehat{\Omega}$};
	\node[left,yshift=2] at(0,0) {\small $0$};
	\node[left] at(0, 4) {\small $1$};
	
	\node[left,cyan] at (0, 2) {\small $\widehat{\Gamma}_1$};
	\node[right,blue] at (4, 2) {\small $\widehat{\Gamma}_2$};
	\node[below,olive] at (2, 0) {\small $\widehat{\Gamma}_3$};
	\node[above,orange] at (2, 4) {\small $\widehat{\Gamma}_4$};
	
	\node at (0,-0.4) {\small $0$};
	\node at (4,-0.4) {\small $1$};
	
	\draw[-stealth] (4, 0) -- (4.8,0) node[right] {$\zeta_1$};
	\draw[-stealth] (0, 4) -- (0, 4.8) node[above] {$\zeta_2$};
	
	\draw [-stealth, very thick] (5.1,2) to [out=30,in=150] (6.9,2) node[above] at (6,2.25) {$\bs F$};
\end{tikzpicture}

%% file: phys_L.tikz
	\begin{tikzpicture}[scale=1.4, baseline=-12mm, line width=0.4mm]	
	\filldraw[fill = gray, fill opacity = 0.3, draw = white] (-1,-1) -- (-1,1) -- (1,1) -- (1,0) -- (0,0) -- (0,-1) -- (-1,-1);
	\draw[blue] (-1,1) -- (1,1) node[midway, below] {$\Gamma_2$};
	\draw[blue] (-1,1) -- (-1,-1); 
	\draw[blue] (-1,-1) -- (0,-1);
	\draw[blue] (1,1) -- (1,0);
	
	\draw[dash pattern=on 2pt off 1pt, line width=0.05mm] (1,1) -- (0, 0);
	\draw[dash pattern=on 2pt off 1pt, line width=0.05mm] (-1,1) -- (0, 0);
	\draw[dash pattern=on 2pt off 1pt, line width=0.05mm] (-1,-1) -- (0, 0);
	
	\draw[olive] (0,0) -- (1, 0) node[midway, below,xshift=5] {$\Gamma_3$};
	\draw[orange] (0,0) -- (0, -1) node[midway, right, yshift=-5] {$\Gamma_4$};
	
	\draw[-stealth] (0.1:0.3) arc (0:270:0.3);
	
	\node[cyan] at (315:0.25) {$\bs P$};
	\node at (90:0.15) {$\omega$};	
	\node[gray] at (180:0.625) {$\Omega$};
	
	\fill[cyan] (0,0) circle (1pt);
\end{tikzpicture}

%% file: grad_mesh_circ_par.tikz
	\begin{tikzpicture}[baseline, scale = 0.675]
	
	\def \N {5}			
	\def \hxi{1/6}		
	
	\def \M {3} 		
	\def \heta{1/4}		
	
	\node at (2,4.5) {$\widehat{\mathcal{M}}^\mu$};
	\node[gray] at (2,-0.5) {$\widehat{\Omega}$};
	
	\fill[gray!30] (0,0) rectangle (4,4);
	
	\draw[-stealth] (4,0) -- (5, 0)  node[below, xshift = -2] {$\zeta_1$};
	\draw[-stealth] (0,4) -- (0, 5) node[left, yshift = -2] {$\zeta_2$};

	\foreach \i in {0,...,3}{
		\draw (\i*\i*\hxi*\hxi*4,0) -- (\i*\i*\hxi*\hxi*4,4);
	}
	
	\foreach \i in {4,5,6}{
		\draw (\i*\i*\hxi*\hxi*4,0) -- (\i*\i*\hxi*\hxi*4,4);
	}
	
	\foreach \i in {1,...,2}{
		\draw[dashed] (0, \i*4/3) -- (1, \i*4/3);
		\draw[dashed] (1, \i*4/3) -- (4, \i*4/3);
	}
	
	\foreach \j in {0,...,2}{
		\foreach \i in {1,...,\M}{
			\draw (0, \i*\heta*4/3 + \j*4/3) -- (1, \i*\heta*4/3 + \j*4/3);
			\draw (1, \i*\heta*4/3 + \j*4/3) -- (4, \i*\heta*4/3 + \j*4/3);
		}
	}
	
	\draw[orange, line width=0.4mm] (0,4) -- (4, 4);
	\draw[cyan, line width=0.4mm] (0,0) -- (0, 4);
	\draw[olive, line width=0.4mm] (0,0) -- (4, 0);
	\draw[blue, line width=0.4mm] (4,0) -- (4, 4);

	\node[left, yshift = 2] at (0,0) {$0$};
	\node[left] at (0,4) {$1$};
	
	\node[below] at (0,0) {$0$};;
	\node[below] at (4,0) {$1$};
	
	\draw [-stealth, very thick] (4.6,2) to [out=30,in=150] (6.1,2) node[above] at (5.35,2.3) {$\bs F$};
\end{tikzpicture}

%% file: grad_mesh_circ_phys.tikz
				\begin{tikzpicture}[baseline=-1.3cm, scale = 1.4]	
	
	\def \N {5}			
	\def \hxi{1/6}		
	\def \M {3} 		
	\def \heta{1/4}		
	
	\filldraw[fill = gray!30, draw = white] (1,0) arc (0:300:1) -- (0,0);
	
	\draw[blue,line width=0.4mm] (1,0) arc (0:300:1);	
	\draw[olive,line width=0.4mm] (0,0.013) -- (1, 0.013);
	\draw[orange,line width=0.4mm] (0,0) -- (0.5,-0.86603);
	\fill[cyan] (0,0) circle (1pt);		
	
	\foreach \i in {0,...,6}{	
		\draw (\i*\i*\hxi*\hxi,0) arc (0:300:\i*\i*\hxi*\hxi);
	}

	\foreach \j in {0,...,2}{
		\foreach \i in {1,...,\M}{
			\draw (0,0) -- (\i*\heta*100 + \j*100:1);
		}
	}		
	
	\draw[dash pattern=on 2pt off 1pt] (0,0) -- (100:1);					
	\draw[dash pattern=on 2pt off 1pt] (0,0) -- (200:1);

	\node at (0.5,-0.2) {$\mathcal{M}^\mu $};			
	\node[gray] at (0,1.2) {$\Omega$};	
\end{tikzpicture}

%% file: grad_mesh_L_par.tikz
\begin{tikzpicture}[baseline, scale = 0.675]
	
	\def \N {5}			
	\def \hxi{1/6}		
	
	\def \M {3} 		
	\def \heta{1/4}		
	
	\node at (2,4.5) {$\widehat{\mathcal{M}}^\mu$};
	\node[gray] at (2,-0.5) {$\widehat{\Omega}$};

	\fill[gray!30] (0,0) rectangle (1,4);
	\fill[gray!30] (1,0) rectangle (4,4);
	
	\draw[-stealth] (4,0) -- (5, 0)  node[below, xshift = -2] {$\zeta_1$};
	\draw[-stealth] (0,4) -- (0, 5) node[left, yshift = -2] {$\zeta_2$};

	\foreach \i in {0,...,3}{
		\draw (\i*\i*\hxi*\hxi*4,0) -- (\i*\i*\hxi*\hxi*4,4);
	}
	
	\foreach \i in {4,5,6}{
		\draw (\i*\i*\hxi*\hxi*4,0) -- (\i*\i*\hxi*\hxi*4,4);
	}
	
	\foreach \i in {1,...,3}{
		\draw[dashed] (0, \i) -- (1, \i);
		\draw[dashed] (1, \i) -- (4, \i);
	}
	
	\foreach \j in {0,...,3}{
		\foreach \i in {1,...,\M}{
			\draw (0, \i*\heta + \j) -- (1, \i*\heta + \j);
			\draw (1, \i*\heta + \j) -- (4, \i*\heta + \j);
		}
	}
	
	\draw[orange, line width=0.4mm] (0,4) -- (4, 4);
	\draw[cyan, line width=0.4mm] (0,0) -- (0, 4);
	\draw[olive, line width=0.4mm] (0,0) -- (4, 0);
	\draw[blue, line width=0.4mm] (4,0) -- (4, 4);
	
	\node[left, yshift = 2] at (0,0) {$0$};
	\node[left] at (0,4) {$1$};
	
	\node[below] at (0,0) {$0$};
	\node[below] at (4,0) {$1$};
	
	\draw [-stealth, very thick] (4.6,2) to [out=30,in=150] (6.1,2) node[above] at (5.35,2.3) {$\bs F$};
\end{tikzpicture}

%% file: grad_mesh_L_phys.tikz
\begin{tikzpicture}[scale=1.4, baseline=-13mm]
	
	\def \N {5}			
	\def \hxi{1/6}		
	
	\def \M {3} 		
	\def \heta{1/4}		
	
	\filldraw[fill = gray!30, draw = white] (-1,-1) -- (-1,1) -- (1,1) -- (1,0) -- (0.25,0) -- (0.25,0.25) -- (-0.25,0.25) -- (-0.25,-0.25) -- (0,-0.25) -- (0,-1);
	\filldraw[fill = gray!30, draw = white] (0.25,0) -- (0.25,0.25) -- (-0.25,0.25) -- (-0.25,-0.25) -- (0,-0.25) -- (0,0);
	\draw[blue,line width =0.4mm] (-1,1) -- (1,1);
	\draw[blue,line width =0.4mm] (-1,1) -- (-1,-1); 
	\draw[blue,line width =0.4mm] (-1,-1) -- (0,-1);
	\draw[blue,line width =0.4mm] (1,1) -- (1,0);
	\draw[orange,line width =0.4mm] (0,-1) -- (0,0);
	\draw[olive,line width =0.4mm] (1,0) -- (0,0);
	
	\draw[dash pattern=on 2pt off 1pt] (1,1) -- (0.25, 0.25);
	\draw[dash pattern=on 2pt off 1pt] (-1,1) -- (-0.25, 0.25);
	\draw[dash pattern=on 2pt off 1pt] (-1,-1) -- (-0.25, -0.25);
	
	\draw[dash pattern=on 2pt off 1pt, blue] (0,0) -- (0.25, 0.25);
	\draw[dash pattern=on 2pt off 1pt, blue] (0,0) -- (-0.25, 0.25);
	\draw[dash pattern=on 2pt off 1pt, blue] (0,0) -- (-0.25, -0.25);

	\foreach \i in {0,...,3}{
		\draw (\i*\i*\hxi*\hxi,0) -- (\i*\i*\hxi*\hxi,\i*\i*\hxi*\hxi);
		\draw (\i*\i*\hxi*\hxi,\i*\i*\hxi*\hxi) -- (-\i*\i*\hxi*\hxi,\i*\i*\hxi*\hxi);
		\draw (-\i*\i*\hxi*\hxi,\i*\i*\hxi*\hxi) -- (-\i*\i*\hxi*\hxi,-\i*\i*\hxi*\hxi);
		\draw (-\i*\i*\hxi*\hxi,-\i*\i*\hxi*\hxi) -- (0,-\i*\i*\hxi*\hxi);
	}
	
	\foreach \i in {4,...,6}{
		\draw (\i*\i*\hxi*\hxi,0) -- (\i*\i*\hxi*\hxi,\i*\i*\hxi*\hxi);
		\draw (\i*\i*\hxi*\hxi,\i*\i*\hxi*\hxi) -- (-\i*\i*\hxi*\hxi,\i*\i*\hxi*\hxi);
		\draw (-\i*\i*\hxi*\hxi,\i*\i*\hxi*\hxi) -- (-\i*\i*\hxi*\hxi,-\i*\i*\hxi*\hxi);
		\draw (-\i*\i*\hxi*\hxi,-\i*\i*\hxi*\hxi) -- (0,-\i*\i*\hxi*\hxi);
	}
	
	\draw (0,0) -- (0.25,0.0497);	
	\draw (0,0) -- (0.25,0.1035);	
	\draw (0,0) -- (0.25,0.1670);
	
	\draw (0.25,0.0497) -- (1,0.1989);	
	\draw (0.25,0.1035) -- (1,0.4142);	
	\draw (0.25,0.1670) -- (1,0.6681);
	
	\draw (0,0) -- (0.1035,0.25);	
	\draw (0,0) -- (0,0.25);	
	\draw (0,0) -- (-0.1035,0.25);
	
	\draw (0.1035,0.25) -- (0.41421,1);	
	\draw (0,0.25) -- (0,1);	
	\draw (-0.1035,0.25) -- (-0.41421,1);
	
	\draw (0,0) -- (-0.25,0.1035);	
	\draw (0,0) -- (-0.25,0);	
	\draw (0,0) -- (-0.25,-0.1035);
	
	\draw (-1,0.41421) -- (-0.25,0.1035);	
	\draw (-1,0) -- (-0.25,0);	
	\draw (-1,-0.41421) -- (-0.25,-0.1035);
	
	\draw (0,0) -- (-0.0497,-0.25);	
	\draw (0,0) -- (-0.1035,-0.25);	
	\draw (0,0) -- (-0.1670,-0.25);
	
	\draw (-0.1989,-1) -- (-0.0497,-0.25);	
	\draw (-0.4142,-1) -- (-0.1035,-0.25);	
	\draw (-0.6681,-1) -- (-0.1670,-0.25);
	
	\node at (0.4,-0.3) {$\mathcal{M}^\mu $};			
	\node[gray] at (0,1.2) {$\Omega$};	
	
\end{tikzpicture}

%% file: splitting_circ_par.tikz
\begin{tikzpicture}[baseline, scale = 0.675]
	
	\def \N {5}			
	\def \hxi{1/6}		
	
	\def \M {3} 		
	\def \heta{1/4}		
	
	\node[blue] at (0.75,4.7) {$\widehat{\mathcal{M}}^\mu_C$};
	\node at (2.75,4.7) {$\widehat{\mathcal{M}}^\mu_R$};
	\node[gray!80] at (3/8*3/8*4,-1.2) {$\widehat{\Omega}_C$};
	\node[orange!80] at (2.75,-1.2) {$\widehat{\Omega}_R$};
	
	\fill[gray!30] (0,0) rectangle (1,4);
	\fill[orange!30] (1,0) rectangle (4,4);
	
	\draw[-stealth] (4,0) -- (5, 0)  node[below, xshift = -2] {$\zeta_1$};
	\draw[-stealth] (0,4) -- (0, 5) node[left, yshift = -2] {$\zeta_2$};

	\draw[blue] (0,0) -- (1, 0);
	\draw[blue] (0,4) -- (1, 4);
	\draw (1,0) -- (4, 0);
	\draw (1,4) -- (4, 4);
	
	\foreach \i in {0,...,3}{
		\draw[blue] (\i*\i*\hxi*\hxi*4,0) -- (\i*\i*\hxi*\hxi*4,4);
	}
	
	\foreach \i in {4,5,6}{
		\draw (\i*\i*\hxi*\hxi*4,0) -- (\i*\i*\hxi*\hxi*4,4);
	}
	
	\foreach \i in {1,...,2}{
		\draw[dashed,blue] (0, \i*4/3) -- (1, \i*4/3);
		\draw[dashed] (1, \i*4/3) -- (4, \i*4/3);
	}
	
	\foreach \j in {0,...,2}{
		\foreach \i in {1,...,\M}{
			\draw[blue] (0, \i*\heta*4/3 + \j*4/3) -- (1, \i*\heta*4/3 + \j*4/3);
			\draw (1, \i*\heta*4/3 + \j*4/3) -- (4, \i*\heta*4/3 + \j*4/3);
		}
	}
	
	\node[left, yshift = 2] at (0,0) {$0$};
	\node[left] at (0,4) {$1$};
	
	\node[below, xshift = -2] at (0,0) {$0$};
	\node[below, yshift = 1] at (\hxi*\hxi*3*3*4,0) {$\zeta_{1,p_1+2}$};
	\node[below] at (4,0) {$1$};
	
	\draw [-stealth, very thick] (4.6,2) to [out=30,in=150] (6.1,2) node[above] at (5.35,2.3) {$\bs F$};
\end{tikzpicture}

%% file: splitting_circ_phys.tikz
	\begin{tikzpicture}[baseline=-1.3cm, scale = 1.4]	
	
	\def \N {5}			
	\def \hxi{1/6}		
	\def \M {3} 		
	\def \heta{1/4}		
	
	\filldraw[fill = orange!30, draw = white] (0:1) arc (0:300:1) -- (300:0.25) arc (300:0:0.25);
	\filldraw[fill = gray!30, draw = white] (0:0.25) arc (0:300:0.25) -- (0:0);
	
	\foreach \i in {0,...,3}{	
		\draw[blue] (\i*\i*\hxi*\hxi,0) arc (0:300:\i*\i*\hxi*\hxi);
	}
	
	\foreach \i in {4,...,6}{	
		\draw (\i*\i*\hxi*\hxi,0) arc (0:300:\i*\i*\hxi*\hxi);
	}
	
	\foreach \j in {0,...,2}{
		\foreach \i in {1,...,\M}{
			\draw (\i*\heta*100 + \j*100:0.25) -- (\i*\heta*100 + \j*100:1);
			\draw[blue] (\i*\heta*100 + \j*100:0.25) -- (\i*\heta*100 + \j*100:0);
		}
	}		
	
	\draw (0:1) -- (0:0.25);		
	\draw (300:1) -- (300:0.25);
	
	\draw[blue] (0:0) -- (0:0.25);		
	\draw[blue] (0:0) -- (300:0.25);
	
	\draw[dash pattern=on 2pt off 1pt] (100:1) -- (100:0.25);		
	\draw[dash pattern=on 2pt off 1pt] (200:1) -- (200:0.25);
	
	\draw[dash pattern=on 2pt off 1pt,blue] (0,0) -- (100:0.25);	
	\draw[dash pattern=on 2pt off 1pt,blue] (0,0) -- (200:0.25);	
	
	\node[blue] at (0.75,-0.4) {$\mathcal{M}^\mu_C $};
	\node at (0.333,1.2) {$\mathcal{M}^\mu_R$};
	
	\node[gray!80] at (0.37,-0.2) {$\Omega_C$};
	\node[orange!80] at (-0.333,1.2) {$\Omega_R$};	
\end{tikzpicture}

%% file: splitting_L_par.tikz
	\begin{tikzpicture}[baseline, scale = 0.675]
	
	\def \N {5}			
	\def \hxi{1/6}		
	
	\def \M {3} 		
	\def \heta{1/4}		
	
	\node[blue] at (0.75,4.7) {$\widehat{\mathcal{M}}^\mu_C$};
	\node at (2.75,4.7) {$\widehat{\mathcal{M}}^\mu_R$};
	\node[gray!60] at (3/8*3/8*4,-1.2) {$\widehat{\Omega}_C$};
	\node[orange!60] at (2.75,-1.2) {$\widehat{\Omega}_R$};
	
	\fill[gray!30] (0,0) rectangle (1,4);
	\fill[orange!30] (1,0) rectangle (4,4);
	
	\draw[-stealth] (4,0) -- (5, 0)  node[below, xshift = -2] {$\zeta_1$};
	\draw[-stealth] (0,4) -- (0, 5) node[left, yshift = -2] {$\zeta_2$};

	\draw[blue] (0,0) -- (1, 0);
	\draw[blue] (0,4) -- (1, 4);
	\draw (1,0) -- (4, 0);
	\draw (1,4) -- (4, 4);
	
	\foreach \i in {0,...,3}{
		\draw[blue] (\i*\i*\hxi*\hxi*4,0) -- (\i*\i*\hxi*\hxi*4,4);
	}
	
	\foreach \i in {4,5,6}{
		\draw (\i*\i*\hxi*\hxi*4,0) -- (\i*\i*\hxi*\hxi*4,4);
	}
	
	\foreach \i in {1,...,3}{
		\draw[dashed,blue] (0, \i) -- (1, \i);
		\draw[dashed] (1, \i) -- (4, \i);
	}
	
	\foreach \j in {0,...,3}{
		\foreach \i in {1,...,\M}{
			\draw[blue] (0, \i*\heta + \j) -- (1, \i*\heta + \j);
			\draw (1, \i*\heta + \j) -- (4, \i*\heta + \j);
		}
	}
	
	\node[left, yshift = 2] at (0,0) {$0$};
	\node[left] at (0,4) {$1$};
	
	\node[below, xshift = -2] at (0,0) {$0$};
	\node[below, yshift = 1] at (\hxi*\hxi*3*3*4,0) {$\zeta_{1,p_1+2}$};
	\node[below] at (4,0) {$1$};
	
	\draw [-stealth, very thick] (4.6,2) to [out=30,in=150] (6.1,2) node[above] at (5.35,2.3) {$\bs F$};
\end{tikzpicture}

%% file: splitting_L_phys.tikz
\begin{tikzpicture}[scale=1.4, baseline=-13mm]
	
	\def \N {5}			
	\def \hxi{1/6}		
	
	\def \M {3} 		
	\def \heta{1/4}		
	
	\filldraw[fill = orange!30, draw = white] (-1,-1) -- (-1,1) -- (1,1) -- (1,0) -- (0.25,0) -- (0.25,0.25) -- (-0.25,0.25) -- (-0.25,-0.25) -- (0,-0.25) -- (0,-1);
	\filldraw[fill = gray!30, draw = white] (0.25,0) -- (0.25,0.25) -- (-0.25,0.25) -- (-0.25,-0.25) -- (0,-0.25) -- (0,0);
	\draw (-1,1) -- (1,1);
	\draw (-1,1) -- (-1,-1); 
	\draw (-1,-1) -- (0,-1);
	\draw (1,1) -- (1,0);
	\draw (0,-1) -- (0,-0.25);
	\draw[blue] (0,0) -- (0,-0.25);
	\draw (1,0) -- (0.25,0);
	\draw[blue] (0,0) -- (0.25,0);
	
	\draw[dash pattern=on 2pt off 1pt] (1,1) -- (0.25, 0.25);
	\draw[dash pattern=on 2pt off 1pt] (-1,1) -- (-0.25, 0.25);
	\draw[dash pattern=on 2pt off 1pt] (-1,-1) -- (-0.25, -0.25);
	
	\draw[dash pattern=on 2pt off 1pt, blue] (0,0) -- (0.25, 0.25);
	\draw[dash pattern=on 2pt off 1pt, blue] (0,0) -- (-0.25, 0.25);
	\draw[dash pattern=on 2pt off 1pt, blue] (0,0) -- (-0.25, -0.25);

	\foreach \i in {0,...,3}{
		\draw[blue] (\i*\i*\hxi*\hxi,0) -- (\i*\i*\hxi*\hxi,\i*\i*\hxi*\hxi);
		\draw[blue] (\i*\i*\hxi*\hxi,\i*\i*\hxi*\hxi) -- (-\i*\i*\hxi*\hxi,\i*\i*\hxi*\hxi);
		\draw[blue] (-\i*\i*\hxi*\hxi,\i*\i*\hxi*\hxi) -- (-\i*\i*\hxi*\hxi,-\i*\i*\hxi*\hxi);
		\draw[blue] (-\i*\i*\hxi*\hxi,-\i*\i*\hxi*\hxi) -- (0,-\i*\i*\hxi*\hxi);
	}
	
	\foreach \i in {4,...,6}{
		\draw (\i*\i*\hxi*\hxi,0) -- (\i*\i*\hxi*\hxi,\i*\i*\hxi*\hxi);
		\draw (\i*\i*\hxi*\hxi,\i*\i*\hxi*\hxi) -- (-\i*\i*\hxi*\hxi,\i*\i*\hxi*\hxi);
		\draw (-\i*\i*\hxi*\hxi,\i*\i*\hxi*\hxi) -- (-\i*\i*\hxi*\hxi,-\i*\i*\hxi*\hxi);
		\draw (-\i*\i*\hxi*\hxi,-\i*\i*\hxi*\hxi) -- (0,-\i*\i*\hxi*\hxi);
	}
	
	\draw[blue] (0,0) -- (0.25,0.0497);	
	\draw[blue] (0,0) -- (0.25,0.1035);	
	\draw[blue] (0,0) -- (0.25,0.1670);
	
	\draw (0.25,0.0497) -- (1,0.1989);	
	\draw (0.25,0.1035) -- (1,0.4142);	
	\draw (0.25,0.1670) -- (1,0.6681);
	
	\draw[blue] (0,0) -- (0.1035,0.25);	
	\draw[blue] (0,0) -- (0,0.25);	
	\draw[blue] (0,0) -- (-0.1035,0.25);
	
	\draw (0.1035,0.25) -- (0.41421,1);	
	\draw (0,0.25) -- (0,1);	
	\draw (-0.1035,0.25) -- (-0.41421,1);
	
	\draw[blue] (0,0) -- (-0.25,0.1035);	
	\draw[blue] (0,0) -- (-0.25,0);	
	\draw[blue] (0,0) -- (-0.25,-0.1035);
	
	\draw (-1,0.41421) -- (-0.25,0.1035);	
	\draw (-1,0) -- (-0.25,0);	
	\draw (-1,-0.41421) -- (-0.25,-0.1035);
	
	\draw[blue] (0,0) -- (-0.0497,-0.25);	
	\draw[blue] (0,0) -- (-0.1035,-0.25);	
	\draw[blue] (0,0) -- (-0.1670,-0.25);
	
	\draw (-0.1989,-1) -- (-0.0497,-0.25);	
	\draw (-0.4142,-1) -- (-0.1035,-0.25);	
	\draw (-0.6681,-1) -- (-0.1670,-0.25);
	
	\node[blue] at (0.5,-0.5) {$\mathcal{M}^\mu_C $};
	\node at (0.333,1.2) {$\mathcal{M}^\mu_R$};
	
	\node[gray!80] at (0.25,-0.2) {$\Omega_C$};
	\node[orange!80] at (-0.333,1.2) {$\Omega_R$};	
	
\end{tikzpicture}

%% file: ref_elem.tikz
	\begin{tikzpicture}[baseline, scale = 0.6]
	
	\def \N {5}			
	\def \hxi{1/6}		
	
	\def \M {3} 		
	\def \heta{1/4}		
	
	\node[gray!80] at (2,4.5) {$\mr{\Omega}_C$};
	\node at (8,4.5) {$\mr{\Omega}$};
	
	\fill[gray!30] (0,0) rectangle (4,4);	
	\fill[orange!10] (4,0) rectangle (16,4);
	
	\draw[-stealth] (16,0) -- (17, 0)  node[below, xshift = -2] {$\mr r$};
	\draw[-stealth] (0,4) -- (0, 5) node[left, yshift = -2] {$\mr \varphi$};

	\draw[blue] (0,0) -- (4, 0);
	\draw[blue] (0,4) -- (4, 4);
	\draw[black!20]  (4,0) -- (4*4, 0);
	\draw[black!20]  (4,4) -- (4*4, 4);
	
	\foreach \i in {0,...,3}{
		\draw[blue] (\i*\i*\hxi*\hxi*4*4,0) -- (\i*\i*\hxi*\hxi*4*4,4);
	}
	
	\foreach \i in {4,5,6}{
		\draw[black!20] (\i*\i*\hxi*\hxi*4*4,0) -- (\i*\i*\hxi*\hxi*4*4,4);
	}
	
	\foreach \i in {1,...,2}{
		\draw[dashed,blue] (0, \i*4/3) -- (4, \i*4/3);
		\draw[dashed,black!20] (4, \i*4/3) -- (4*4, \i*4/3);
	}
	
	\foreach \j in {0,...,2}{
		\foreach \i in {1,...,\M}{
			\draw[blue] (0, \i*\heta*4/3 + \j*4/3) -- (4,\i*\heta*4/3 + \j*4/3);
			\draw[black!20] (4, \i*\heta*4/3 + \j*4/3) -- (4*4, \i*\heta*4/3 + \j*4/3);
		}
	}
	
	\node[below,orange!40] at (15.5,0) {$1/h_{r, \widehat{\Omega}_C}$};
	
	\node[left, yshift = 2] at (0,0) {$0$};
	\node[left] at (0,4) {$1$};
	
	\node[below, xshift = -2] at (0,0) {$0$};
	\node[below] at (4,0) {$1$};
	
	\draw [-stealth, very thick] (17,2) to [out=30,in=150] (18.5,2) node[above] at (17.75,2.3) {$\mr {\bs \Phi}$};
\end{tikzpicture}

%% file: ref_elem_ori.tikz
\begin{tikzpicture}[baseline, scale = 0.6]
	
	\def \N {5}			
	\def \hxi{1/6}		
	
	\def \M {3} 		
	\def \heta{1/4}		

	\node[gray!80] at (3/8*3/8*4,4.5) {$\widehat{\Omega}_C$};
	\node at (2,4.5) {$\widehat{\Omega}$};
	
	
	\fill[gray!30] (0,0) rectangle (1,4);
	\fill[orange!10] (1,0) rectangle (4,4);
	
	\draw[-stealth] (4,0) -- (5, 0)  node[below, xshift = -2] {$r$};
	\draw[-stealth] (0,4) -- (0, 5) node[left, yshift = -2] {$\varphi$};

	\draw[blue] (0,0) -- (1, 0);
	\draw[blue] (0,4) -- (1, 4);
	\draw[black!20]  (1,0) -- (4, 0);
	\draw[black!20]  (1,4) -- (4, 4);
	
	\foreach \i in {0,...,3}{
		\draw[blue] (\i*\i*\hxi*\hxi*4,0) -- (\i*\i*\hxi*\hxi*4,4);
	}
	
	\foreach \i in {4,5,6}{
		\draw[black!20] (\i*\i*\hxi*\hxi*4,0) -- (\i*\i*\hxi*\hxi*4,4);
	}
	
	\foreach \i in {1,...,2}{
		\draw[dashed,blue] (0, \i*4/3) -- (1, \i*4/3);
		\draw[dashed,black!20] (1, \i*4/3) -- (4, \i*4/3);
	}
	
	\foreach \j in {0,...,2}{
		\foreach \i in {1,...,\M}{
			\draw[blue] (0, \i*\heta*4/3 + \j*4/3) -- (1, \i*\heta*4/3 + \j*4/3);
			\draw[black!20]  (1, \i*\heta*4/3 + \j*4/3) -- (4, \i*\heta*4/3 + \j*4/3);
		}
	}
	
	\node[left, yshift = 2] at (0,0) {$0$};
	\node[left] at (0,4) {$1$};
	
	\node[below, xshift = -2] at (0,0) {$0$};
	\node[below, yshift = 1,blue] at (\hxi*\hxi*3*3*4,0) {$h_{r,\widehat{\Omega}_C}$};
	\node[below] at (4,0) {$1$};	
\end{tikzpicture}

%% file: ref_elem_ext.tikz
	\begin{tikzpicture}[baseline, scale = 0.6]		
	\def \N {5}			
	\def \hxi{1/6}		
	
	\def \M {3} 		
	\def \heta{1/4}		
	
	\node[green!80] at (4,4.5) {$\widetilde{\mr{\Omega}}_C$};
	\node[gray!80] at (2,4.5) {$\mr{\Omega}_C$};
	\node at (8,4.5) {$\mr{\Omega}$};

	\fill[gray!30] (0,0) rectangle (4,4);
	\draw[pattern = north west lines, pattern color = green] (0,0) rectangle (7.11,4);
	\fill[orange!10] (7.11,0) rectangle (16,4);
	
	\draw[-stealth] (16,0) -- (17, 0)  node[above, xshift = -2] {$\mr r$};
	\draw[-stealth] (0,4) -- (0, 5) node[left, yshift = -2] {$\mr \varphi$};

	\draw[blue] (0,0) -- (4, 0);
	\draw[blue] (0,4) -- (4, 4);
	\draw[red] (4,0) -- (7.11, 0);
	\draw[red] (4,4) -- (7.11, 4);
	\draw[red] (7.11,0) -- (7.11, 4);
	\draw[black!20]  (7.11,0) -- (4*4, 0);
	\draw[black!20]  (7.11,4) -- (4*4, 4);
	
	\foreach \i in {0,...,3}{
		\draw[blue] (\i*\i*\hxi*\hxi*4*4,0) -- (\i*\i*\hxi*\hxi*4*4,4);
	}
	
	\foreach \i in {5,6}{
		\draw[black!20] (\i*\i*\hxi*\hxi*4*4,0) -- (\i*\i*\hxi*\hxi*4*4,4);
	}
	
	\foreach \i in {1,...,2}{
		\draw[dashed,blue] (0, \i*4/3) -- (4, \i*4/3);
		\draw[dashed,red] (4, \i*4/3) -- (7.11, \i*4/3);
		\draw[dashed,black!20] (7.11, \i*4/3) -- (4*4, \i*4/3);
	}
	
	\foreach \j in {0,...,2}{
		\foreach \i in {1,...,\M}{
			\draw[blue] (0, \i*\heta*4/3 + \j*4/3) -- (4, \i*\heta*4/3 + \j*4/3);
			\draw[red] (4, \i*\heta*4/3 + \j*4/3) -- (7.11, \i*\heta*4/3 + \j*4/3);
			\draw[black!20] (7.11, \i*\heta*4/3 + \j*4/3) -- (4*4, \i*\heta*4/3 + \j*4/3);
		}
	}
	
	\node[left, yshift = 2] at (0,0) {$0$};
	\node[left] at (0,4) {$1$};

	\node[below, xshift = -2] at (0,0) {$0$};
	\node[below] at (4,0) {$1$};
	\node[below,red] at (7.11,0) {$h_{r,\widetilde{\widehat{\Omega}}_{C}} / h_{r, \widehat{\Omega}_C}$};
	\node[below,orange!40] at (15.5,0) {$1/h_{r, \widehat{\Omega}_C}$};
	
	\draw [-stealth, very thick] (17,2) to [out=30,in=150] (18.5,2) node[above] at (17.75,2.3) {$\mr{\bs\Phi}$};
\end{tikzpicture}

%% file: ref_elem_ext_ori.tikz
	\begin{tikzpicture}[baseline, scale = 0.6]
	
	\def \N {5}			
	\def \hxi{1/6}		
	
	\def \M {3} 		
	\def \heta{1/4}		

	\node[green!80] at (1.4,4.5) {$\widetilde{\widehat{\Omega}}_C$};
	\node[gray!80] at (3/8*3/8*4,4.5) {$\widehat{\Omega}_C$};
	\node at (2.5,4.5) {$\widehat{\Omega}$};
	
	\fill[gray!30] (0,0) rectangle (1,4);
	\fill[orange!10] (1,0) rectangle (4,4);
	\draw[pattern = north west lines, pattern color = green] (0,0) rectangle (4*4*1/36*4,4);
	
	\draw[-stealth] (4,0) -- (5, 0)  node[below, xshift = -2] {$r$};
	\draw[-stealth] (0,4) -- (0, 5) node[left, yshift = -2] {$\varphi$};

	\draw[blue] (0,0) -- (1, 0);
	\draw[blue] (0,4) -- (1, 4);
	\draw[red] (1,0) -- (4*4*1/36*4, 0);
	\draw[red] (1,4) -- (4*4*1/36*4, 4);
	\draw[black!20] (4*4*1/36*4,0) -- (4, 0);
	\draw[black!20] (4*4*1/36*4,4) -- (4, 4);
	\draw[red] (4*4*1/36*4,0) -- (4*4*1/36*4, 4);
	
	\foreach \i in {0,...,3}{
		\draw[blue] (\i*\i*\hxi*\hxi*4,0) -- (\i*\i*\hxi*\hxi*4,4);
	}
	
	\foreach \i in {5,6}{
		\draw[black!20] (\i*\i*\hxi*\hxi*4,0) -- (\i*\i*\hxi*\hxi*4,4);
	}
	
	\foreach \i in {1,...,2}{
		\draw[dashed,blue] (0, \i*4/3) -- (1, \i*4/3);
		\draw[dashed,red] (1, \i*4/3) -- (4*4*1/36*4, \i*4/3);
		\draw[dashed,black!20] (4*4*1/36*4, \i*4/3) -- (4, \i*4/3);
	}
	
	\foreach \j in {0,...,2}{
		\foreach \i in {1,...,\M}{
			\draw[blue] (0, \i*\heta*4/3 + \j*4/3) -- (1, \i*\heta*4/3 + \j*4/3);
			\draw[red]  (1, \i*\heta*4/3 + \j*4/3) -- (4*4*1/36*4, \i*\heta*4/3 + \j*4/3);
			\draw[black!20]  (4*4*1/36*4, \i*\heta*4/3 + \j*4/3) -- (4, \i*\heta*4/3 + \j*4/3);
		}
	}
	
	\node[left, yshift = 2] at (0,0) {$0$};
	\node[left] at (0,4) {$1$};
	
	\node[below, xshift = -2] at (0,0) {$0$};
	\node[below, yshift = 1, blue] at (\hxi*\hxi*3*3*4,0) {$h_{r,\widehat{\Omega}_C}$};
	\node[below, yshift = 1, xshift =10,red] at (4*4*1/36*4,0) {$h_{r,\widetilde{\widehat{\Omega}}_C}$};
	\node[below] at (4,0) {$1$};
	
	\draw [-stealth, very thick] (2,-1.5) to [out=270,in=60] (1,-4) node at (2.2,-2.75) {$\bs F$};
\end{tikzpicture}

%% file: ref_elem_ext_phys.tikz
	\begin{tikzpicture}[baseline=-1.3cm, scale = 2.5]	
	
	\def \N {5}			
	\def \hxi{1/6}		
	\def \M {3} 		
	\def \heta{1/4}		
	
	\filldraw[fill = gray!30, draw = white] (0:0.25) arc (0:300:0.25) -- (0:0);
	\filldraw[fill = orange!10, draw = white] (0:1) arc (0:300:1) -- (300:0.25) arc (300:0:0.25);
	\draw[pattern = north west lines, pattern color = green, draw = white] (0:4*4*\hxi*\hxi) arc (0:300:4*4*\hxi*\hxi) -- (0:0);
	
	\foreach \i in {0,...,3}{	
		\draw[blue] (\i*\i*\hxi*\hxi,0) arc (0:300:\i*\i*\hxi*\hxi);
	}
	
	\draw[red]  (4*4*\hxi*\hxi,0) arc (0:300:4*4*\hxi*\hxi);
	
	\foreach \i in {5,...,6}{	
		\draw[black!20]  (\i*\i*\hxi*\hxi,0) arc (0:300:\i*\i*\hxi*\hxi);
	}
	
	\foreach \j in {0,...,2}{
		\foreach \i in {1,...,\M}{
			\draw[blue] (\i*\heta*100 + \j*100:0) -- (\i*\heta*100 + \j*100:0.25);
			\draw[red]  (\i*\heta*100 + \j*100:0.25) -- (\i*\heta*100 + \j*100:4*4*\hxi*\hxi);
			\draw[black!20]  (\i*\heta*100 + \j*100:4*4*\hxi*\hxi) -- (\i*\heta*100 + \j*100:1);
		}
	}		
	
	\draw[black!20]  (0:1) -- (0:4*4*\hxi*\hxi);		
	\draw[black!20]  (300:1) -- (300:4*4*\hxi*\hxi);
	
	\draw[red] (0:0.25) -- (0:4*4*\hxi*\hxi);		
	\draw[red] (300:0.25) -- (300:4*4*\hxi*\hxi);
	
	\draw[blue] (0:0) -- (0:0.25);		
	\draw[blue] (0:0) -- (300:0.25);

	\draw[dash pattern=on 2pt off 1pt,black!20] (100:1) -- (100:0.25);		
	\draw[dash pattern=on 2pt off 1pt,black!20] (200:1) -- (200:0.25);
	
	\draw[dash pattern=on 2pt off 1pt,red] (100:0.25) -- (100:4*4*\hxi*\hxi);	
	\draw[dash pattern=on 2pt off 1pt,red] (200:0.25) -- (200:4*4*\hxi*\hxi);	
	
	\draw[dash pattern=on 2pt off 1pt,blue] (0:0) -- (100:0.25);	
	\draw[dash pattern=on 2pt off 1pt,blue] (0:0) -- (200:0.25);	
	
	\node at (270:1.12) {$\breve{\Omega}$};
	\node[gray!80] at (335:0.3) {$\breve{\Omega}_C$};
	\node[green] at (335:0.6) {$\widetilde{\breve\Omega}_C$};
	\draw [-stealth, very thick] (0:1.4) to [out=30,in=150] (0:2) node[above, yshift = 5] at (0:1.7) {$\breve {\bs \Phi}$};	
\end{tikzpicture}

%% file: ref_elem_ext_ori_phys.tikz
		\begin{tikzpicture}[baseline=-1.3cm, scale = 1.4]	
	
	\def \N {5}			
	\def \hxi{1/6}		
	\def \M {3} 		
	\def \heta{1/4}		
	
	\filldraw[fill = gray!30, draw = white] (0:0.25) arc (0:300:0.25) -- (0:0);
	\filldraw[fill = orange!10, draw = white] (0:1) arc (0:300:1) -- (300:0.25) arc (300:0:0.25);
	\draw[pattern = north west lines, pattern color = green, draw = white] (0:4*4*\hxi*\hxi) arc (0:300:4*4*\hxi*\hxi) -- (0:0);
	
	\foreach \i in {0,...,3}{	
		\draw[blue] (\i*\i*\hxi*\hxi,0) arc (0:300:\i*\i*\hxi*\hxi);
	}
	
	\draw[red]  (4*4*\hxi*\hxi,0) arc (0:300:4*4*\hxi*\hxi);
	
	\foreach \i in {5,...,6}{	
		\draw[black!20]  (\i*\i*\hxi*\hxi,0) arc (0:300:\i*\i*\hxi*\hxi);
	}
	
	\foreach \j in {0,...,2}{
		\foreach \i in {1,...,\M}{
			\draw[blue] (\i*\heta*100 + \j*100:0) -- (\i*\heta*100 + \j*100:0.25);
			\draw[red]  (\i*\heta*100 + \j*100:0.25) -- (\i*\heta*100 + \j*100:4*4*\hxi*\hxi);
			\draw[black!20]  (\i*\heta*100 + \j*100:4*4*\hxi*\hxi) -- (\i*\heta*100 + \j*100:1);
		}
	}		
	
	\draw[black!20]  (0:1) -- (0:4*4*\hxi*\hxi);		
	\draw[black!20]  (300:1) -- (300:4*4*\hxi*\hxi);
	
	\draw[red] (0:0.25) -- (0:4*4*\hxi*\hxi);		
	\draw[red] (300:0.25) -- (300:4*4*\hxi*\hxi);
	
	\draw[blue] (0:0) -- (0:0.25);		
	\draw[blue] (0:0) -- (300:0.25);

	\draw[dash pattern=on 2pt off 1pt,black!20] (100:1) -- (100:0.25);		
	\draw[dash pattern=on 2pt off 1pt,black!20] (200:1) -- (200:0.25);
	
	\draw[dash pattern=on 2pt off 1pt,red] (100:0.25) -- (100:4*4*\hxi*\hxi);	
	\draw[dash pattern=on 2pt off 1pt,red] (200:0.25) -- (200:4*4*\hxi*\hxi);	
	
	\draw[dash pattern=on 2pt off 1pt,blue] (0:0) -- (100:0.25);	
	\draw[dash pattern=on 2pt off 1pt,blue] (0:0) -- (200:0.25);	
	
	\node at (270:1.15) {$\Omega$};	
	\node[gray!80] at (335:0.4) {$\Omega_C$};
	\node[green] at (335:0.85) {$\widetilde{\Omega}_C$};	
\end{tikzpicture}

%% file: laplace_unitdisk300_mg_neumann_hom_H1.tikz
%
%
\definecolor{mycolor1}{rgb}{0.00000,0.44700,0.74100}%
\definecolor{mycolor2}{rgb}{0.85000,0.32500,0.09800}%
\definecolor{mycolor3}{rgb}{0.92900,0.69400,0.12500}%
\definecolor{mycolor4}{rgb}{0.49400,0.18400,0.55600}%
\definecolor{mycolor5}{rgb}{0.46600,0.67400,0.18800}%
\definecolor{mycolor6}{rgb}{0.30100,0.74500,0.93300}%
\definecolor{mycolor7}{rgb}{0.63500,0.07800,0.18400}%
\begin{tikzpicture}
\tikzstyle{every node}=[font=\large]
\begin{axis}[%
width=4.65in,
height=3.861in,
at={(0.78in,0.521in)},
scale only axis,
xmode=log,
xmin=0.00390625,
xmax=0.25,
xminorticks=true,
xlabel style={font=\color{white!15!black}},
xlabel={\Large $h$},
ymode=log,
ymin=8.94323494649383e-13,
ymax=0.0282516659997199,
yminorticks=true,
ylabel style={font=\color{white!15!black}},
ylabel={\Large $||u-u_h||_{L^2(\Omega)}$},
axis background/.style={fill=white},
legend style={at={(0.97,0.03)}, anchor=south east, legend cell align=left, align=left, draw=white!15!black}
]
\addplot [color=mycolor1, line width=2.0pt, mark size=3.5pt, mark=square, mark options={solid, mycolor1}]
  table[row sep=crcr]{%
0.25	0.0212135094676677\\
0.125	0.00532242648887976\\
0.0625	0.00133804138890919\\
0.03125	0.000336390886455377\\
0.015625	8.45344921481829e-05\\
0.0078125	2.12309282461713e-05\\
0.00390625	5.32903570298247e-06\\
};
\addlegendentry{$p=1, k=0, \mu=0.54$}

\addplot [color=mycolor2, dotted, line width=2.0pt, mark size=3.3pt, mark=triangle, mark options={solid, mycolor2}]
  table[row sep=crcr]{%
0.25	0.0282516659997199\\
0.125	0.0100146687441624\\
0.0625	0.00380397601794746\\
0.03125	0.00151298154463825\\
0.015625	0.000621095372889056\\
0.0078125	0.000260459173026761\\
0.00390625	0.000110762299050118\\
};
\addlegendentry{$p=1, k=0, \mu=1$}

\addplot [color=mycolor3, line width=2.0pt, mark size=3.5pt, mark=square, mark options={solid, mycolor3}]
  table[row sep=crcr]{%
0.25	0.00452017932188995\\
0.125	0.000609965334405823\\
0.0625	7.59601717116672e-05\\
0.03125	9.34304040331861e-06\\
0.015625	1.16323094824584e-06\\
0.0078125	1.45255655347731e-07\\
0.00390625	1.81522293764666e-08\\
};
\addlegendentry{$p=2, k=1, \mu=0.27$}

\addplot [color=mycolor4, dotted, line width=2.0pt, mark size=3.3pt, mark=triangle, mark options={solid, mycolor4}]
  table[row sep=crcr]{%
0.25	0.00741826526980411\\
0.125	0.00253390462373358\\
0.0625	0.000997496654340443\\
0.03125	0.000412300852237635\\
0.015625	0.000174223633049253\\
0.0078125	7.45407288073167e-05\\
0.00390625	3.21243777456237e-05\\
};
\addlegendentry{$p=2, k=1, \mu=1$}

\addplot [color=mycolor5, line width=2.0pt, mark size=3.5pt, mark=square, mark options={solid, mycolor5}]
  table[row sep=crcr]{%
0.25	0.00158049152007544\\
0.125	0.000173713293153119\\
0.0625	1.49725468656479e-05\\
0.03125	1.06335750432993e-06\\
0.015625	7.02604289553734e-08\\
0.0078125	4.50000369205461e-09\\
0.00390625	2.84441966264726e-10\\
};
\addlegendentry{$p=3, k=2, \mu=0.18$}

\addplot [color=mycolor6, dotted, line width=2.0pt, mark size=3.3pt, mark=triangle, mark options={solid, mycolor6}]
  table[row sep=crcr]{%
0.25	0.00336582005997648\\
0.125	0.00130605636109783\\
0.0625	0.000526561223686412\\
0.03125	0.000218546795844632\\
0.015625	9.24253475051244e-05\\
0.0078125	3.95491384060394e-05\\
0.00390625	1.70438292072048e-05\\
};
\addlegendentry{$p=3, k=2, \mu=1$}

\addplot [color=mycolor7, line width=2.0pt, mark size=3.5pt, mark=square, mark options={solid, mycolor7}]
  table[row sep=crcr]{%
0.25	0.00113057464211417\\
0.125	0.000112483216406817\\
0.0625	8.52557185494585e-06\\
0.03125	3.58962568717999e-07\\
0.015625	1.22752860656131e-08\\
0.0078125	3.94597415440279e-10\\
0.00390625	1.24502344070233e-11\\
};
\addlegendentry{$p=4, k=3, \mu=0.135$}

\addplot [color=mycolor1, dotted, line width=2.0pt, mark size=3.3pt, mark=triangle, mark options={solid, mycolor1}]
  table[row sep=crcr]{%
0.25	0.00212183518931365\\
0.125	0.00083506993477436\\
0.0625	0.00033677114729115\\
0.03125	0.000139650781912087\\
0.015625	5.90109544014336e-05\\
0.0078125	2.5235929984045e-05\\
0.00390625	1.08712271112973e-05\\
};
\addlegendentry{$p=4, k=3, \mu=1$}

\addplot [color=mycolor2, line width=2.0pt, mark size=3.5pt, mark=square, mark options={solid, mycolor2}]
  table[row sep=crcr]{%
0.25	0.000818891620131171\\
0.125	9.08478158512378e-05\\
0.0625	7.353336158629e-06\\
0.03125	1.90757668456449e-07\\
0.015625	3.44847095477327e-09\\
0.0078125	5.63649566673493e-11\\
0.00390625	8.94323494649383e-13\\
};
\addlegendentry{$p=5, k=4, \mu=0.108$}

\addplot [color=mycolor3, dotted, line width=2.0pt, mark size=3.3pt, mark=triangle, mark options={solid, mycolor3}]
  table[row sep=crcr]{%
0.25	0.00147030867757567\\
0.125	0.000590119632712504\\
0.0625	0.000237973960385099\\
0.03125	9.85887846272814e-05\\
0.015625	4.16285938666703e-05\\
0.0078125	1.77932456288748e-05\\
0.00390625	7.66252167593269e-06\\
};
\addlegendentry{$p=5, k=4, \mu=1$}

\logLogSlopeTriangleNeg{0.08}{0.20}{0.855}{-1.2}{gray, line width = 2pt}{1.2};
\logLogSlopeTriangleUpsideDownNeg{0.21}{0.15}{0.635}{-2}{gray, line width = 2pt}{2};
\logLogSlopeTriangleNeg{0.045}{0.15}{0.54}{-3}{gray, line width = 2pt}{3};
\logLogSlopeTriangleNeg{0.045}{0.12}{0.385}{-4}{gray, line width = 2pt}{4};
\logLogSlopeTriangleNeg{0.045}{0.1}{0.25}{-5}{gray, line width = 2pt}{5};
\logLogSlopeTriangleUpsideDownNeg{0.22}{0.12}{0.065}{-6}{gray, line width = 2pt}{6};

\end{axis}
\end{tikzpicture}%

%% file: laplace_unitdisk300_mg_neumann_hom_L2.tikz
%
%
\definecolor{mycolor1}{rgb}{0.00000,0.44700,0.74100}%
\definecolor{mycolor2}{rgb}{0.85000,0.32500,0.09800}%
\definecolor{mycolor3}{rgb}{0.92900,0.69400,0.12500}%
\definecolor{mycolor4}{rgb}{0.49400,0.18400,0.55600}%
\definecolor{mycolor5}{rgb}{0.46600,0.67400,0.18800}%
\definecolor{mycolor6}{rgb}{0.30100,0.74500,0.93300}%
\definecolor{mycolor7}{rgb}{0.63500,0.07800,0.18400}%
\begin{tikzpicture}
\tikzstyle{every node}=[font=\large]
\begin{axis}[%
width=4.568in,
height=3.565in,
at={(0.766in,0.524in)},
scale only axis,
xmode=log,
xmin=0.00390625,
xmax=0.25,
xminorticks=true,
xlabel style={font=\color{white!15!black}},
xlabel={\Large $h$},
ymode=log,
ymin=5.10901078386652e-10,
ymax=0.281519920871668,
yminorticks=true,
ylabel style={font=\Large\color{white!15!black}},
ylabel={\Large $||u-u_h||_{H^1(\Omega)}$},
axis background/.style={fill=white},
legend style={at={(0.98,0.02)}, anchor=south east, legend cell align=left, align=left, draw=white!15!black}
]
\addplot [color=mycolor1, line width=2.0pt, mark size=3.5pt, mark=square, mark options={solid, mycolor1}]
  table[row sep=crcr]{%
0.25	0.239015993386356\\
0.125	0.119092799597025\\
0.0625	0.0597812496615245\\
0.03125	0.0300471153736351\\
0.015625	0.0150985558705997\\
0.0078125	0.00758245816711956\\
0.00390625	0.00380556016993941\\
};
\addlegendentry{$p=1, k=0, \mu=0.54$}

\addplot [color=mycolor2, dotted, line width=2.0pt, mark size=3.3pt, mark=triangle, mark options={solid, mycolor2}]
  table[row sep=crcr]{%
0.25	0.281519920871668\\
0.125	0.164788855363779\\
0.0625	0.10038741515964\\
0.03125	0.0628895183957252\\
0.015625	0.0401653798011313\\
0.0078125	0.025979157708927\\
0.00390625	0.0169378286560218\\
};
\addlegendentry{$p=1, k=0, \mu=1$}

\addplot [color=mycolor3, line width=2.0pt, mark size=3.5pt, mark=square, mark options={solid, mycolor3}]
  table[row sep=crcr]{%
0.25	0.0553971462026986\\
0.125	0.0146611878558397\\
0.0625	0.00380791683535774\\
0.03125	0.000977105231324648\\
0.015625	0.000249223255940375\\
0.0078125	6.32274528905883e-05\\
0.00390625	1.59758728966215e-05\\
};
\addlegendentry{$p=2, k=1, \mu=0.27$}

\addplot [color=mycolor4, dotted, line width=2.0pt, mark size=3.3pt, mark=triangle, mark options={solid, mycolor4}]
  table[row sep=crcr]{%
0.25	0.121033298441414\\
0.125	0.074961644650429\\
0.0625	0.048527619586381\\
0.03125	0.0317803472497654\\
0.015625	0.020896594457155\\
0.0078125	0.0137640943872862\\
0.00390625	0.00907358757045541\\
};
\addlegendentry{$p=2, k=1, \mu=1$}

\addplot [color=mycolor5, line width=2.0pt, mark size=3.5pt, mark=square, mark options={solid, mycolor5}]
  table[row sep=crcr]{%
0.25	0.0252292115493342\\
0.125	0.00494738497386289\\
0.0625	0.000767837546647352\\
0.03125	0.000106809845657862\\
0.015625	1.41580447018831e-05\\
0.0078125	1.83021925684159e-06\\
0.00390625	2.33363828644649e-07\\
};
\addlegendentry{$p=3, k=2, \mu=0.18$}

\addplot [color=mycolor6, dotted, line width=2.0pt, mark size=3.3pt, mark=triangle, mark options={solid, mycolor6}]
  table[row sep=crcr]{%
0.25	0.0818883150674496\\
0.125	0.0532540170326016\\
0.0625	0.0348821461780438\\
0.03125	0.0229328144804416\\
0.015625	0.0151036280099091\\
0.0078125	0.00995603122644309\\
0.00390625	0.00656569312704284\\
};
\addlegendentry{$p=3, k=2, \mu=1$}

\addplot [color=mycolor7, line width=2.0pt, mark size=3.5pt, mark=square, mark options={solid, mycolor7}]
  table[row sep=crcr]{%
0.25	0.0214619222576052\\
0.125	0.00367971144820587\\
0.0625	0.0003762138822617\\
0.03125	2.897408188854e-05\\
0.015625	1.99515019259084e-06\\
0.0078125	1.30822823223141e-07\\
0.00390625	8.38029022117391e-09\\
};
\addlegendentry{$p=4, k=3, \mu=0.135$}

\addplot [color=mycolor1, dotted, line width=2.0pt, mark size=3.3pt, mark=triangle, mark options={solid, mycolor1}]
  table[row sep=crcr]{%
0.25	0.0637949340171673\\
0.125	0.0417888094148538\\
0.0625	0.0274368286190352\\
0.03125	0.0180579114210243\\
0.015625	0.0118994977617107\\
0.0078125	0.00784605788564863\\
0.00390625	0.00517493020816426\\
};
\addlegendentry{$p=4, k=3, \mu=1$}

\addplot [color=mycolor2, line width=2.0pt, mark size=3.5pt, mark=square, mark options={solid, mycolor2}]
  table[row sep=crcr]{%
0.25	0.0200469769296328\\
0.125	0.00367229436219096\\
0.0625	0.000288306489849215\\
0.03125	1.29797242122536e-05\\
0.015625	4.73246153312748e-07\\
0.0078125	1.58418250612684e-08\\
0.00390625	5.10901078386652e-10\\
};
\addlegendentry{$p=5, k=4, \mu=0.108$}

\addplot [color=mycolor3, dotted, line width=2.0pt, mark size=3.3pt, mark=triangle, mark options={solid, mycolor3}]
  table[row sep=crcr]{%
0.25	0.0519129114614354\\
0.125	0.0342594875998191\\
0.0625	0.0225218345280173\\
0.03125	0.0148314684124985\\
0.015625	0.00977614790110057\\
0.0078125	0.00644690881081315\\
0.00390625	0.00425240678764475\\
};
\addlegendentry{$p=5, k=4, \mu=1$}

\logLogSlopeTriangleNeg{0.08}{0.2}{0.925}{-0.6}{gray, line width = 2pt}{0.6};
\logLogSlopeTriangleUpsideDownNeg{0.2}{0.15}{0.76}{-1}{gray, line width = 2pt}{1};
\logLogSlopeTriangleNeg{0.045}{0.15}{0.63}{-2}{gray, line width = 2pt}{2};
\logLogSlopeTriangleNeg{0.045}{0.15}{0.46}{-3}{gray, line width = 2pt}{3};
\logLogSlopeTriangleNeg{0.045}{0.12}{0.31}{-4}{gray, line width = 2pt}{4};
\logLogSlopeTriangleUpsideDownNeg{0.22}{0.12}{0.06}{-5}{gray, line width = 2pt}{5};

\end{axis}
\end{tikzpicture}%

%% file: laplace_Lshaped_mixed_hom_p=1to5_n=6_H1_err_h_slopes.tikz
%
%
\definecolor{mycolor1}{rgb}{0.00000,0.44700,0.74100}%
\definecolor{mycolor2}{rgb}{0.85000,0.32500,0.09800}%
\definecolor{mycolor3}{rgb}{0.92900,0.69400,0.12500}%
\definecolor{mycolor4}{rgb}{0.49400,0.18400,0.55600}%
\definecolor{mycolor5}{rgb}{0.46600,0.67400,0.18800}%
\definecolor{mycolor6}{rgb}{0.30100,0.74500,0.93300}%
\definecolor{mycolor7}{rgb}{0.63500,0.07800,0.18400}%
\begin{tikzpicture}
\tikzstyle{every node}=[font=\large]
\begin{axis}[%
width=4.568in,
height=3.603in,
at={(0.766in,0.486in)},
scale only axis,
xmode=log,
xmin=0.00390625,
xmax=0.25,
xminorticks=true,
xlabel style={font=\color{white!15!black}},
xlabel={\Large $h$},
ymode=log,
ymin=1.8877939912865e-08,
ymax=0.682718688410901,
yminorticks=true,
ylabel style={font=\color{white!15!black}},
ylabel={\Large $||u-u_h||_{H^1(\Omega)}$},
axis background/.style={fill=white},
legend style={at={(0.97,0.03)}, anchor=south east, legend cell align=left, align=left, draw=white!15!black}
]
\addplot [color=mycolor1, line width=2.0pt, mark size=3.5pt, mark=square, mark options={solid, mycolor1}]
  table[row sep=crcr]{%
0.25	0.672026371640505\\
0.125	0.381348366425487\\
0.0625	0.200859256834002\\
0.03125	0.102728699251798\\
0.015625	0.0518682924109016\\
0.0078125	0.0260764434395335\\
0.00390625	0.0130893068838009\\
};
\addlegendentry{$p=1, k=0, \mu=0.3$}

\addplot [color=mycolor2, dotted, line width=2.0pt, mark size=3.3pt, mark=triangle, mark options={solid, mycolor2}]
  table[row sep=crcr]{%
0.25	0.682718688410901\\
0.125	0.439163316936274\\
0.0625	0.303766036580537\\
0.03125	0.224435083627358\\
0.015625	0.172218380713083\\
0.0078125	0.13459077710127\\
0.00390625	0.106058591020759\\
};
\addlegendentry{$p=1, k=0, \mu=1$}

\addplot [color=mycolor3, line width=2.0pt, mark size=3.5pt, mark=square, mark options={solid, mycolor3}]
  table[row sep=crcr]{%
0.25	0.340945066585529\\
0.125	0.169481375406107\\
0.0625	0.0582040086918835\\
0.03125	0.0156033685114532\\
0.015625	0.00391594840688262\\
0.0078125	0.000977624320222386\\
0.00390625	0.000244454962792321\\
};
\addlegendentry{$p=2, k=1, \mu=0.15$}

\addplot [color=mycolor4, dotted, line width=2.0pt, mark size=3.3pt, mark=triangle, mark options={solid, mycolor4}]
  table[row sep=crcr]{%
0.25	0.331700791483333\\
0.125	0.244513853675242\\
0.0625	0.189451514704371\\
0.03125	0.149404065692206\\
0.015625	0.118157476856647\\
0.0078125	0.0935717211255763\\
0.00390625	0.0741630499368871\\
};
\addlegendentry{$p=2, k=1, \mu=1$}

\addplot [color=mycolor5, line width=2.0pt, mark size=3.5pt, mark=square, mark options={solid, mycolor5}]
  table[row sep=crcr]{%
0.25	0.112397075965167\\
0.125	0.0476700204117239\\
0.0625	0.0124417427918611\\
0.03125	0.00193786641051676\\
0.015625	0.000267200553749455\\
0.0078125	3.62718117767082e-05\\
0.00390625	4.77793486755253e-06\\
};
\addlegendentry{$p=3, k=2, \mu=0.1$}

\addplot [color=mycolor6, dotted, line width=2.0pt, mark size=3.3pt, mark=triangle, mark options={solid, mycolor6}]
  table[row sep=crcr]{%
0.25	0.252407331842546\\
0.125	0.195899596456269\\
0.0625	0.15431646907969\\
0.03125	0.121965752185187\\
0.015625	0.0965519668593673\\
0.0078125	0.076507454347093\\
0.00390625	0.0606611608322241\\
};
\addlegendentry{$p=3, k=2, \mu=1$}

\addplot [color=mycolor7, line width=2.0pt, mark size=3.5pt, mark=square, mark options={solid, mycolor7}]
  table[row sep=crcr]{%
0.25	0.11457346198794\\
0.125	0.0337277807523036\\
0.0625	0.005119832738104\\
0.03125	0.000477738940738975\\
0.015625	3.63342826608659e-05\\
0.0078125	2.50398473982696e-06\\
0.00390625	1.64358879813407e-07\\
};
\addlegendentry{$p=4, k=3, \mu=0.075$}

\addplot [color=mycolor1, dotted, line width=2.0pt, mark size=3.3pt, mark=triangle, mark options={solid, mycolor1}]
  table[row sep=crcr]{%
0.25	0.212537241590212\\
0.125	0.167600665611599\\
0.0625	0.132243431045029\\
0.03125	0.104604390965877\\
0.015625	0.0828487312794479\\
0.0078125	0.0656693806213353\\
0.00390625	0.0520780460192825\\
};
\addlegendentry{$p=4, k=3, \mu=1$}

\addplot [color=mycolor2, line width=2.0pt, mark size=3.5pt, mark=square, mark options={solid, mycolor2}]
  table[row sep=crcr]{%
0.25	0.108583657160326\\
0.125	0.0372403203580947\\
0.0625	0.00563352788888548\\
0.03125	0.00036596284415336\\
0.015625	1.57222372653923e-05\\
0.0078125	5.66280001013689e-07\\
0.00390625	1.8877939912865e-08\\
};
\addlegendentry{$p=5, k=4, \mu=0.06$}

\addplot [color=mycolor3, dotted, line width=2.0pt, mark size=3.3pt, mark=triangle, mark options={solid, mycolor3}]
  table[row sep=crcr]{%
0.25	0.185681094207607\\
0.125	0.147111080731\\
0.0625	0.116235597817571\\
0.03125	0.0919883701698068\\
0.015625	0.072878121419394\\
0.0078125	0.0577770095672102\\
0.00390625	0.045824497988099\\
};
\addlegendentry{$p=5, k=4, \mu=1$}

\logLogSlopeTriangleNeg{0.09}{0.2}{0.945}{-0.33}{gray, line width = 2pt}{0.33};
\logLogSlopeTriangleUpsideDownNeg{0.22}{0.15}{0.76}{-1}{gray, line width = 2pt}{1};
\logLogSlopeTriangleNeg{0.045}{0.15}{0.67}{-2}{gray, line width = 2pt}{2};
\logLogSlopeTriangleNeg{0.045}{0.15}{0.49}{-3}{gray, line width = 2pt}{3};
\logLogSlopeTriangleNeg{0.045}{0.125}{0.32}{-4}{gray, line width = 2pt}{4};
\logLogSlopeTriangleUpsideDownNeg{0.22}{0.125}{0.065}{-5}{gray, line width = 2pt}{5};

\end{axis}
\end{tikzpicture}%

%% file: laplace_Lshaped_mixed_hom_p=1to5_n=6_L2_err_h_slopes.tikz
%
%
\definecolor{mycolor1}{rgb}{0.00000,0.44700,0.74100}%
\definecolor{mycolor2}{rgb}{0.85000,0.32500,0.09800}%
\definecolor{mycolor3}{rgb}{0.92900,0.69400,0.12500}%
\definecolor{mycolor4}{rgb}{0.49400,0.18400,0.55600}%
\definecolor{mycolor5}{rgb}{0.46600,0.67400,0.18800}%
\definecolor{mycolor6}{rgb}{0.30100,0.74500,0.93300}%
\definecolor{mycolor7}{rgb}{0.63500,0.07800,0.18400}%
\begin{tikzpicture}
\tikzstyle{every node}=[font=\large]
\begin{axis}[%
width=4.65in,
height=3.861in,
at={(0.78in,0.521in)},
scale only axis,
xmode=log,
xmin=0.00390625,
xmax=0.25,
xminorticks=true,
xlabel style={font=\color{white!15!black}},
xlabel={\Large $h$},
ymode=log,
ymin=3.56958327094639e-11,
ymax=0.142346881340735,
yminorticks=true,
ylabel style={font=\color{white!15!black}},
ylabel={\Large $||u-u_h||_{L^2(\Omega)}$},
axis background/.style={fill=white},
legend style={at={(0.97,0.03)}, anchor=south east, legend cell align=left, align=left, draw=white!15!black}
]
\addplot [color=mycolor1, line width=2.0pt, mark size=3.5pt, mark=square, mark options={solid, mycolor1}]
  table[row sep=crcr]{%
0.25	0.121148767527469\\
0.125	0.0440642232361699\\
0.0625	0.0127357285826844\\
0.03125	0.00337991559759156\\
0.015625	0.000861926633320303\\
0.0078125	0.000217284012973081\\
0.00390625	5.45960445778544e-05\\
};
\addlegendentry{$p=1, k=0, \mu=0.3$}

\addplot [color=mycolor2, dotted, line width=2.0pt, mark size=3.3pt, mark=triangle, mark options={solid, mycolor2}]
  table[row sep=crcr]{%
0.25	0.142346881340735\\
0.125	0.0654474846846619\\
0.0625	0.032682579181296\\
0.03125	0.0184854860685843\\
0.015625	0.0110930899111152\\
0.0078125	0.00683373227617706\\
0.00390625	0.00425937554980156\\
};
\addlegendentry{$p=1, k=0, \mu=1$}

\addplot [color=mycolor3, line width=2.0pt, mark size=3.5pt, mark=square, mark options={solid, mycolor3}]
  table[row sep=crcr]{%
0.25	0.046705601153808\\
0.125	0.0170028884093574\\
0.0625	0.00378207727134886\\
0.03125	0.00050323428291148\\
0.015625	5.86372417920254e-05\\
0.0078125	7.03972991611173e-06\\
0.00390625	8.68367716081826e-07\\
};
\addlegendentry{$p=2, k=1, \mu=0.15$}

\addplot [color=mycolor4, dotted, line width=2.0pt, mark size=3.3pt, mark=triangle, mark options={solid, mycolor4}]
  table[row sep=crcr]{%
0.25	0.0451054773529105\\
0.125	0.0236744462511736\\
0.0625	0.014309537492909\\
0.03125	0.00888356517841622\\
0.015625	0.00555481115297951\\
0.0078125	0.00348469158463207\\
0.00390625	0.00218984721421057\\
};
\addlegendentry{$p=2, k=1, \mu=1$}

\addplot [color=mycolor5, line width=2.0pt, mark size=3.5pt, mark=square, mark options={solid, mycolor5}]
  table[row sep=crcr]{%
0.25	0.00746927893435205\\
0.125	0.00313537973864161\\
0.0625	0.00078059911397498\\
0.03125	6.64961734117262e-05\\
0.015625	4.61264087282087e-06\\
0.0078125	3.39596980553819e-07\\
0.00390625	2.3337045427044e-08\\
};
\addlegendentry{$p=3, k=2, \mu=0.1$}

\addplot [color=mycolor6, dotted, line width=2.0pt, mark size=3.3pt, mark=triangle, mark options={solid, mycolor6}]
  table[row sep=crcr]{%
0.25	0.0267646939759221\\
0.125	0.0161474411576558\\
0.0625	0.00995812001400725\\
0.03125	0.00621382914693119\\
0.015625	0.00389453359857744\\
0.0078125	0.00244628272748764\\
0.00390625	0.00153843102076257\\
};
\addlegendentry{$p=3, k=2, \mu=1$}

\addplot [color=mycolor7, line width=2.0pt, mark size=3.5pt, mark=square, mark options={solid, mycolor7}]
  table[row sep=crcr]{%
0.25	0.00718637734264031\\
0.125	0.00110546679666313\\
0.0625	0.000110292213093007\\
0.03125	6.24557060767255e-06\\
0.015625	2.78063568248993e-07\\
0.0078125	1.05108199801407e-08\\
0.00390625	3.52519815901708e-10\\
};
\addlegendentry{$p=4, k=3, \mu=0.075$}

\addplot [color=mycolor1, dotted, line width=2.0pt, mark size=3.3pt, mark=triangle, mark options={solid, mycolor1}]
  table[row sep=crcr]{%
0.25	0.0202330593274469\\
0.125	0.0124242283665208\\
0.0625	0.00771323868144214\\
0.03125	0.00482337648917527\\
0.015625	0.00302643729498409\\
0.0078125	0.00190221478406505\\
0.00390625	0.00119672449767131\\
};
\addlegendentry{$p=4, k=3, \mu=1$}

\addplot [color=mycolor2, line width=2.0pt, mark size=3.5pt, mark=square, mark options={solid, mycolor2}]
  table[row sep=crcr]{%
0.25	0.00493631295482084\\
0.125	0.00104675297950858\\
0.0625	0.000116888764019734\\
0.03125	5.76532546986815e-06\\
0.015625	1.26460468369425e-07\\
0.0078125	2.20380871545781e-09\\
0.00390625	3.56958327094639e-11\\
};
\addlegendentry{$p=5, k=4, \mu=0.06$}

\addplot [color=mycolor3, dotted, line width=2.0pt, mark size=3.3pt, mark=triangle, mark options={solid, mycolor3}]
  table[row sep=crcr]{%
0.25	0.0163498394391707\\
0.125	0.0101605678307338\\
0.0625	0.00632722319874956\\
0.03125	0.00396175002613188\\
0.015625	0.00248751909435678\\
0.0078125	0.00156410398064433\\
0.00390625	0.000984241123455484\\
};
\addlegendentry{$p=5, k=4, \mu=1$}

\logLogSlopeTriangleNeg{0.09}{0.2}{0.905}{-0.66}{gray, line width = 2pt}{0.66};
\logLogSlopeTriangleNeg{0.04}{0.15}{0.74}{-2}{gray, line width = 2pt}{2};
\logLogSlopeTriangleNeg{0.045}{0.15}{0.595}{-3}{gray, line width = 2pt}{3};
\logLogSlopeTriangleNeg{0.045}{0.15}{0.465}{-4}{gray, line width = 2pt}{4};
\logLogSlopeTriangleNeg{0.045}{0.125}{0.295}{-5}{gray, line width = 2pt}{5};
\logLogSlopeTriangleUpsideDownNeg{0.22}{0.125}{0.065}{-6}{gray, line width = 2pt}{6};

\end{axis}
\end{tikzpicture}%